\documentclass[10pt,reqno]{amsart}
\usepackage{graphicx}
\usepackage{epstopdf}
\usepackage{pdfpages}
\usepackage{epsfig,amsthm,amsmath}
\usepackage{color,cases,amssymb}
\usepackage{tabularx}
\usepackage{hyperref}

\newtheorem{theorem}{Theorem}[section]

\theoremstyle{definition}

\theoremstyle{remark}
\newtheorem{remark}[theorem]{Remark}

\newcommand{\fe}{\mathrm{e}}

\newcommand{\eps}{\varepsilon}
\newcommand{\bR}{{\mathbb R}}
\newcommand{\bC}{{\mathbb C}}

\newcommand{\bN}{{\mathbb N}}

\newcommand{\re}{{\mathrm{Re}}}

\newcommand{\bx}{\mathbf{x}}

\numberwithin{equation}{section}

\begin{document}

\title[PML for Klein-Gordon equation]{Pseudospectral methods with PML for nonlinear Klein-Gordon equations in classical and non-relativistic regimes}


\author[X. Antoine]{Xavier Antoine}
\address{\hspace*{-12pt}X.~Antoine:  Universit\'e de Lorraine, CNRS, Inria, IECL, F-54000 Nancy, France.}
\email{xavier.antoine@univ-lorraine.fr}

\author[X. Zhao]{Xiaofei Zhao}
\address{\hspace*{-12pt}X.~Zhao: School of Mathematics and Statistics \& Computational Sciences Hubei Key Laboratory, Wuhan University, 430072 Wuhan, China}
\email{matzhxf@whu.edu.cn}


\date{}

\dedicatory{}

\begin{abstract}
Two different Perfectly Matched Layer (PML) formulations with efficient pseudo-spectral numerical schemes are derived for the standard and non-relativistic nonlinear Klein-Gordon equations (NKGE). A pseudo-spectral explicit exponential integrator scheme for a first-order formulation and a linearly implicit preconditioned finite-difference scheme for a second-order  formulation are proposed and analyzed. To obtain a high spatial accuracy, new regularized Berm\'udez type absorption profiles are introduced for the PML. It is shown that the two schemes are efficient, but the linearly implicit scheme should be preferred for accuracy purpose
when used within the framework of pseudo-spectral methods combined with the regularized Berm\'udez type functions. In addition, in the non-relativistic regime, numerical examples lead to the conclusion that the error related to regularized Berm\'udez type absorption functions is insensitive to the small parameter $\varepsilon$ involved in the NKGE. The paper ends by a two-dimensional example  showing that the strategy
 extends to the rotating NKGE 
where the vortex dynamics  is very well-reproduced.


\end{abstract}

\maketitle


\section{Introduction}
The Klein-Gordon equation was originally proposed to describe the dynamics of spinless particles, and is known as the relativistic version of the Schr\"odinger equation. It  has been widely applied in the studies of quantum field theory, cosmology and plasma physics \cite{BaoZhao-KGZ,book,book2,RKG-PRD}.  
In this work, we consider
the following $d$-dimensional ($d=1,2,3$) nonlinear Klein-Gordon  equation (NKGE) \cite{BaoZhao,book,Nonrelat2}:
\begin{equation}\label{KG model}
\left\{\begin{split}
&\frac{1}{c^{2}}\partial_{tt} u(\bx,t)-\Delta u(\bx,t)+\frac{m^2c^2}{\hbar^2}u(\bx,t)+\lambda |u(\bx,t)|^2u(\bx,t)=0,\quad t>0,\ \bx\in\bR^d,\\
&u(\bx,0)=u_0(\bx),\quad \partial_tu(\bx,0)=v_0(\bx),\quad \bx\in\bR^d,
\end{split}\right.
\end{equation}
where $c>0$ is the speed of light, $\hbar>0$ denotes the Planck constant, $m>0$ is the particle mass and $\lambda>0$ is a given constant describing the strength of the defocusing interaction. In addition,   the unknown function is $u=u(\bx,t):\bR^d\times[0,\infty)\to\bC$, where $u_0,v_0:\bR^d\to\bC$ are two given initial data. Here, the global well-posedness of the defocusing NKGE (\ref{KG model}) is guaranteed by \cite{wellpose}, while for the focusing case $\lambda<0$, there could be finite time blow-up \cite{blow-up}. The following energy or Hamiltonian of (\ref{KG model}) is also conserved
\begin{equation}\label{energy def}
H(t):=\int_{\bR^d}\left[\frac{1}{c^2}|\partial_tu(\bx,t)|^2+|\nabla u(\bx,t)|^2+\frac{m^2c^2}{\hbar^2}|u(\bx,t)|^2+
\frac{\lambda}{2}|u(\bx,t)|^4\right]d\bx\equiv H(0),\quad t\geq0.
\end{equation}

Since the initial-value problem  (\ref{KG model}) is set in $\bR^d$, the spatial domain has to be truncated to use a standard numerical discretization scheme, e.g. finite-difference, finite-element or pseudo-spectral scheme \cite{BaoZhao}.
In the present paper, our aim is to apply the Fourier pseudo-spectral discretization scheme \cite{Shen} which leads to highly
accurate numerical solutions of PDEs with smooth coefficients
and is widely used in computations of quantum mechanics \cite{ABB,Antoine-lib1,Antoine-lib2,Antoine-lib3,AD2015,BC,BCZ,Dong,BDZ,BaoZhao,BaoZhao-KGZ,EWI,RKG}. As a consequence,  one needs to impose periodic boundary conditions on the fictitious boundary delimiting the (rectangular) computational domain.  Because of this constraint, any truncation technique based
on non-reflecting/artificial/absorbing boundary condition for nonlinear PDEs
\cite{ReviewCICP,Antoine-review} cannot be applied. Since its introduction in the seminal paper by B\'erenger \cite{Classical-sigma1},
the Perfectly Matched Layer (PML) approach  provides an alternative powerful tool to simulate the numerical solution of  PDEs in unbounded domains \cite{PML-siap,PML-wave,ReviewCICP,Antoine-review,Antoine-RNLS,Wu,Bermudez,BermudezJCP,BermudezSISC,Classical-sigma2,Chen,PML-notes,PML-NLS}.
Concerning nonlinear PDEs, the PML technique has been 
 applied for example to the nonlinear Schr\"odinger equations \cite{ReviewCICP,Antoine-review,Antoine-RNLS,PML-NLS}, Euler and Navier-Stokes equations \cite{Hu} or two-fluid plasma equations \cite{two-fluid}. Concerning the nonlinear wave equations,
Appel\"o and Kreiss \cite{PML-wave} proposed and studied a first-order
PML formulation with Dirichlet/Neumann boundary conditions related to Hagstrom techniques \cite{PML-siap,Hagstrom}. In addition, 
some first numerical experiments were presented to analyze the 
potentiality of the approach which was also tested for the one-dimensional linear Klein-Gordon equation in \cite{Antoine-review} again using finite-differences. 
The aim of the present paper is to address the application and
assessment of the PML approach to solve complex NKGE, both in the classical and non-relativistic regimes, considering Fourier pseudo-spectral  approximation schemes.

In Section \ref{sec:method}, we introduce the first-order PML formulation in the classical scaling, inspired by the works by 
Appel\"o and Kreiss in \cite{PML-wave}. To obtain an explicit
pseudo-spectral scheme, we use an exponential-wave integrator
discretization. In addition, we also derive
a second-order PML formulation of the NKGE, related to the standard
developments for time-harmonic wave 
and Schr\"odinger equations. The discretization is now
based on a Crank-Nicolson scheme combined with
the pseudo-spectral method and an efficient 
preconditioned Krylov (GMRES) subspace iterative solver. This
yields
a linearly implicit scheme with a cost similar to
the one for the  explicit first-order PML formulation. Since we use a pseudo-spectral approach,
the PDE with PML that has to be discretized requires some
smooth coefficients  but also some stability of the PML  for 
the tuning parameters. To this end, we introduce
some generalized singular Berm\'udez-type absorption functions
which are locally smoothed at the inner PML boundary. 
 This allows us to achieve  simultaneously a high (near-spectral) accuracy of 
the scheme and stable PML layers that are less sensitive to the tuning parameters,
for the second-order formulation. The first-order PML
formulation is shown to be less accurate  and is not recommended
in the pseudo-spectral framework. The conclusions are supported
by a thorough numerical study.
In Section \ref{SectionNR}, we extend and evaluate the numerical methods
for the non-relativistic scaling. This shows that the 
second-order formulation combined with the pseudo-spectral scheme
is again an efficient and accurate method for solving the NKGE, and most particularly thanks to the small $\varepsilon$ parameter appearing in the non-relativistic regime.
In Section \ref{Sec:rotating}, we extend the method to a two-dimensional
rotating NKGE that models the dynamics of the cosmic superfluid 
set in a rotating frame \cite{RKG,RKG-PRD}. We show in particular
that  the numerical method in the bounded domain is able
to simulate very accurately the dynamics of the vortices.
Finally, we conclude in Section \ref{sec:conclusion}.

\section{PML in classical scaling}\label{sec:method}
 To simplify the presentation, we first consider the one-dimensional real-valued case of the NKGE, i.e. $d=1$, $\bx=x$ and $u\in\bR$ in (\ref{KG model}), and we begin by looking at its classical dimensionless form \cite{wellpose,classical NKG}:
\begin{equation}\label{KG model 1d}
\left\{\begin{split}
&\partial_{tt} u(x,t)-\partial_{xx}u(x,t)+u(x,t)+\lambda u(x,t)^3=0,\quad t>0,\ x\in\bR,\\
&u(x,0)=u_0(x),\quad \partial_tu(x,0)=v_0(x),\quad x\in\bR.
\end{split}\right.
\end{equation}
Such classical scaling describes  the physical system with wave speed at the same order of the speed of light.
We shall consider two types of PML formulations for (\ref{KG model 1d}) in the sequel. One is a first-order formulation proposed by Appel\"o and Kreiss in \cite{PML-wave} for general nonlinear wave equations, and the other one is a second-order analogy of the PML widely applied e.g. for nonlinear Schr\"odinger equations  \cite{Antoine-review,Antoine-RNLS,PML-NLS}.
\subsection{PML-I}\label{sec:classical PMLI}
Let us begin by the modal ansatz construction with postulation
\begin{equation}\label{ansatz}
\displaystyle
u(x,t)=\fe^{k x+st}\fe^{\frac{k}{s+\alpha}
\int_0^x\sigma(\rho)d\rho},\quad
\partial_x\to\frac{s+\alpha}{s+\alpha+\sigma}\partial_x,
\end{equation}
where $\sigma=\sigma(x)\geq0$ is some chosen function known as the \emph{absorption function}, $s$ is interpreted as the variable in the Laplace-transform domain, and
$\alpha\geq0$ is a chosen  parameter. By assuming that $\re(s)\geq0$ and $k\in\bC$ in (\ref{ansatz}) satisfy the dispersion relation for the linear version of the model, e.g.
$\partial_{tt} u-\partial_{xx}u+u=0$,
a formulation of the PML equations as a  first-order system  has been proposed for a class of general nonlinear wave equations in \cite{PML-wave}. Here, we can directly apply such PML formulation (called PML-I) to the NKGE (\ref{KG model 1d}), leading to
\begin{equation}\label{KG model pml}
\left\{\begin{split}
&\partial_{tt} u-\partial_{xx}u+u+\lambda u^3=\sigma
\left[\eta_2-\partial_tu+\alpha u\right]+\partial_x(\sigma \eta_1),\quad t>0,\ x\in I^*,\\
&\partial_t\eta_1+(\sigma+\alpha)\eta_1+\partial_xu=0,\\
&\partial_t\eta_2+\alpha\eta_2+(\alpha^2+1)u+\lambda u^3=0,\\
&u(x,0)=u_0(x),\quad \partial_tu(x,0)=v_0(x),\quad\eta_1(x,0)=\eta_2(x,0)=0, \quad x\in I^*,\\
&u(-L^*,t)=u(L^*,t),\quad \eta_1(-L^*,t)=\eta_1(L^*,t),\quad \eta_2(-L^*,t)=\eta_2(L^*,t),\quad t\geq0,
\end{split}\right.
\end{equation}
where  $$I^*=(-L^*,L^*)\quad \mbox{with}\quad L^*=L+\delta$$ is a bounded interval for some $\delta>0$, and $\eta_1=\eta_1(x,t),\ \eta_2=\eta_2(x,t)$ are two auxiliary functions introduced in order to localize the layer equations in time.
The initial conditions for the two auxiliary functions in (\ref{KG model pml}) are suggested to be zero in \cite{PML-NLS}.

The absorption  function $\sigma(x)$ is required to satisfy
\begin{equation*}
\sigma(x)=0,\quad -L\leq x\leq L,\quad \mbox{and}\quad \sigma(x)>0,\quad L<|x|\leq L^*,
\end{equation*}
such that inside the physical domain $[-L,L]$ the NKGE and its solution remain the same. The damping effect takes place only inside the layer, i.e. $L<|x|\leq L^*$.
 The smoothness of $\sigma(x)$ in the layer region $L\leq|x|\leq L^*$, particularly at the interface $|x|=L$, determines the regularity of the PML solution $u(x,t)$. When the NKGE (\ref{KG model 1d}) is linear ($\lambda=0$), the absorbing layer (\ref{KG model pml}) is perfectly matching \cite{PML-wave}.
Similar PMLs of the first-order formulation have  also been derived for other hyperbolic systems \cite{PML-siap,PML-notes}.  The boundary conditions for (\ref{KG model pml}) at $|x|=L^*$ 
can be for example the homogeneous Dirichlet or Neumann boundary condition which  were considered in \cite{PML-wave}. In either case,  the reflected waves at the boundary  enter again the layer region and so are damped twice.  Here we adopt  periodic boundary conditions as suggested in \cite{Antoine-review,Antoine-RNLS}, so that the  Fourier pseudo-spectral method \cite{Shen} can be easily applied for the spatial discretization of the PML equations. Note that the periodicity of the boundary ensures that the waves hitting the outer boundary  enter the layer from the other side of the domain which is the same as the `round-way' damping.

\begin{figure}[hbt!]
$$\begin{array}{c}
\psfig{figure=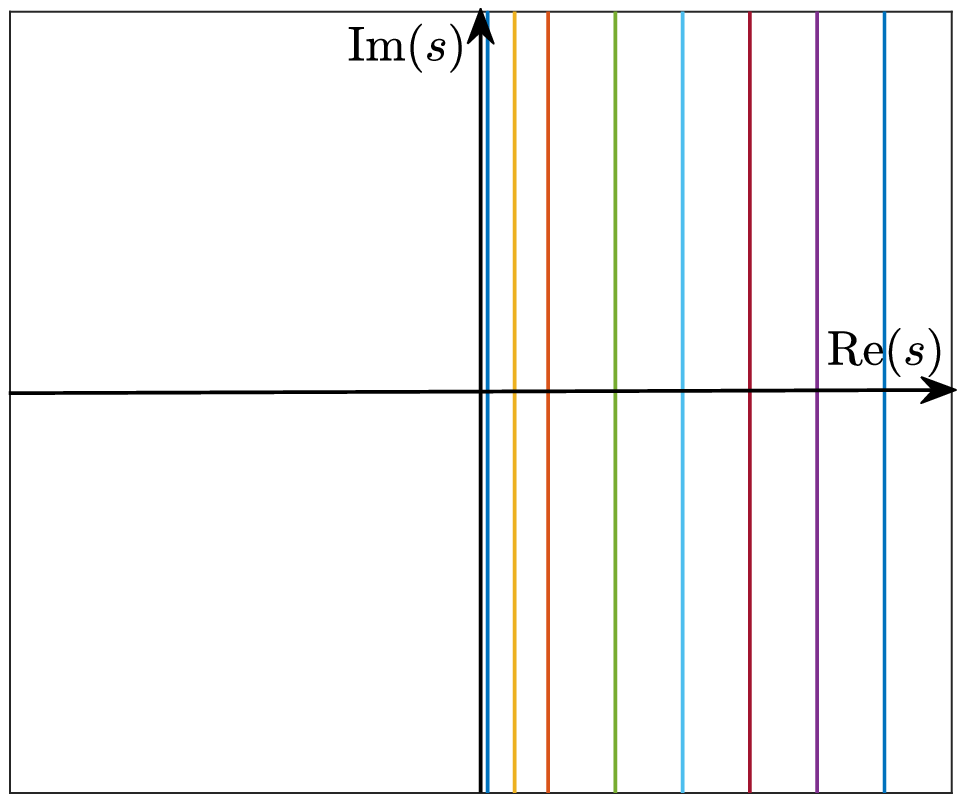,height=4cm,width=6cm}
\psfig{figure=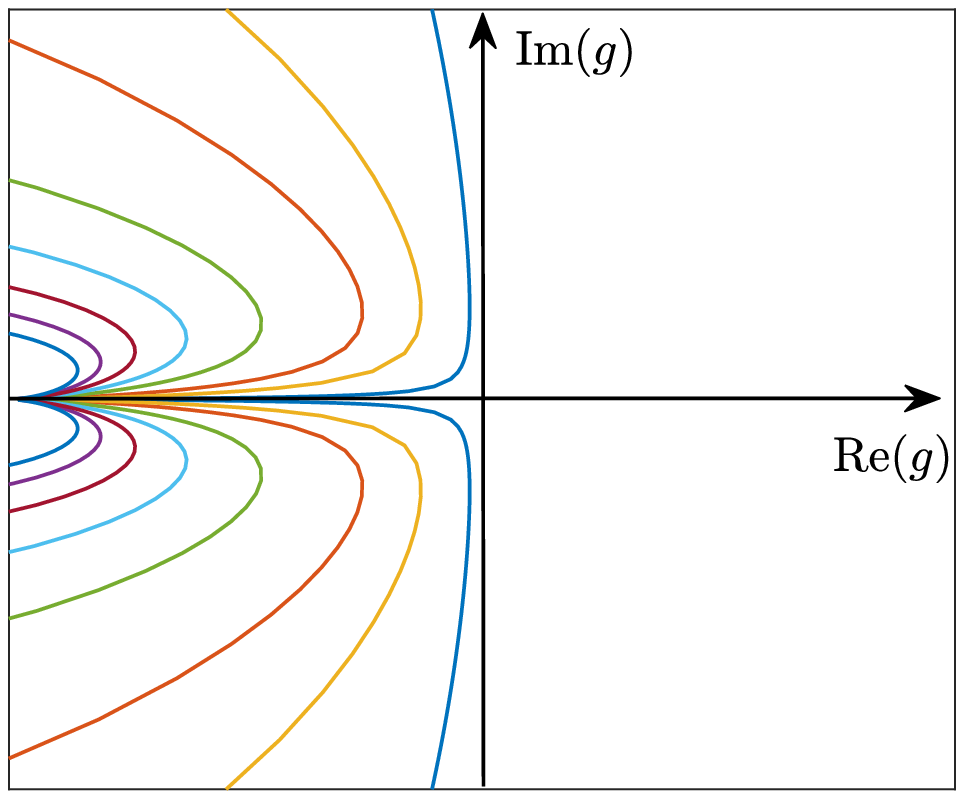,height=4cm,width=6cm}\\
\psfig{figure=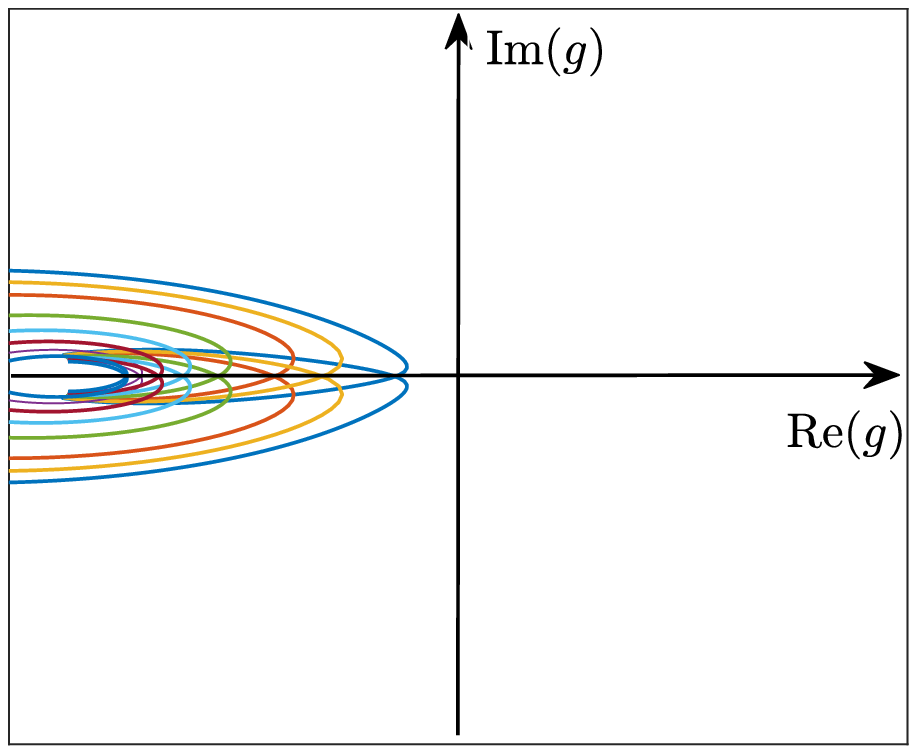,height=4cm,width=6cm}
\end{array}$$
\caption{The right half complex plane $\re(s)>0$ (top left) under mapping $g=-\sqrt{s^2+1}/(s+\alpha)$: $\alpha=0$ (top right) and $\alpha=0.5$ (bottom).}
\label{fig:PML-I R}
\end{figure}

Let us discuss briefly the damping effect brought by (\ref{ansatz}) for the Klein-Gordon equation and the choice of $\alpha$.  Firstly, we remark that with $\alpha=0$, the damping effect
$$\frac{k}{s}\int_0^x\sigma(\rho)d\rho$$ in (\ref{ansatz}) is a typical choice for
dispersionless materials, where the phase velocity $k/s$  is constant so that all waves  attenuate at the same rate \cite{PML-notes}. For the PML-I formulation for the NKGE (\ref{KG model pml}), we  consider the linear case, i.e. $\lambda=0$ in (\ref{KG model 1d}) or (\ref{KG model pml}) to illustrate the damping effect. By taking the modal solution $u=\fe^{kx+st}$, we have the dispersion relation for the original linear Klein-Gordon equation:
\begin{equation}\label{dispersion}
s^2=k^2-1\ \Rightarrow\ k=-\sqrt{s^2+1},\end{equation}
and so (\ref{ansatz}) reads
$$u(x,t)=\fe^{k x+st}\fe^{-\frac{\sqrt{s^2+1}}{s+\alpha}
\int_0^x\sigma(\rho)d\rho}.$$
This means that the modal solution $\fe^{k x+st}$ for PML-I (\ref{KG model pml}) is damped inside the layer $L<|x|\leq L^*$ by the absorption function $\sigma$ modulated by the factor $g(s)$:
\begin{equation}\label{gs}
g(s)=-\frac{\sqrt{s^2+1}}{s+\alpha},\quad \mbox{for}\quad \re(s)>0.
\end{equation}
We show in Figure \ref{fig:PML-I R} the mapping of the right half plane of the complex domain $\re(s)>0$ under the function $g(s)$ for $\alpha=0$ and $\alpha=0.5$. We can see that for $\re(s)>0$ and $\alpha\geq0$, we always have
$\re(g)<0$. However, it is difficult to say  if there exists or not an optimal choice of the shifting parameter $\alpha$ even for such a linear equation. The reason  is  that there could always be some $g(s)$ getting arbitrarily close to the imaginary axis for some $s$, with $\re(s)$ arbitrarily small.

To numerically solve the PML-I equations (\ref{KG model pml}), instead of the finite difference time-domain discretization in \cite{PML-wave}, our periodic boundary setup in (\ref{KG model pml}) is natural to apply the Fourier spectral discretization which is directly diagonalized in space, and we integrate in time under the framework of the exponential wave integrator  \cite{BDZ,EWI,Act}.  
Denoting $\tau=\Delta t>0$ as the time step and $t_n=n\tau$ for $n\in\bN$ as the discrete times, we shall derive in the following a fully explicit exponential wave integrator Fourier pseudo-spectral method for (\ref{KG model pml}).  For simplicity of notations, we will omit in the following the spatial variable $x$ in a space-time dependent  function, e.g. $u(t)=u(x,t)$.
Let $$v(t):=\partial_t u(t),\quad t\geq0.$$ Then, by applying the Duhamel's formula to (\ref{KG model pml}), we have
\begin{align*}
 u(t_{n+1})=&\cos(\langle\partial_x\rangle\tau)u(t_n)+\frac{\sin(\langle\partial_x\rangle\tau)}{
 \langle\partial_x\rangle}v(t_n)+\int_0^\tau\frac{
 \sin(\langle\partial_x\rangle(\tau-\rho))}{\langle\partial_x\rangle}f(t_n+\rho)d\rho,\\
 v(t_{n+1})=&-\langle\partial_x\rangle\sin(\langle\partial_x\rangle\tau)u(t_n)
 +\cos(\langle\partial_x\rangle\tau)v(t_n)+\int_0^\tau
 \cos(\langle\partial_x\rangle(\tau-\rho))f(t_n+\rho)d\rho,\\
 \eta_1(t_{n+1})=&\fe^{-(\sigma+\alpha)\tau}\eta_1(t_n)-\int_0^\tau\fe^{(\sigma+\alpha)(\rho-\tau)}
 \partial_x u(t_n+\rho)d\rho,\\
 \eta_2(t_{n+1})=&\fe^{-\alpha\tau}\eta_2(t_n)-\int_0^\tau\fe^{\alpha(\rho-\tau)}\left[(\alpha^2+1)u(t_n+\rho)
 +\lambda u(t_n+\rho)^3\right]d\rho,\quad n\geq0, \  x\in I^*,
\end{align*}
where we set $\langle\partial_x\rangle=
\sqrt{1-\partial_{xx}}$ and
\begin{align*}
  f(t_n+\rho):=\sigma(
 \eta_2(t_n+\rho)-v(t_n+\rho)+\alpha u(t_n+\rho))+\partial_x(\sigma \eta_1(t_n+\rho))-\lambda u(t_n+\rho)^3.
\end{align*}
Applying the trapezoidal rule to the above integrals, we obtain the following exponential-wave integrator Fourier pseudo-spectral (EWI-FP) scheme: by introducing $u^n\approx u(t_n), v^n\approx v(t_n),\eta_1^n\approx\eta_1(t_n),\eta_2^n\approx \eta_2(t_n)$, for $n\geq0$, we have
\begin{subequations}\label{EWI-FP}
\begin{align}
 u^{n+1}=&\cos(\langle\partial_x\rangle\tau)u^n+\frac{\sin(\langle\partial_x\rangle\tau)}{
 \langle\partial_x\rangle}v^n+\frac{
 \tau\sin(\langle\partial_x\rangle\tau)}{2\langle\partial_x\rangle}f^n,\label{EWI-FP a}\\
 \eta_1^{n+1}=&\fe^{-(\sigma+\alpha)\tau}\eta_1^n-\frac{\tau}{2}\left[\fe^{-(\sigma+\alpha)\tau}
 \partial_x u^{n}+\partial_x u^{n+1}\right],\label{EWI-FP b}\\
 \eta_2^{n+1}=&\fe^{-\alpha\tau}\eta_2^n-\frac{\tau}{2}\left[\fe^{-\alpha\tau}\left(\left(\alpha^2+1\right)u^n
 +\lambda (u^n)^3\right)+\left(\alpha^2+1\right)u^{n+1}
 +\lambda \left(u^{n+1}\right)^3\right],\label{EWI-FP c}\\
 v^{n+1}=&-\langle\partial_x\rangle\sin(\langle\partial_x\rangle\tau)u^n
 +\cos(\langle\partial_x\rangle\tau)v^n+
 \frac{\tau}{2}\left[\cos(\langle\partial_x\rangle\tau)f^n+f^{n+1}\right],\label{EWI-FP d}
\end{align}
\end{subequations}
with
$$
f^n=\sigma\left(
 \eta_2^n-v^n+\alpha u^n\right)+\partial_x(\sigma \eta_1^n)-\lambda (u^n)^3.
$$
Here, the scheme (\ref{EWI-FP}) and the numerical solution  $(u^n, v^n,\eta_1^n,\eta_2^n)$ are defined on the spatial grids: $x_j=-L^*+jh$ for $j=0,1,\ldots,N-1$ with some even integer  $N>0$ and mesh size $h=2L^*/N$.

The above EWI-FP scheme is fully explicit in time by computing (\ref{EWI-FP a}) to (\ref{EWI-FP d}).
The involved spatial differentiation operators can all be implemented efficiently under Fourier pseudo-spectral method \cite{Shen} by fast Fourier transform (FFT), resulting in a  computational cost  $O(N\log N)$. Concerning the stability of the  EWI-FP scheme, noting the term $\partial_xu$ in (\ref{EWI-FP b}) and the term $\partial_x(\sigma\eta_1^n)$ appearing in (\ref{EWI-FP a}) and (\ref{EWI-FP d}), we  have some unbalanced energy norms on both sides of (\ref{EWI-FP}). Therefore, we are expecting a CFL stability constraint:
$$\tau<C h,\quad \mbox{for some}\quad C>0.$$
For the accuracy, noticing the quadrature error from the trapezoidal rule and the interpolation error from the Fourier pseudo-spectral method, it can be analyzed similarly as in \cite{BDZ,EWI} that the global error of the above EWI-FP method for (\ref{KG model pml}) up to some finite time reads
$$O(\tau^2)+O(h^{m_0}),$$
for some $m_0>0$ which depends on the smoothness of the PML solution of (\ref{KG model pml}). We shall show numerically that by using a very smooth absorption function $\sigma(x)$ particularly with
high-order continuous derivatives at the interface  point $|x|=L$, the EWI-FP could reach a spectral-like accuracy for spatial discretization. Thanks to the spectral accuracy, we can use a rather large $h>0$ in practice so that the above CFL condition is not too restrictive.

\subsection{PML-II}\label{sec2 PML2}
For dispersive models like the nonlinear Schr\"odinger  equations (NSE), the phase velocity is naturally non-constant with respect to the wavelength. Therefore, instead of (\ref{ansatz}), an average damping effect has been widely considered in the postulation of a modal solution:
\begin{equation}\label{ansatz2}
u(x,t)=\fe^{k x+st}\fe^{kR
\int_0^x\sigma(\rho)d\rho},
\end{equation}
for some chosen $R\in\bC$ as for the NSE. Using  (\ref{ansatz2}) for the NKGE (\ref{KG model 1d}),  by directly
changing
$$\partial_x\rightarrow\frac{1}{1+R\sigma}\partial_x,$$
we consider the following second PML formulation (PML-II):
\begin{equation}\label{KG model pml2}
\left\{\begin{split}
&\partial_{tt} u-\frac{1}{1+R\sigma}\partial_x\left(\frac{1}{1+R\sigma}\partial_xu\right) +u+\lambda u^3=0,\quad t>0,\ x\in I^*,\\
&u(x,0)=u_0(x),\quad \partial_tu(x,0)=v_0(x), \quad x\in I^*,\\
&u(-L^*,t)=u(L^*,t),\quad t\geq0. 
\end{split}\right.
\end{equation}
Here, we also impose the periodic boundary condition in (\ref{KG model pml2}) for the convenience of spatial discretization.

\begin{figure}[hbt!]
$$\begin{array}{cc}
\psfig{figure=PMLs.eps,height=4cm,width=6cm}&
\psfig{figure=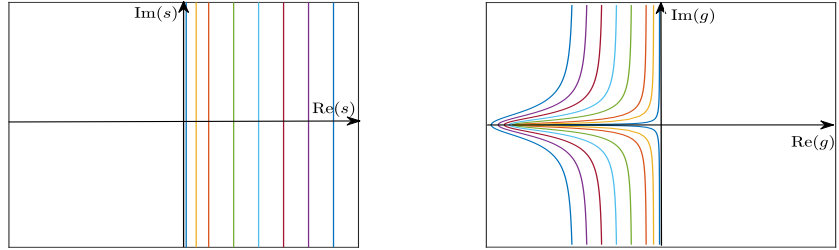,height=4cm,width=6cm}
\end{array}$$
\caption{The right half complex plane $\re(s)>0$ (left) under the mapping $g=-\sqrt{s^2+1}$ (right).}
\label{fig:PML-II R}
\end{figure}

An obvious advantage of PML-II (\ref{KG model pml2}) over PML-I (\ref{KG model pml}) is that PML-II does not require any auxiliary variables. We only need to solve one equation in (\ref{KG model pml2}), instead of three equations in (\ref{KG model pml}).
The PML-II can be considered as a second-order formulation, which is an analogy of the PML e.g. for NSE \cite{Antoine-review,Antoine-RNLS,PML-NLS}.
Note that for the NSE, the shifting parameter $R$ is usually chosen as $R=\fe^{i\pi/4}$ \cite{PML-NLS}. However, for the NKGE (\ref{KG model pml2}), this choice generates some instabilities, and we need to consider $R>0$. Indeed, assuming again $\lambda=0$ in (\ref{KG model 1d}) with the modal solution $u=\fe^{kx+st}$, by the dispersion relation (\ref{dispersion}), we have in  (\ref{ansatz2}):
$$u(x,t)=\fe^{k x+st}\fe^{-R\sqrt{s^2+1}
\int_0^x\sigma(\rho)d\rho}.$$
Then in order to make the solution decay inside the layer, we ask for
$$\re\left(-R\sqrt{s^2+1}\right)<0,\quad \forall\ \re(s)>0.$$
We illustrate in Figure \ref{fig:PML-II R} the mapping $g(s)=-\sqrt{s^2+1}$ of the right half complex plane. It can be clearly seen that the choice of $R$ for the PML (\ref{KG model pml2}) is restricted  to a positive real number, otherwise any non-zero imaginary part of $R$  leads to  some $\re\left(-R\sqrt{s^2+1}\right)>0$ which triggers instability. This will be verified by numerical tests later.
Since $R$ must be positive, it can be viewed  in (\ref{KG model pml2})  only as a modification of  $\sigma_0>0$, namely the strength of the absorption function $\sigma(x)$ (see e.g. (\ref{sigma}) or (\ref{Bermudez})). Therefore, we will only consider the impact from the choice of $\sigma_0>0$ on the performance of the PML-II with $R>0$ fixed in this section.

For the numerical discretization of the PML-II formulation (\ref{KG model pml2}), we use the following direct semi-implicit finite-difference Fourier pseudo-spectral (FD-FP) method:
\begin{subequations}\label{FD-FP}
\begin{align}
  &\frac{u^{n+1}-2u^n+u^{n-1}}{\tau^2}+\frac{A}{2}
  (u^{n+1}+u^{n-1})+\frac{1}{2}
  (u^{n+1}+u^{n-1})+\lambda (u^n)^3=0,\quad n\geq1,
  \label{FD-FP a}\\
  &u^1=u_0+\tau v_0-\frac{\tau^2}{2}\left[Au_0+u_0+\lambda u_0^3\right],\label{FD-FP b}
\end{align}
\end{subequations}
where $u^n\in\bR^N$ and the matrix
$$A=-d_0D_1d_0D_1,\quad \mbox{with}\quad d_0=\mathrm{diag}(1/(1+R\sigma)),$$
and $D_1\in \bC^{N\times N}$ is the matrix representation \cite{Shen} of the approximate first-order differentiation operator $\partial_x$ in the periodic domain $I^*$. The time averaging in (\ref{FD-FP a}) indeed avoids any CFL type stability condition.
The starting value (\ref{FD-FP b}) is obtained directly from the Taylor  expansion for the three-level scheme (\ref{FD-FP a}).
Note that (\ref{FD-FP a}) reads explicitly
\begin{align}\label{FD-FP c}
  u^{n+1}=-u^{n-1}
  +G^{-1}\left[\frac{2}{\tau^2}u^n
  -\lambda (u^n)^3\right],\quad n\geq1,\quad
  G:=\left[\left(\frac{1}{\tau^{2}}+\frac{1}{2}\right)\mathbb{I}
  +\frac{1}{2}A\right],
\end{align}
where $\mathbb{I}$ is the identity matrix in $\bR^{N\times N}$.
The matrix $G^{-1}$ is dense here, but it can be pre-computed with the chosen $\tau$ and $N$ once for all.
The computational cost of FD-FP at each time level is therefore $O(N^2)$ if one performs the previous matrix-vector product. However, by borrowing the idea from \cite{Antoine-RNLS}, the scheme (\ref{FD-FP c}) can be implemented efficiently by incorporating the GMRES solver with FFT. Here we briefly describe such a strategy for evaluating (\ref{FD-FP c}). The equation (\ref{FD-FP c}) equivalently reads
\begin{equation}\label{FD-FP gmres}
\left\{
  \begin{array}{l}
  \displaystyle
  u^{n+1}=-u^{n-1}
  +w^n,\quad n\geq1,\\
   \displaystyle  G w^n   =\frac{2}{\tau^2}u^n -\lambda (u^n)^3. 
  \end{array}
  \right.
  \end{equation}
The second equation of (\ref{FD-FP gmres}) can be treated as a linear system with  unknown $w^n$, which can be solved iteratively e.g. by  GMRES \cite{saad} since it is a matrix-free  solver. Within the GMRES iteration, the involved spatial derivatives in  $G$  can be straightforwardly implemented by applying the  FFT to  discretize the operators in the iteration process. The GMRES convergence can be strongly improved by introducing a Fourier diagonalizable preconditioner for $\sigma=0$, similarly as in \cite{Antoine-lib1,Antoine-lib2,Antoine-lib3}, i.e.
$$\mathcal{P}G w^n
  =\mathcal{P}\left[\frac{2}{\tau^2}u^n
  -\lambda (u^n)^3\right],$$
  with
\begin{equation}\label{preconditioner}\mathcal{P}=\left(\frac{1}{\tau^{2}}
  +\frac{1}{2}-\frac{\partial_{xx}}{2}\right)^{-1}.\end{equation}
Let us remark that the preconditioning  operator $\mathcal{P}$ can also be directly  implemented at low cost by FFT, since it has constant coefficients and is therefore diagonalizable. The computational cost of the FD-FP method is then  nearly $O(N\log N)$.
The global error  of FD-FP is similar to EWI-FP, which is  $O(\tau^2)+O(h^{m_0})$.

\subsection{Absorption function} In this subsection, we discuss the choices for the absorption function.   To enhance the efficiency of the PMLs, we look for spectral accuracy in the related spatial discretization. Let us recall that getting a smooth absorption function is crucial to get the spectral accuracy, where the parameter $m_0>0$ specifies the accuracy induced by the PML absorption function on the error estimates. We shall consider two types of absorption functions $\sigma(x)$ for the PMLs: a classical polynomial choice and a singular function case.

The first possibility is to choose a high degree \emph{polynomial type} function as proposed in \cite{PML-wave} for wave equations
\begin{equation}\label{sigma}
  \sigma_{P}(x):=\left\{\begin{split}
  &\sigma_0\left[1-\left(\frac{x-L^*}{\delta}\right)^2\right]^8,\quad
  L\leq |x|\leq L^*=L+\delta,\\
  &0,\qquad\qquad\qquad\qquad\qquad\  \mbox{else}.\end{split}\right.
\end{equation}
It can be viewed as a smoothed version of the commonly used quadratic or cubic absorbing functions for PMLs \cite{Antoine-RNLS,Classical-sigma1,Classical-sigma2,PML-NLS}. Thanks to the large power, (\ref{sigma}) is very smooth at the interfaces  $|x|=L$ so that one can maintain the high order accuracy from the Fourier spectral discretization in space.

Nevertheless, it is well-known since the work by
Berm\'{u}dez \textit{et al.} \cite{Bermudez,BermudezJCP,BermudezSISC}
that singular type functions are much more adapted than the polynomial functions for PMLs applied to wave-like equations. More specifically, the $\sigma$ functions must be such that
\begin{equation}
\int_L^{L^*}\sigma(x)dx=+\infty.
\end{equation}
In particular, they are less subject to parameter tuning problems  than in the polynomial case (\ref{sigma}) to fix the values of $\sigma_0$ and $\delta$,
in particular thanks to the speed $c$. This is crucial  in our case when prospecting the non-relativistic limit of the NKGE in the next section.

The concrete examples provided by Berm\'{u}dez \textit{et al.} do not correspond however to smooth functions, and therefore they lead to a limited accuracy for the pseudo-spectral approach. Here, we propose   some simple locally corrected smoother
 \emph{Berm\'{u}dez type}  functions as follows: for $k=-1,0,1,\ldots$, we define the absorption functions
\begin{equation}\label{Bermudez}
  \sigma_{B_k}(x):=\left\{\begin{split}
  &\sigma_0\beta_k(|x|-L^*),\quad
  L\leq |x|\leq L^*=L+\delta,\\
  &0,\qquad\qquad\qquad \quad  \mbox{else},\end{split}\right.
\end{equation}
where we set
\begin{align*}
  \beta_k(|x|-L^*)=\beta_{-1}(|x|-L^*)-\sum_{j=0}^{k}\frac{1}{j!}
  \frac{d^j \beta_{-1}(z)}{dz^j}\bigg|_{z=-\delta}(|x|-L)^j,\quad \beta_{-1}(z)=-\frac{1}{z},\quad z\in\bR.
\end{align*}
The function $\sigma_{B_{-1}}$ (i.e. $k=-1$ in (\ref{Bermudez})) is the original absorption function proposed   in \cite{Bermudez,BermudezJCP,BermudezSISC}. It contains both discontinuity and singularity over the domain. For  $k\geq0$ as above, we have introduced a truncated Taylor series expansion of $\beta_{-1}$ in $\beta_k$ so that the function  $\sigma_{B_k}(x)$ has $k$-th order continuous derivatives on the interval $|x|<L^*$, which is expected for the spectral scheme.

The polynomial type choice $\sigma(x)=\sigma_P(x)$ is the classical absorption function which is bounded on the whole interval $|x|\leq L+\delta$. Here, $\sigma_0=\|\sigma_{P}\|_{\infty}>0$ and $\delta>0$ are known as the \emph{strength} and \emph{thickness} of the PML, respectively.
The function $\sigma_P$ can be directly applied to both  PML formulations, i.e. PML-I (\ref{KG model pml}) and PML-II (\ref{KG model pml2}).
The Berm\'{u}dez type choice $\sigma(x)=\sigma_{B_k}(x)$ has singularities at the two boundary points $|x|=L^*$. The strength of such absorption function can be considered as infinite in some sense. This indeed is a problem for the PML-I formulation (\ref{KG model pml}), where in the equations  we can see that  the infinite function value of $\sigma_{B_k}(x)$ at the boundary is inconsistent with the periodic boundary condition mathematically. Moreover,  when the propagating waves are getting close to the boundary, the large values of $\sigma_{B_k}(x)$ near the boundary could make the right-hand-side of the first equation in (\ref{KG model pml}) stiff which causes numerical instability problems in schemes like EWI-FP. Thus, for PML-I (\ref{KG model pml}), we shall only consider the polynomial absorption function $\sigma_P$ (\ref{sigma}). As for the PML-II formulation (\ref{KG model pml2}), in contrast, $\sigma_{B_k}(x)$ can be directly applied in the equation, where we have a well-defined function $S(x)$ on  $I^*$ with zero boundary values, i.e. in  (\ref{KG model pml2}) with $\sigma=\sigma_{B_k}$ for any $k\in\bN$,
$$S(x):=\frac{1}{1+R\sigma_{B_k}(x)},\quad |x|\leq L^*,\quad\mbox{where}\quad S(\pm L^*)=0.$$
Therefore, we shall consider both types of absorption functions for PML-II (\ref{KG model pml2}).

\subsection{Numerical results} In this subsection, we conduct some numerical experiments to illustrate the performance of the two presented types of PML formulations as well as the proposed numerical schemes.
We shall denote in the following $u=u(x,t)$ as the exact solution of the NKGE (\ref{KG model 1d}), $u_{\textrm{pml}}=u_{\textrm{pml}}(x,t)$ as the exact solution of the PML-I (\ref{KG model pml}) or PML-II (\ref{KG model pml2}), and $u_{\textrm{pml}}^n\approx u_{\textrm{pml}}(x,t_n)$ as the corresponding numerical solution from EWI-FP (\ref{EWI-FP}) or FD-FP (\ref{FD-FP}).

\begin{figure}[hbt!]
$$\begin{array}{cc}
\psfig{figure=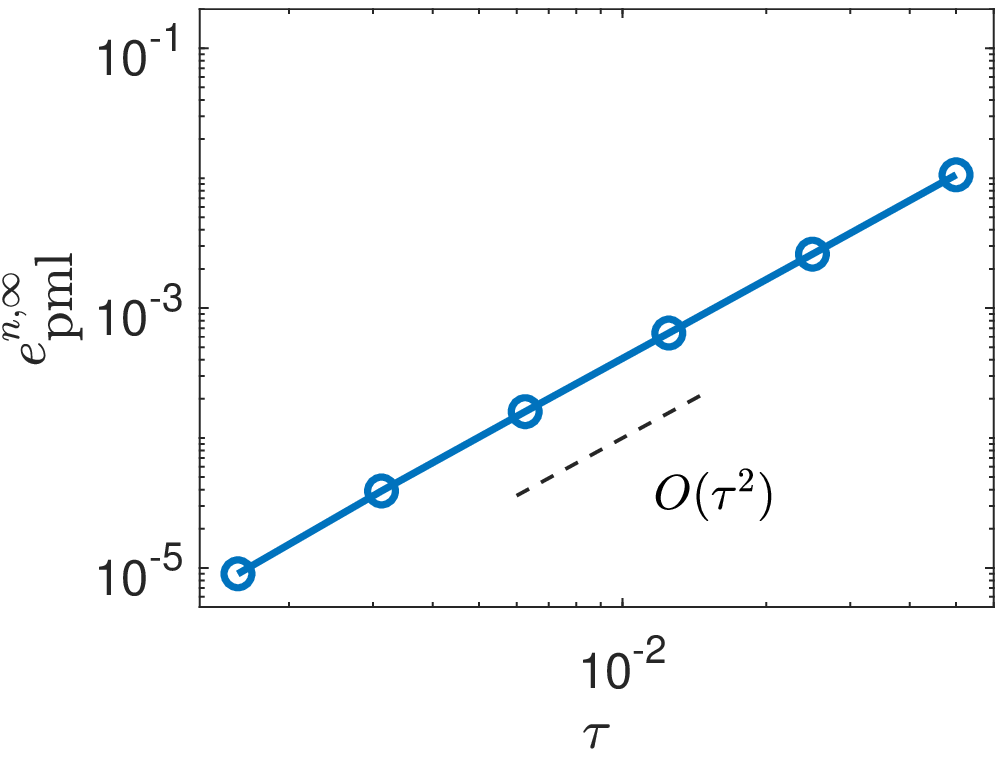,height=4cm,width=6cm}&
\psfig{figure=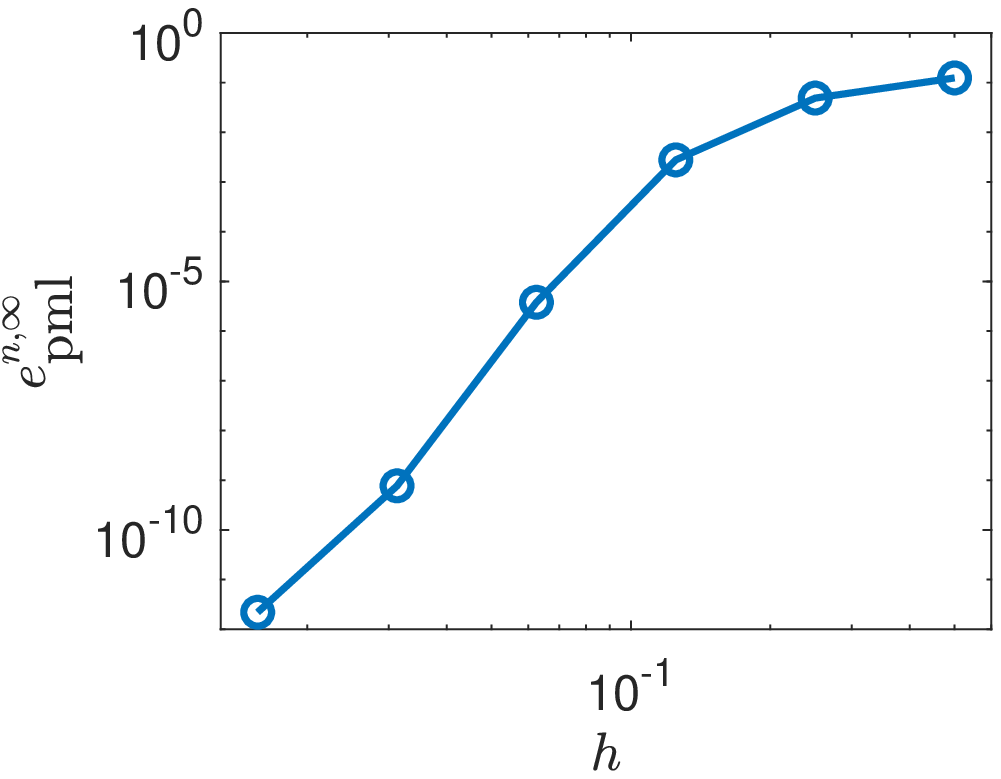,height=4cm,width=6cm}
\end{array}$$
\caption{The temporal (left) and spatial (right) $e^{n,\infty}_{\textrm{pml}}$-errors of EWI-FP for PML-I.}
\label{fig:convergencePML-I}
\end{figure}

\noindent \textbf{PML-I.}
We begin with the  PML-I formulation (\ref{KG model pml}). Within these tests, we   fix the absorption function as the high order polynomial type function (\ref{sigma}), i.e.
$$\sigma(x)=\sigma_P(x),\quad |x|\leq L^*.$$
First of all, we verify the accuracy of our numerical solver EWI-FP (\ref{EWI-FP}) for approximating (\ref{KG model pml}) with
\begin{equation}\label{example}
 \lambda=1, \quad u_0=5\fe^{-x^2},\quad v_0=\frac{1}{2}\mathrm{sech}(x^2),\end{equation}
and fix $\alpha=0$, $\sigma_0=8,\,\delta=0.5,\,L=4$ for the layer. Then we compute the relative maximum  error  \begin{equation}\label{errorinfty}
e^{n,\infty}_{\textrm{pml}}:=\frac{\|u_{\textrm{pml}}^n-u_{\textrm{pml}}\|_{L^\infty(I)}}{
\|u_{\textrm{pml}}\|_{L^\infty(I)}}
\end{equation}
on the physical domain $I=(-L,L)$ at some time $t=t_n$, where the reference solution $u_{\textrm{pml}}$ of (\ref{KG model pml}) is obtained here numerically by using EWI-FP under a fine mesh with $\tau=1\times 10^{-4}$ and $h=1/128$.
The temporal and spatial errors of the EWI-FP scheme (\ref{EWI-FP}) at $t=4$ are shown in Figure \ref{fig:convergencePML-I}, where we fix $h=1/32$ for the temporal error (left) and  $\tau=1\times 10^{-4}$ for the spatial error (right).
From the numerical results in Figure \ref{fig:convergencePML-I}, we can clearly see that the proposed EWI-FP (\ref{EWI-FP}) for solving (\ref{KG model pml}) with (\ref{sigma}) converges in time with a second-order accuracy and converges in space with a near spectral accuracy.

\begin{figure}[hbt!]
\centerline{Linear case: $\lambda=0$}
$$\begin{array}{cc}
\psfig{figure=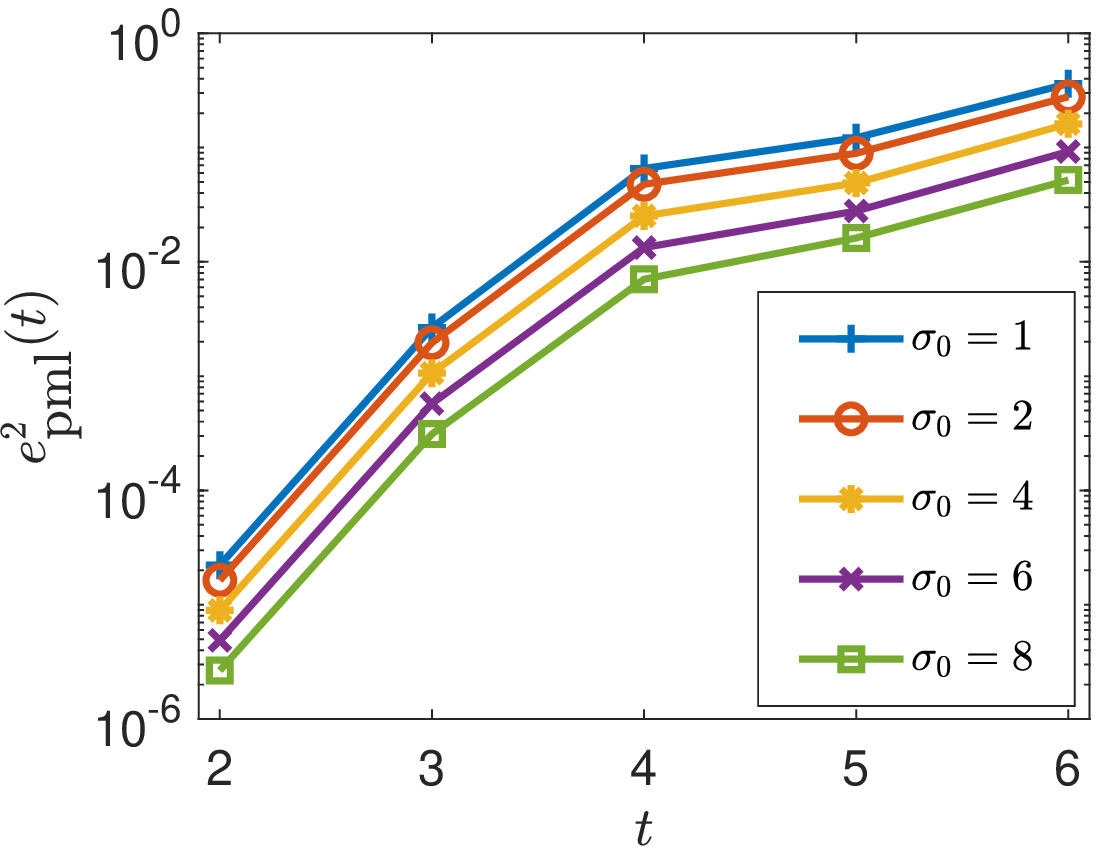,height=4.0cm,width=6cm}&
\psfig{figure=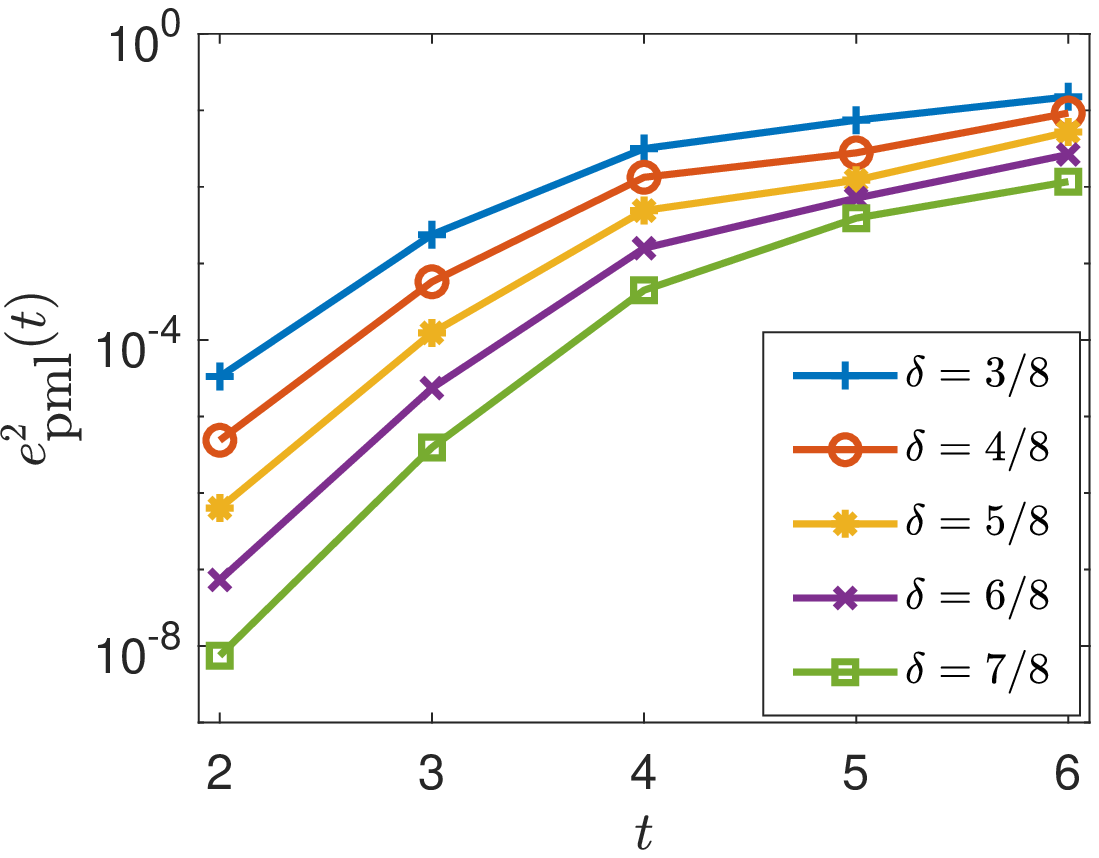,height=4.0cm,width=6cm}
\end{array}$$
\centerline{Nonlinear case: $\lambda=1$}
$$\begin{array}{cc}
\psfig{figure=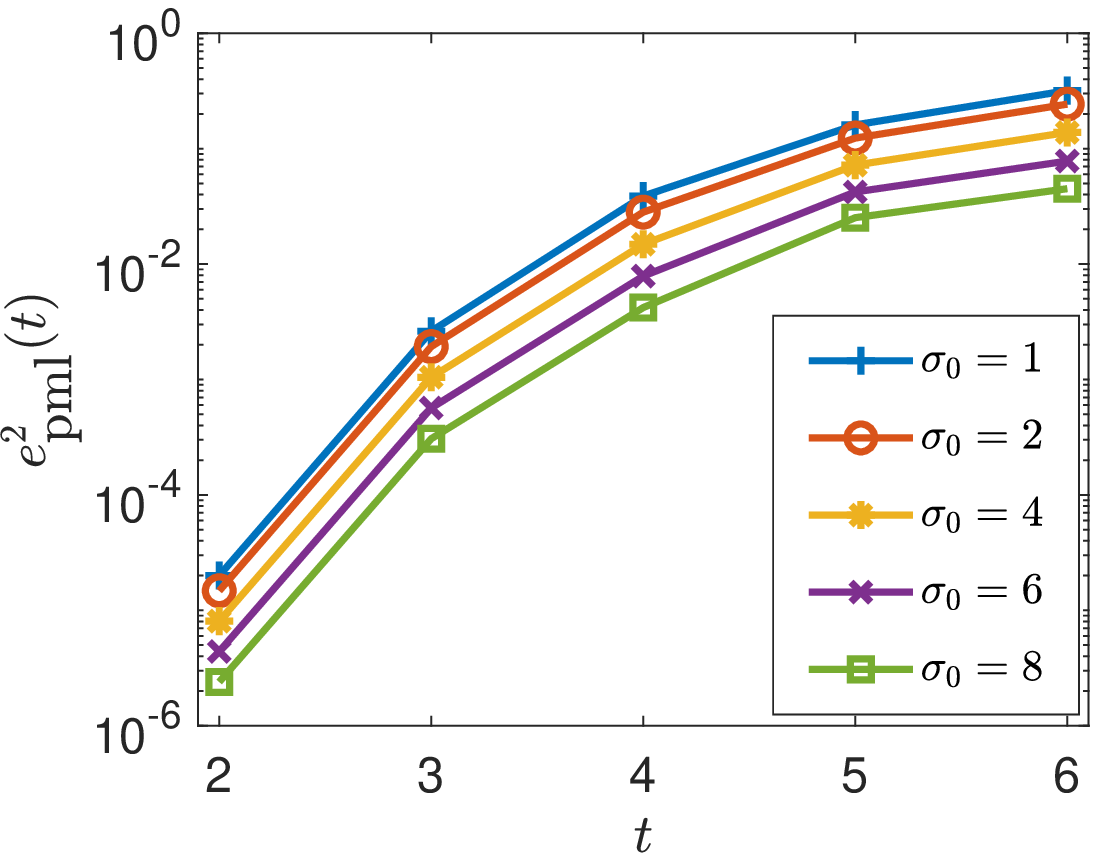,height=4cm,width=6cm}&
\psfig{figure=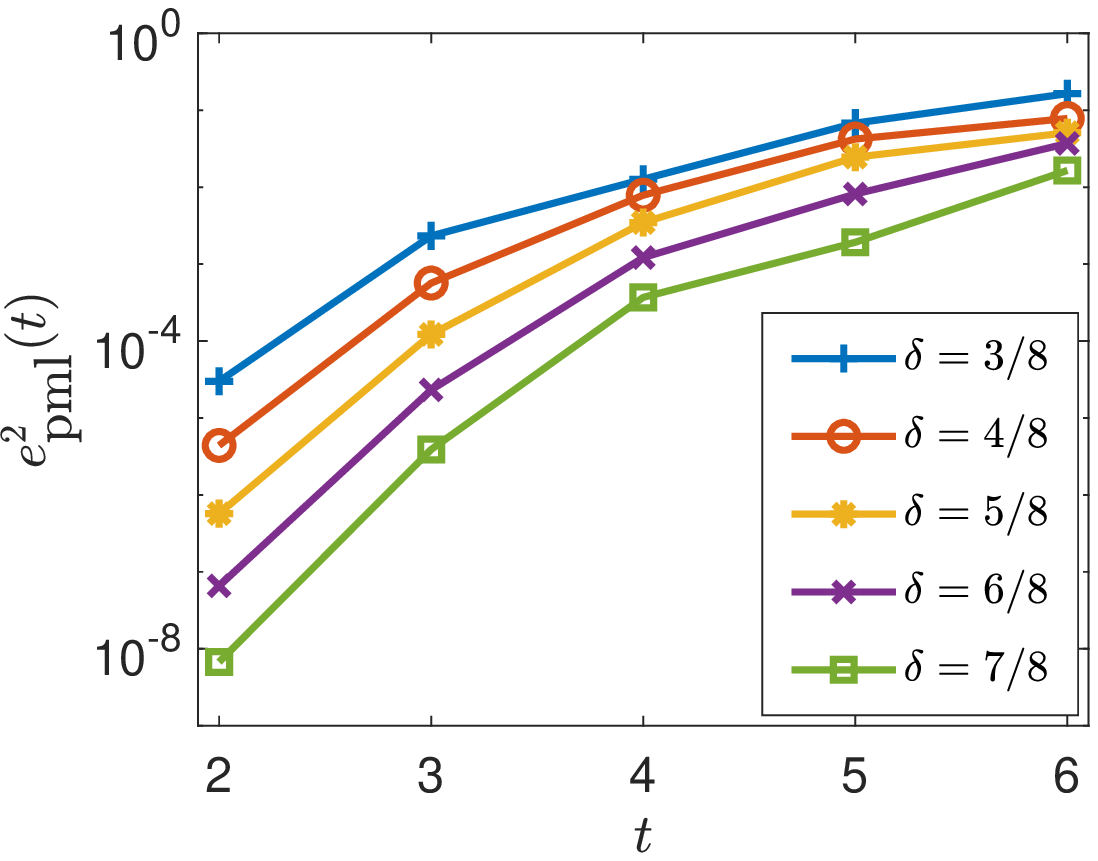,height=4cm,width=6cm}
\end{array}$$
\caption{The PML-I $e^2_\textrm{pml}$-error  as a function of time for $\sigma=\sigma_P$ in the linear and nonlinear cases: $\delta=1/2$ and varying strength $\sigma_0$ (left),
 and $\sigma_0=6$ and varying thickness $\delta$ (right).}
\label{fig:testNKG}
\end{figure}

\begin{figure}[hbt!]
$$\begin{array}{cc}
\psfig{figure=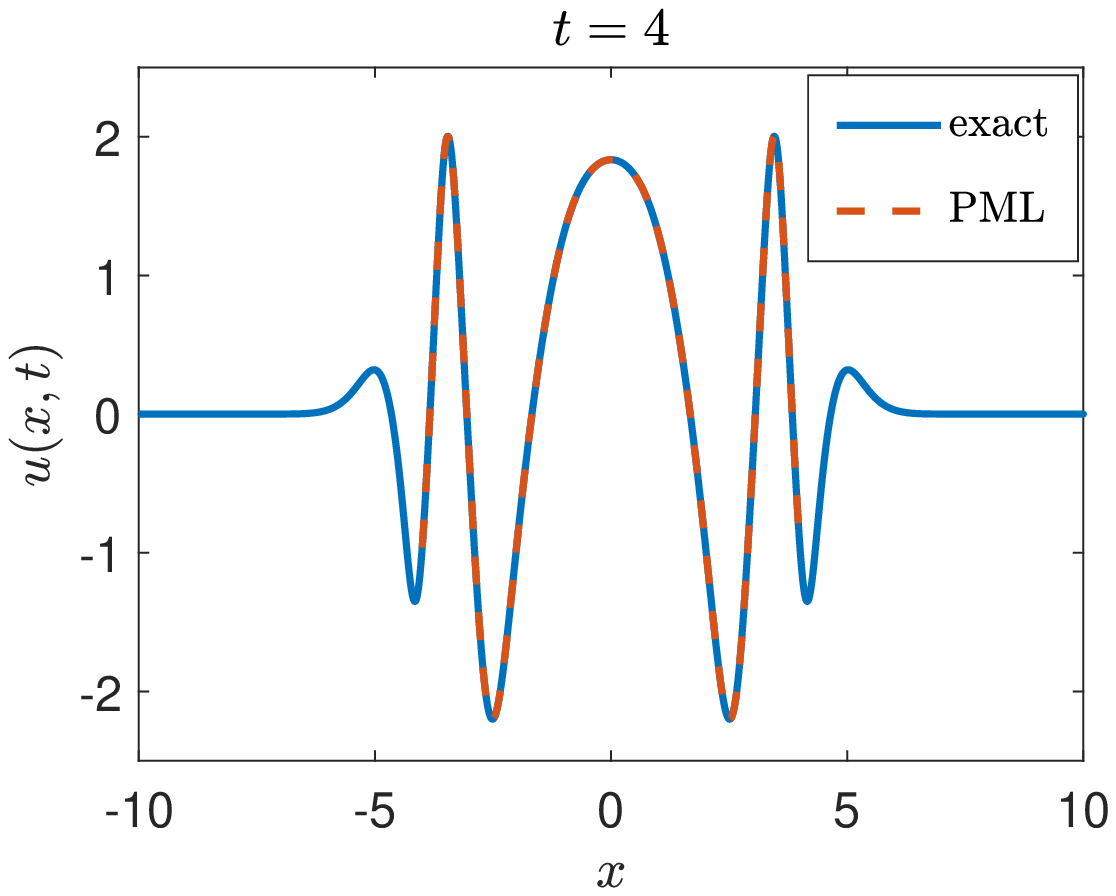,height=4.0cm,width=6cm}&
\psfig{figure=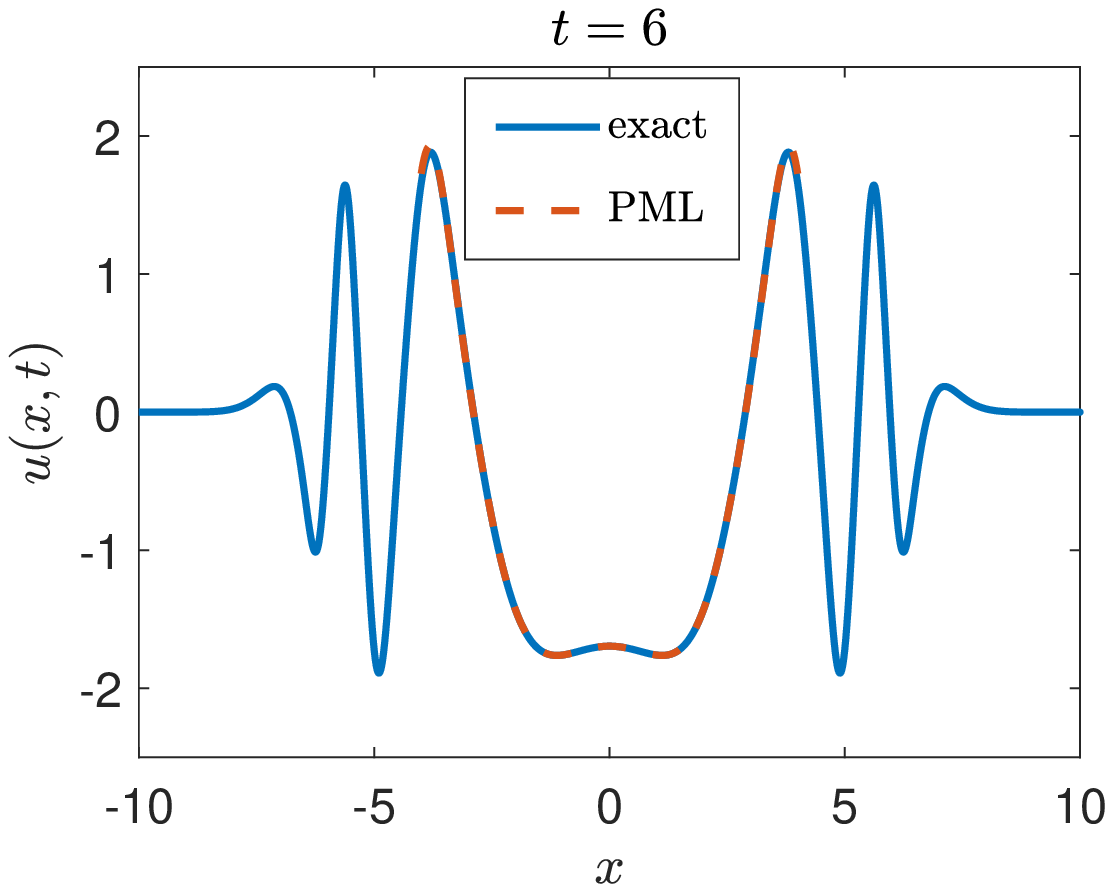,height=4.0cm,width=6cm}
\end{array}$$
\caption{The exact solution and the PML-I solution at $t=4$ (left) and $t=6$ (right).}
\label{fig:T6}
\end{figure}

Next, we test the error of the PML-I formulation (\ref{KG model pml}) for approximating the solution to the NKGE (\ref{KG model 1d}). We  use EWI-FP with a very fine mesh, e.g. $\tau=1\times 10^{-4}$ and $h=1/128$, such that the numerical discretization error is rather negligible. We compute the following relative PML error as a function of time $t>0$
\begin{equation}\label{err}
e^2_\textrm{pml}(t):=
\frac{\|u(\cdot,t)-u_{\textrm{pml}}(\cdot,t)\|_{L^2(I)}}{
\|u(\cdot,t)\|_{L^2(I)}},\quad I=(-L,L),
\end{equation}
under some different strength $\sigma_0$ or thickness $\delta$ parameters for the absorption function (\ref{sigma}).
We use the same numerical example (\ref{example}) as above, and we measure the error (\ref{err}) on the physical domain
$I=(-4,4)$.
The exact solution $u(x,t)$ of (\ref{KG model 1d}) is obtained by using the EWI-FP scheme (for $u$ only) to solve (\ref{KG model 1d}) accurately on a large enough interval, e.g. $(-16,16)$, where the solution within the computational time is still away from the boundary.
The relative errors (\ref{err}) as a function of time under different parameters for (\ref{sigma}) are plotted in Figure \ref{fig:testNKG}.
To compare the performance of the PML-I (\ref{KG model pml}) in the nonlinear case of the Klein-Gordon equation with its linear case, we also include in Figure \ref{fig:testNKG} the error result (\ref{err}) under the same example (\ref{example}) but with $\lambda=0$ in (\ref{KG model 1d}) and (\ref{KG model pml}).
For (\ref{example}), the profiles of the exact solution of (\ref{KG model 1d}) and the PML solution of (\ref{KG model pml}) at time $t=4$ and $t=6$ are reported in Figure \ref{fig:T6} for the PML parameters $\sigma_0=8$ and $\delta=7/8$.

From the numerical results in Figures \ref{fig:testNKG} and \ref{fig:T6}, we can see that the PML-I (\ref{KG model pml}) is effectively approximating the solution of NKGE (\ref{KG model 1d}) within the physical domain (see Figure \ref{fig:T6}). The PML-I formulation in the linear case, i.e. $\lambda=0$ in (\ref{KG model pml}), was proposed to be perfectly matched for the linear Klein-Gordon equation in \cite{PML-wave}. In the nonlinear case, i.e. $\lambda\neq0$ in NKGE (\ref{KG model 1d}), although PML-I (\ref{KG model pml}) is not theoretically perfectly matched, its performance is very close to the linear case (cf. the first row and the second row in Figure \ref{fig:testNKG}). With increasing strength of the absorption function and thickness of the fictitious layer, the accuracy of PML-I could be improved (see Figure \ref{fig:testNKG}).

\begin{remark}\label{remark0}
  With a usual low order polynomial absorption function $\sigma$, e.g. the quadratic or cubic polynomial in \cite{Antoine-RNLS,Classical-sigma1,Classical-sigma2,PML-NLS}, the performance of the PML-I (\ref{KG model pml}) is similar to that in Figure \ref{fig:testNKG}. However, the spatial convergence rate of the EWI-FP method would only be at the second- or third-order. These numerical results are omitted here for brevity.
\end{remark}

\begin{remark}
The choice of $\alpha\geq0$ barely impacts the performance of PML-I (\ref{KG model pml}), where the PML errors (\ref{err}) when using different $\alpha$ are very close. This is not surprising within the context of the analysis of the damping effect in (\ref{gs}).  Therefore,  we only present the numerical results for $\alpha=0$.
\end{remark}

\medskip
\noindent\textbf{PML-II.}
Under the same numerical example (\ref{example}),
we present the numerical experiments for PML-II given by (\ref{KG model pml2}). The physical domain is fixed as $I=(-4,4)$. For these tests, we will consider both the polynomial choice (\ref{sigma}) and the Berm\'{u}dez type function (\ref{Bermudez}) as the absorption function for PML-II.

\begin{figure}[hbt!]
$$\begin{array}{cc}
\psfig{figure=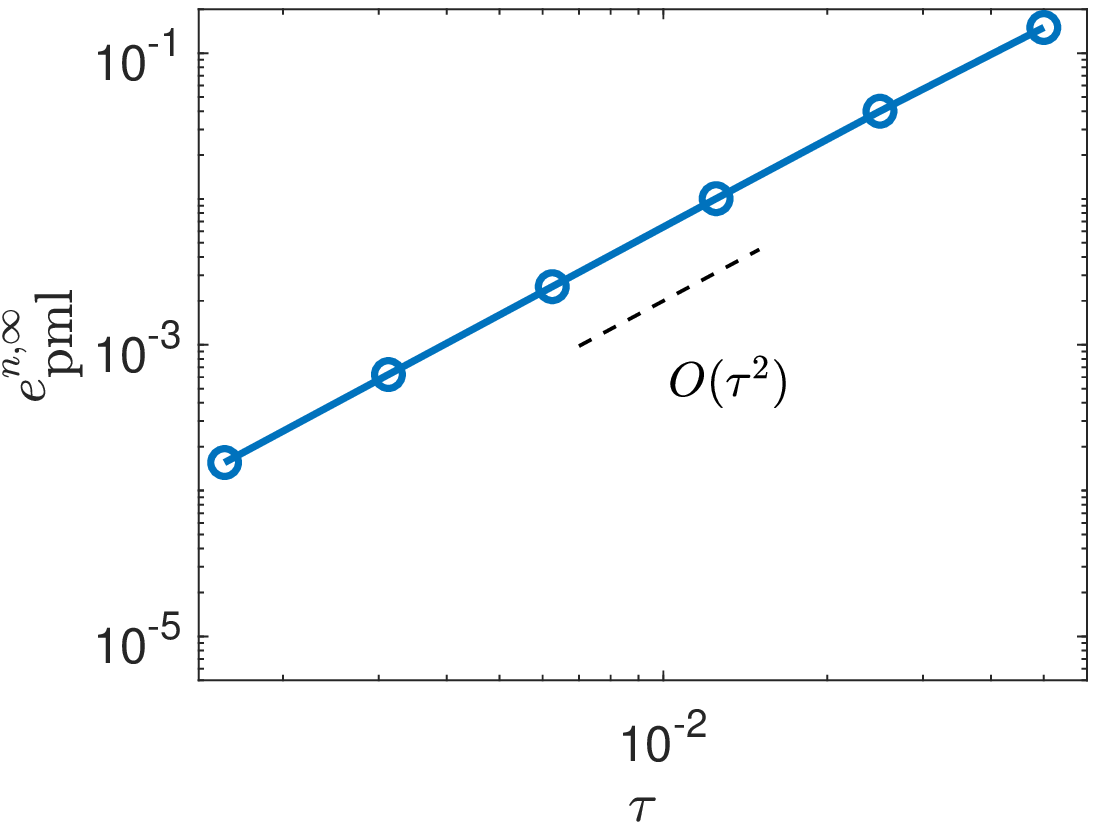,height=4cm,width=6cm}&
\psfig{figure=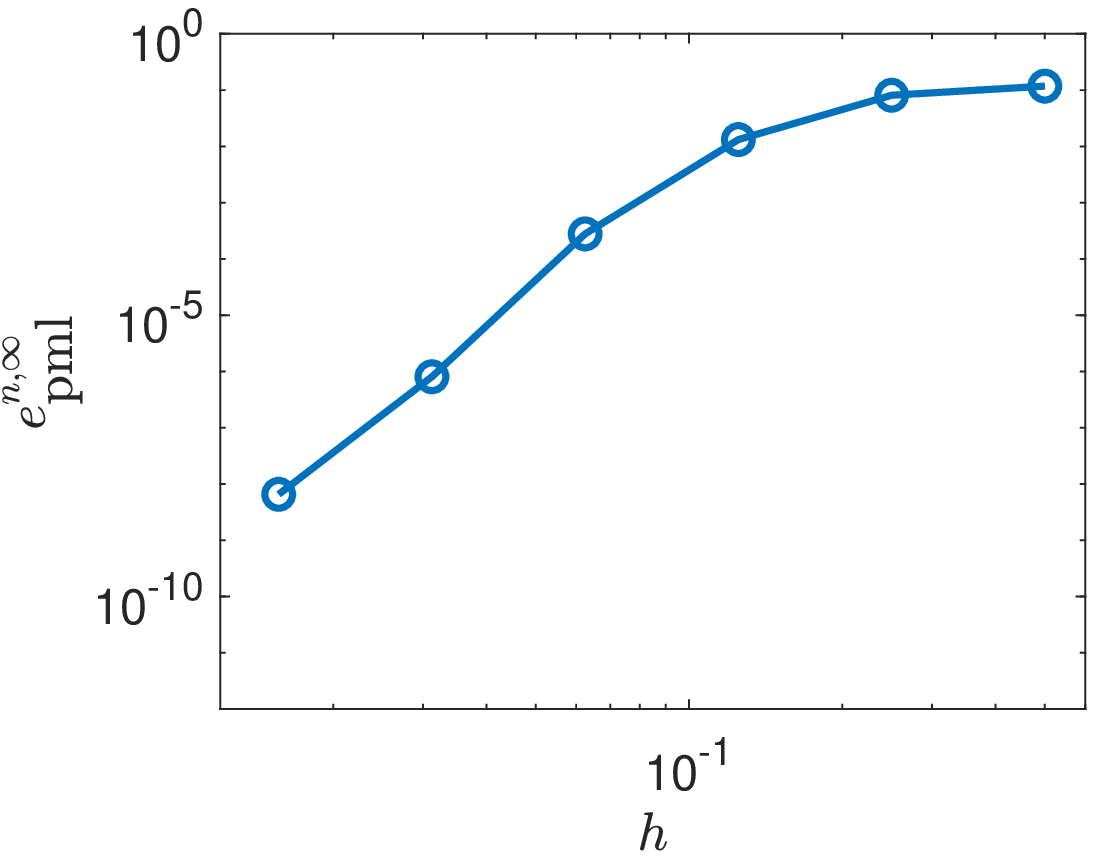,height=4cm,width=6cm}
\end{array}$$
\caption{The temporal (left) and spatial (right) $e^{n,\infty}_{\textrm{pml}}$-errors of FD-FP for PML-II  with $\sigma=\sigma_P$.}
\label{fig:convergencePML-II}
\end{figure}

\begin{figure}[hbt!]
$$\begin{array}{cc}
\psfig{figure=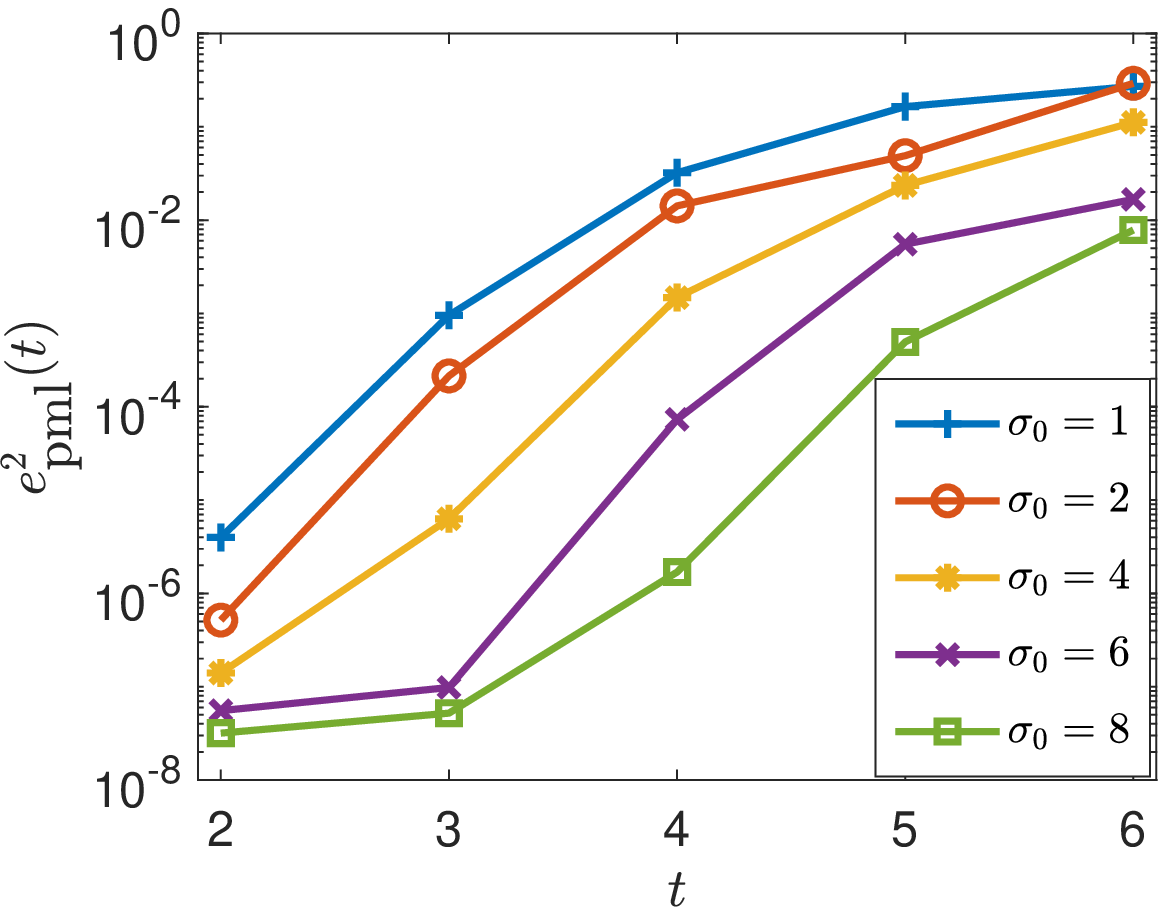,height=4cm,width=6cm}&
\psfig{figure=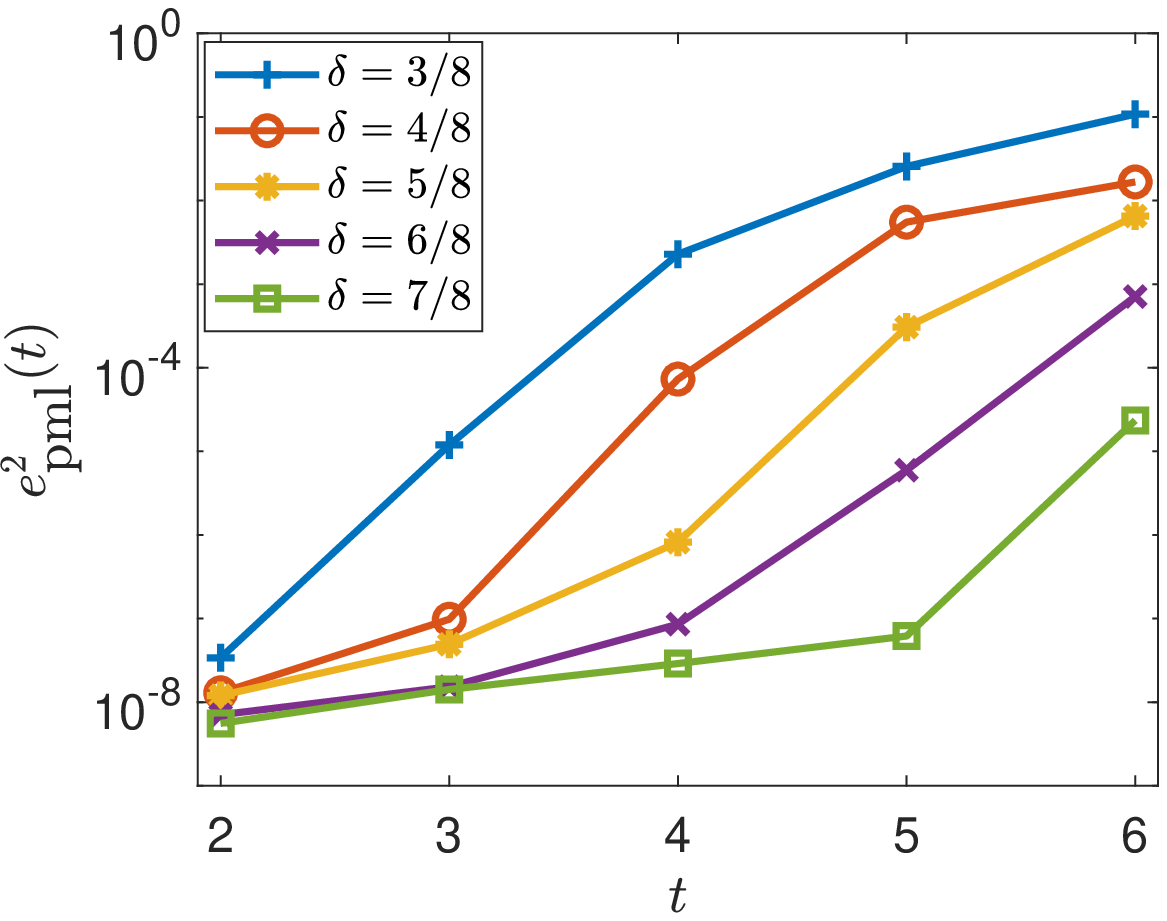,height=4cm,width=6cm}
\end{array}$$
\caption{The PML-II $e^2_\textrm{pml}$-error  as a function of time under $\sigma=\sigma_P$: for  $\delta=1/2$ and varying strength $\sigma_0$ (left),
 and for $\sigma_0=6$ and varying thickness $\delta$ (right).}
\label{fig:testNKG2}
\end{figure}

  Firstly, we take the polynomial choice (\ref{sigma}) as the absorption function, i.e.  $\sigma(x)=\sigma_P(x)$ for PML-II.
Similarly as before, with $\sigma_0=8$ and $\delta=0.5$ fixed for the absorption function, we first show the temporal and spatial accuracy of the numerical method, i.e. FD-FP (\ref{FD-FP}) at $t=4$ in Figure \ref{fig:convergencePML-II}. Next, by taking in FD-FP a very fine mesh ($\tau= 10^{-4}$ and $h=1/256$) so that the numerical discretization error is negligible, the corresponding PML error (\ref{err}) of the PML-II formulation (\ref{KG model pml2}) is shown in Figure \ref{fig:testNKG2}  as a function of time for different values of $\sigma_0$ and $\delta$ (for $R=1$ in (\ref{KG model pml2})).

\begin{figure}[hbt!]
$$\begin{array}{c}
\psfig{figure=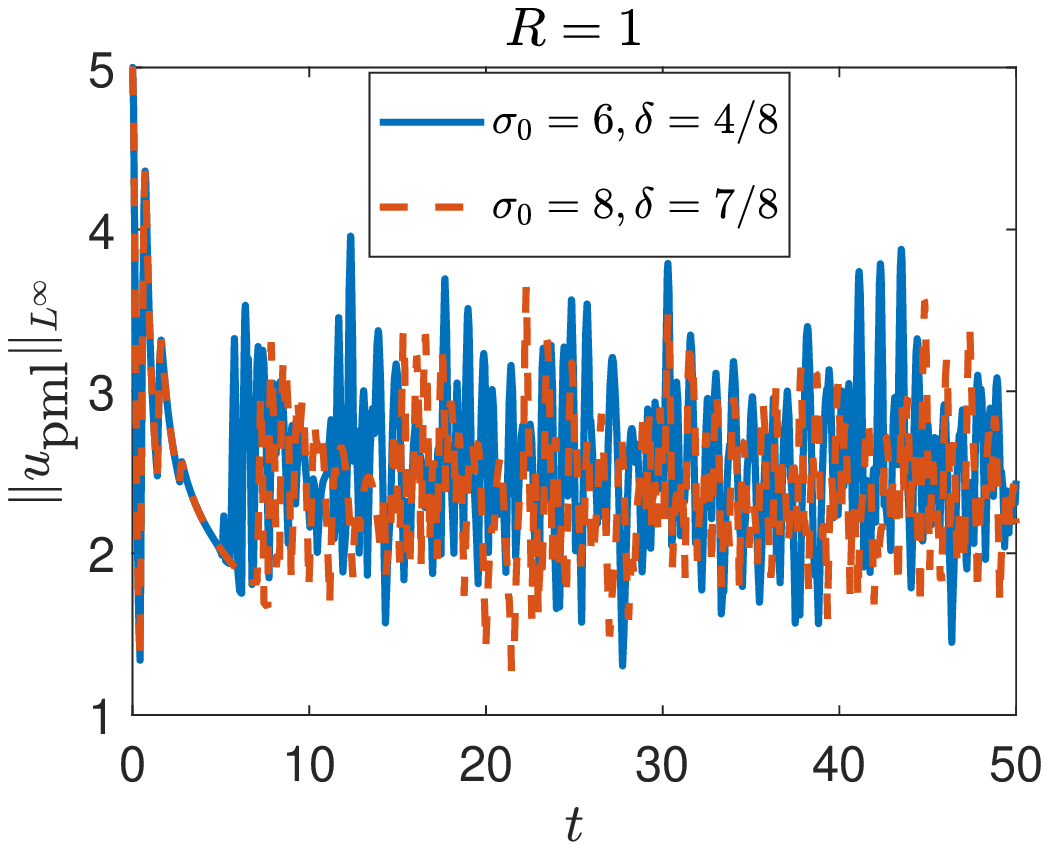,height=4cm,width=5.5cm}
\psfig{figure=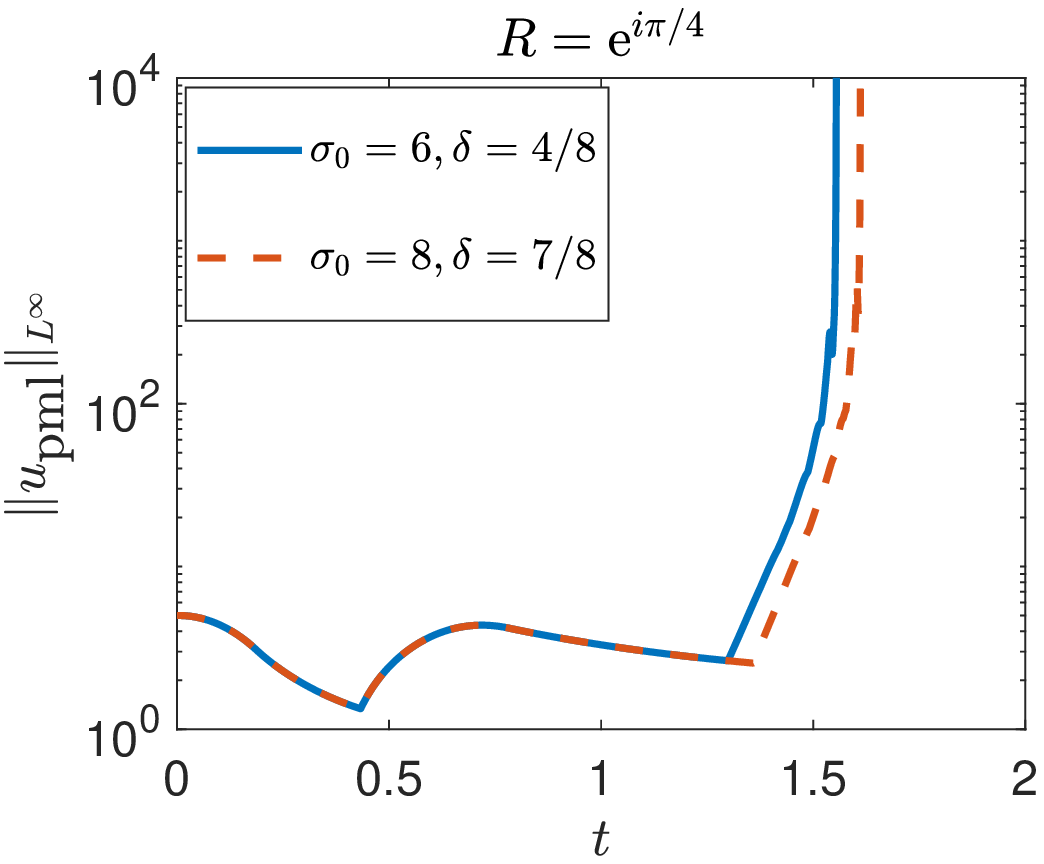,height=4cm,width=5.5cm}\\
\psfig{figure=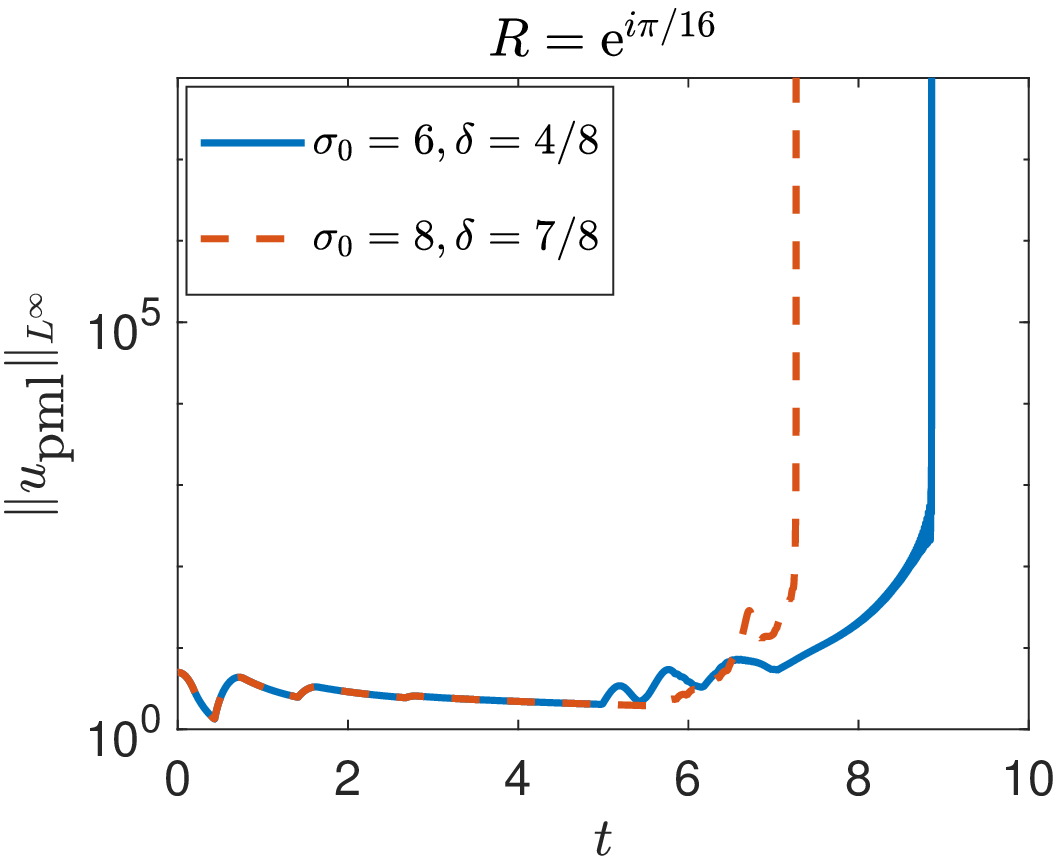,height=4cm,width=5.5cm}
\end{array}$$
\caption{Stability of PML-II: $\|u_{\textrm{pml}}\|_{L^\infty}$ as a function of time for $R=1$ (top left), $R=\fe^{i\pi/4}$ (top right), $R=\fe^{i\pi/16}$ (bottom).}
\label{fig:PML2 stability}
\end{figure}

As  explained in Section \ref{sec2 PML2} by the dispersion relation in the linear case, the factor $R$ must be positive to avoid any instability.
To confirm the formal analysis  and illustrate the stability issue in the PML-II equation (\ref{KG model pml2}) with respect to $R$, we show in Figure \ref{fig:PML2 stability}  the evolution of the maximum norm of the PML solution $u_{\textrm{pml}}$ of (\ref{KG model pml2}) for $R>0$ or $R\in\bC$.

From the numerical results in Figures \ref{fig:convergencePML-II}-\ref{fig:PML2 stability}, we can conclude that:
\begin{itemize}
\item[1)] The PML-II formulation (\ref{KG model pml2}) with $R>0$ is stable (see Figure \ref{fig:PML2 stability}), while, for $R$ with any non-zero imaginary part, instability will occur in (\ref{KG model pml2}).

\item[2)] The FD-FP method (\ref{FD-FP}) is also second-order accurate in time and nearly spectrally accurate in space. Under the same mesh size and time step, it is slightly  less accurate than the EWI-FP method (cf. Figures \ref{fig:convergencePML-II} and  \ref{fig:convergencePML-I}). Therefore, from the numerical discretization  point of view, PML-II is a little bit less efficient than the PML-I.

\item[3)] PML-II effectively  approximates the exact solution of NKGE (\ref{KG model 1d}) (see Figure \ref{fig:testNKG2}). Moreover, under the same absorption function $\sigma(x)$ with the same strength and thickness, the error from the PML-II is smaller than PML-I particularly as time increases (cf. Figures \ref{fig:testNKG2} and  \ref{fig:testNKG}). Therefore, from the modelling truncation point of view, PML-II is more effective than PML-I. We remark that the results of using the usual low order polynomial absorption function for PML-II is similar.
\end{itemize}

\begin{figure}[hbt!]
$$\begin{array}{c}
\psfig{figure=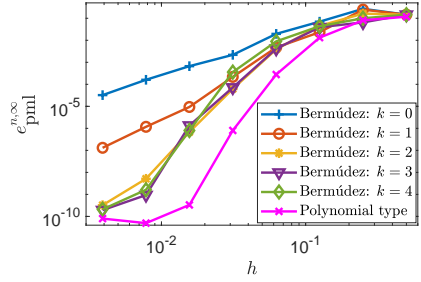,height=4cm,width=6cm}
\end{array}$$
\caption{Spatial $e^{n,\infty}_{\textrm{pml}}$-error of FD-FP for PML-II with $\sigma=\sigma_{B_k}$, for  different values of $k$.}
\label{fig:PML2 Bermudez spatial}
\end{figure}

\begin{figure}[hbt!]
$$\begin{array}{cc}
\psfig{figure=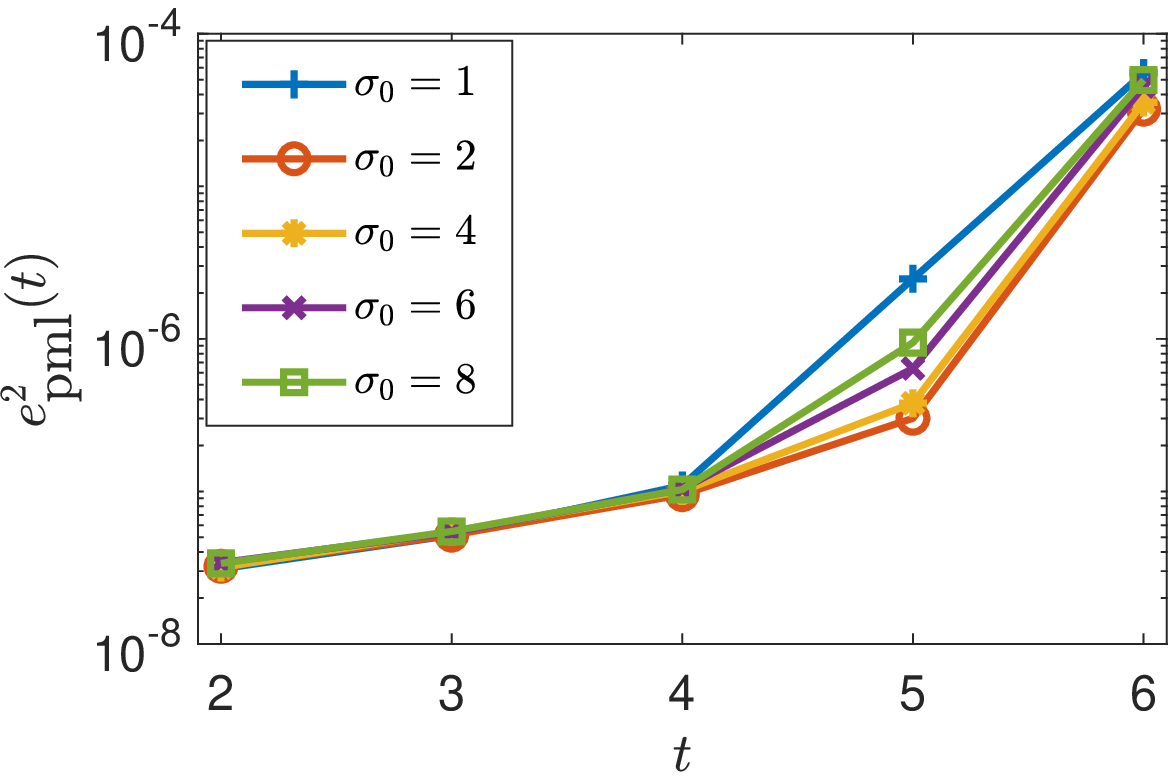,height=4cm,width=6cm}&
\psfig{figure=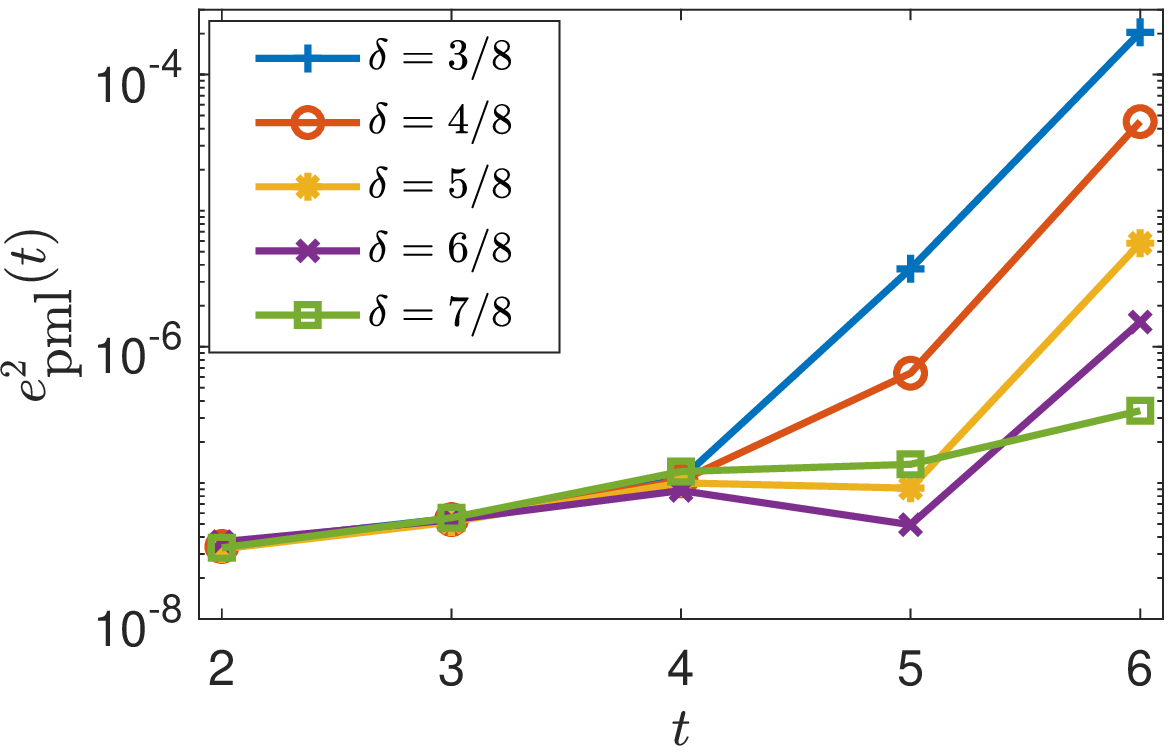,height=4cm,width=6cm}
\end{array}$$
\caption{The PML-II $e^2_\textrm{pml}$-error as a function of time for $\sigma=\sigma_{B_2}$: for  $\delta=1/2$ and varying strength $\sigma_0$ (left), 
 and for  $\sigma_0=6$ and varying thickness $\delta$ (right).}
\label{fig:Bermudez PML2 error}
\end{figure}

  Next, we choose the Berm\'{u}dez type function (\ref{Bermudez}) as the absorption function, i.e.  $\sigma(x)=\sigma_{B_k}(x)$ for   PML-II (\ref{KG model pml2}).
We use the same example (\ref{example}) with $\sigma_0=8,\delta=0.5$ in  (\ref{Bermudez}), and we first test  the numerical discretization error of the  FD-FP method (\ref{FD-FP}) for PML-II (\ref{KG model pml2}) with $R=1$. The spatial error at $t=4$ is presented in Figure \ref{fig:PML2 Bermudez spatial} under different smoothing orders $k\geq0$ in the Berm\'{u}dez function (\ref{Bermudez}). The corresponding temporal error is totally the same as in the previous example with the polynomial choice (\ref{sigma}), so it is omitted here for brevity. On the other hand, to illustrate the efficiency brought by the preconditioner (\ref{preconditioner}) for solving (\ref{FD-FP gmres}), we fix in this example $\sigma=\sigma_{B_2},\ \tau=0.02$ and consider the GMRES solver for the time step $n=1$ in  (\ref{FD-FP gmres}). We show the number of iterations needed by the GMRES solver (without restart) to converge at a threshold $\epsilon>0$ with or without using the preconditioner (\ref{preconditioner}) under different values of $h$ in Table \ref{tab:gmres}. This confirms that the preconditioner leads to a convergence independent of the mesh refinement, which is well-adapted to the pseudo-spectral method. The corresponding PML error (\ref{err}) with respect to time in this example is presented in Figure  \ref{fig:Bermudez PML2 error} for different values of the parameters $\sigma_0$ and $\delta$.

{\footnotesize
\begin{table}[hbt!]
\tabcolsep 0pt \caption{Number of iterations needed by  GMRES under threshold $\epsilon$ for (\ref{FD-FP gmres}) at $n=1$: with (pre) or without (non) the preconditioner. \label{tab:gmres}}
\begin{center}\vspace{-1em}
	\def\temptablewidth{1\textwidth}
{\rule{\temptablewidth}{1pt}}
\begin{tabularx}{\temptablewidth}{@{\extracolsep{\fill}}
p{1.50cm}llllllll}
 $\epsilon=10^{-10}$  & $h=1/128$ &$h=1/256$ &$h=1/512$ &\vline & $\epsilon=10^{-13}$ & $h=1/128$ &$h=1/256$ &$h=1/512$ \\
  \hline
pre &2&2&2&\vline &pre &8 &7&7\\
non  &11&20&35&\vline &non &24 &46&92
\end{tabularx}
{\rule{\temptablewidth}{1pt}}
\end{center}
\end{table}
}

Based on the numerical results in Table \ref{tab:gmres}, Figures \ref{fig:PML2 Bermudez spatial} and \ref{fig:Bermudez PML2 error}, we can draw the following observations:
\begin{itemize}
\item[1)] For a high enough regularization parameter $k\geq0$, the spatial error of FD-FP method converges fast. For $k\geq2$, the Berm\'udez function $\sigma_{B_k}$ is able to offer a near spectral accuracy for the Fourier pseudo-spectral method, although the error is a little bit larger than that for the polynomial choice (\ref{sigma}) (see Figure \ref{fig:PML2 Bermudez spatial}). Practically, taking $k=2$ or $3$ is enough, since  increasing $k$ does not give any more improvements.

\item[2)]  The proposed  GMRES solver (\ref{FD-FP gmres}) for FD-FP method works very well with the help of the preconditioner (\ref{preconditioner}). The number of iterations needed to reach a threshold has been significantly reduced to $O(1)$ for all mesh sizes.  This makes the practical efficiency of the FD-FP method comparable to the EWI-FP method.

\item[3)] The PML error of PML-II with Berm\'udez function  (\ref{Bermudez}) is much smaller than that of the polynomial choice (\ref{sigma}) under the same parameters $\sigma_0$ and $\delta$ (cf. Figures \ref{fig:Bermudez PML2 error} and \ref{fig:testNKG2}). Moreover, the accuracy of the PML-II with Berm\'udez function is much less sensitive to the choice of the parameters $\sigma_0$ and $\delta$ (see Figure  \ref{fig:Bermudez PML2 error}). Therefore, one does not need to tune the parameters for the layer  in a practical computation, which makes it suitable for concrete applications and to later investigate the non-relativistic regime of the NKGE.
\end{itemize}

\medskip
\noindent\textbf{Energy decay and comparison.}
To end, we test and compare the energy behaviour of the two PML formulations.   We define the following energy functional for some function $w=w(x,t)$
\begin{equation*}
    H_I(t;w):=\int_{I}\left[|\partial_tw(x,t)|^2+|\partial_xw(x,t)|^2
    +|w(x,t)|^2+\frac{\lambda}{2}|w(x,t)|^4\right]dx,\quad t\geq0.
\end{equation*}
With the exact solution $w=u$ of (\ref{KG model 1d}), the above $H_I(t;u)$ denotes the part of the total energy (\ref{energy def}) inside the physical domain $I=(-L,L)$ for the NKGE (\ref{KG model 1d}).
Note that the quantity $H_I(t;u)$ is not conserved by (\ref{KG model 1d}), and in fact it should be decaying with respect to time, since the waves in the solution $u(x,t)$ keep propagating to the far field as time evolves.

\begin{figure}[hbt!]
$$\begin{array}{c}
\psfig{figure=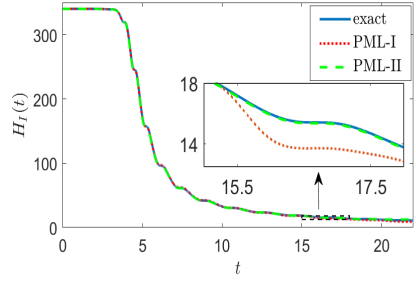,height=4cm,width=6cm}
\end{array}$$
\caption{Energy decay of the PMLs in physical domain: $H_I(t)$ obtained from PML-I, PML-II and the NKGE.}
\label{fig:energy}
\end{figure}

We take the example (\ref{example}) with $I=(-4,4)$. For PML-I given by (\ref{KG model pml}), we consider $\sigma=\sigma_P$ with $\sigma_0=8,\,\delta=6/8,\,\alpha=0$. For PML-II as defined by (\ref{KG model pml2}), we take $\sigma=\sigma_{B_2}$ with $\sigma_0=3,\,\delta=6/8,\,R=1$. We compute $H_I$ with the exact solution $w=u$ of (\ref{KG model 1d}), and with the exact solution $w=u_{\textrm{pml}}$ of PML-I (\ref{KG model pml}) or  PML-II (\ref{KG model pml2}). The ``exact'' solutions here are obtained numerically as before. The numerical results are reported in Figure \ref{fig:energy} until the final time $t=22$. We see that both PML-I and PML-II follow the exact energy decay very well in long time, and the approximation from PML-II is however more accurate than PML-I for all times.

Overall, based on the presented numerical experiments in this section  for approximating the NKGE (\ref{KG model 1d}) in the classical scaling, we suggest to rather use the PML-II formulation (\ref{KG model pml2}) with the regularized Berm\'udez absorption function (\ref{Bermudez})
for $k\geq 2$, not only from the point of view of the temporal and spatial accuracy, 
but also for the stability to fix the tuning parameters of the PML.


\section{PML for non-relativistic scaling}\label{SectionNR}
Physically when $c\to\infty$, the NKGE (\ref{KG model}) is said to be in the non-relativistic limit regime \cite{Nonrelat1,Nonrelat2,RKG-PRD}. In such a case, it is often convenient to introduce the scaling \cite{Dong,RKG}: 
$$x\to\frac{x}{x_s},\quad t\to\frac{t}{t_s},\quad \lambda\to\lambda mx_s^2,\quad \mbox{with}\quad
t_s=\frac{mx_s^2}{\hbar},\quad \eps=\frac{\hbar}{mcx_s},
$$
where $t_s$ and $x_s$ are the time and   length units, respectively, and $\eps>0$ is a dimensionless parameter. Then, the formulation of
the one-dimensional real-valued case of (\ref{KG model}) in the non-relativistic scaling  reads \cite{BCZ,BaoZhao,Mehats}:
\begin{equation}\label{KG model eps 1d}
\left\{\begin{split}
&\eps^2\partial_{tt} u(x,t)-\partial_{xx}u(x,t)+\frac{1}{\eps^2}u(x,t)+\lambda u(x,t)^3=0,\quad t>0,\ x\in\bR,\\
&u(x,0)=u_0(x),\quad \partial_tu(x,0)=v_0(x),\quad x\in\bR.
\end{split}\right.
\end{equation}
In the following, we shall consider $\eps\in(0,1]$. The NKGE (\ref{KG model eps 1d}) in such scaling could describe the physical system  where the wave speed is smaller than the speed of light but with non-negligible relativistic effects \cite{BaoZhao,RKG,RKG-PRD}. It also appears in the high-plasma-frequency limit regime for plasma physics \cite{BDZ,BaoZhao-KGZ}. For $\eps=1,$ (\ref{KG model eps 1d}) leads to the classical NKGE (\ref{KG model 1d}), and as $\eps\to0$, we are in the non-relativistic limit regime.

\subsection{PML-II} We adopt the same notations as introduced in the previous section.
Now in the same spirit as (\ref{KG model pml2}), the PML-II in the non-relativistic scaling reads:
\begin{equation}\label{KG limit pml2}
\left\{\begin{split}
&\eps^2\partial_{tt} u-\frac{1}{1+R\sigma}\partial_x\left(\frac{1}{1+R\sigma}\partial_xu\right) +\frac{1}{\eps^2}u+\lambda u^3=0,\quad t>0,\ x\in I^*,\\
&u(x,0)=u_0(x),\quad \partial_tu(x,0)=v_0(x), \quad x\in I^*,\\
&u(-L^*,t)=u(L^*,t),\quad t\geq0. 
\end{split}\right.
\end{equation}

For the linear case of (\ref{KG limit pml2}), i.e. $\lambda=0$, let us now consider a plane wave solution: $u(x,t)=\fe^{i(k x-\omega t)}$ with $k,\omega\in\bR$ satisfying the dispersion relation:
$$\omega=\pm\frac{\sqrt{1+\eps^2k^2}}{\eps^2},$$
which provides the phase velocity
$$v_p=\frac{\omega}{k}=\pm \frac{\sqrt{k^{-2}+\eps^2}}{\eps^2}.$$
From the above, it is clear to see that waves with all wavelengths are travelling to infinity at the speed $O(\eps^{-2})$ as $\eps\to0$. Therefore, when $\eps$ is small, the waves in (\ref{KG model eps 1d}) enter and pass through the damping layer $L\leq |x|\leq L+\delta$ very quickly, and then get reflected at the opposite outer boundary $|x|=L+\delta$ because of the periodic boundary conditions. This fact may leave the layer not enough time to effectively absorb all the waves. In fact, we will show later by numerical tests, for a fixed PML setup under the classical polynomial absorption function (\ref{sigma}), i.e. fixed $\sigma,\,\delta,\,R$ as $\eps$ decreases, that the  PML-II (\ref{KG limit pml2}) will soon fail to approximate the exact solution of NKGE in the physical domain.

In order to get a stable absorber with rather uniform damping effect for $\eps\in(0,1]$, we will verify later by numerical experiments that  we need to choose
\begin{equation}\label{Reps}
R=R^\eps=O(1/\eps^2)\in\bR^{+},\quad \mbox{if}\quad \sigma(x)=\sigma_P(x),
\end{equation}
for the PML-II (\ref{KG limit pml2}) with fixed $\sigma_P(x)$. This is clearly equivalent to say that with $R>0$ fixed as $\eps\to0$, we need to choose the strength of the absorption (\ref{sigma})  $\sigma_0=O(\eps^{-2})$. Based on the studies from the previous section, we know that the accuracy of the Berm\'{u}dez absorption function (\ref{Bermudez}) is not very sensitive to the choice of the strength parameter $\sigma_0$. 
Therefore, it is hopeful and reasonable that with the Berm\'{u}dez's absorption function in the PML-II (\ref{KG limit pml2}), the damping effect is rather uniform  for $\eps\in(0,1]$ by using
$$R=O(1)\in\bR^{+},\quad \mbox{if}\quad \sigma(x)=\sigma_{B_k}(x).$$

\begin{figure}[hbt!]
$$\begin{array}{c}
\psfig{figure=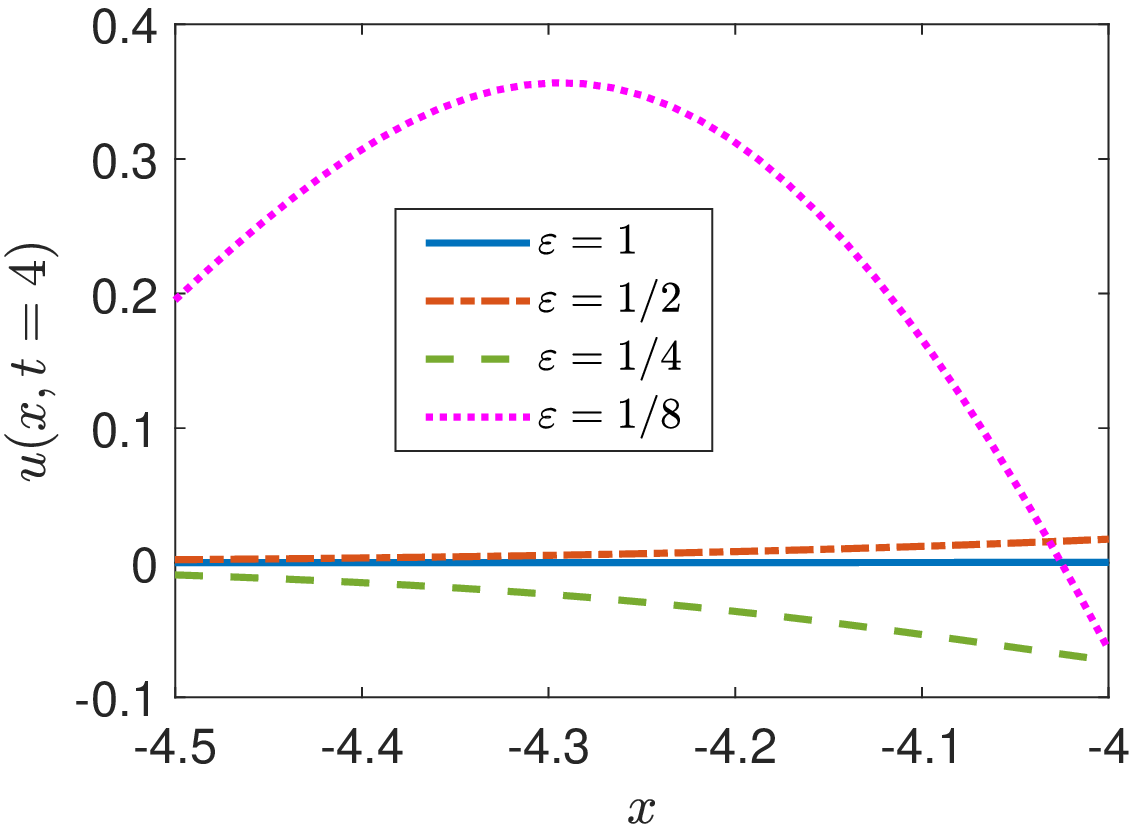,height=4cm,width=6cm}
\psfig{figure=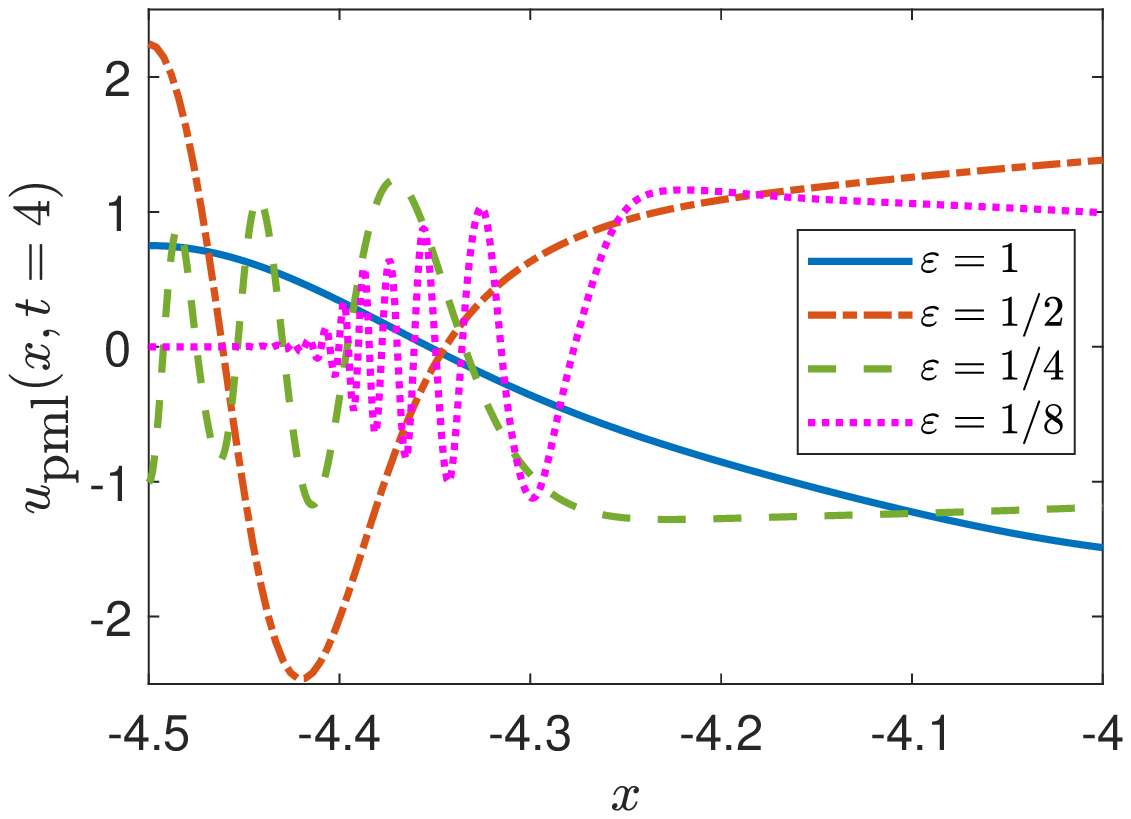,height=4cm,width=6cm}\\
\psfig{figure=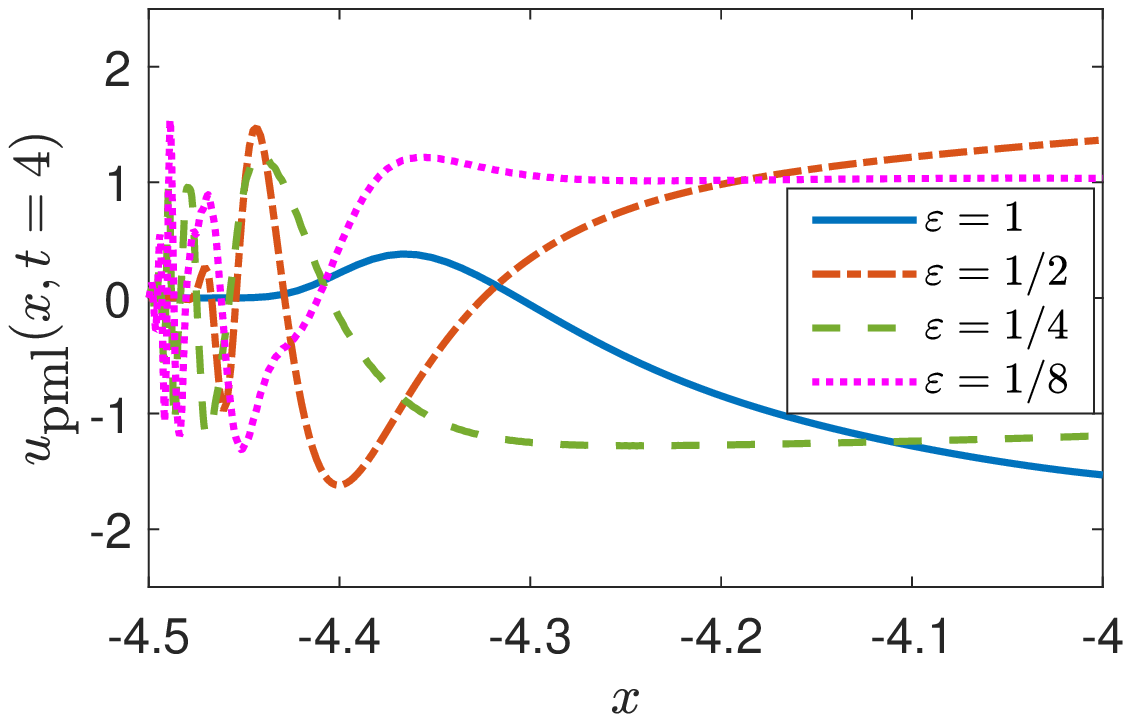,height=4cm,width=6cm}
\end{array}$$
\caption{Profiles of the exact solution $u(x,t)$ of NKGE  (top left) and the PML-II solution $u_{\textrm{pml}}(x,t)$ with $R=\eps^{-2}$ and $\sigma=\sigma_P$ (top right), or with $R=1$ and $\sigma=\sigma_{B_2}$ (bottom).}
\label{fig:pml2insidelayer}
\end{figure}

It is well-known that as $\eps\to0$, the solution of the NKGE (\ref{KG model eps 1d}) contains rapid oscillations in time \cite{BaoZhao,Nonrelat1,Nonrelat2,KG-SZ}:
\begin{equation}\label{mfe}
u(x,t)=\fe^{it/\eps^2}z(x,t)+\fe^{-it/\eps^2}\overline{z}(x,t)+O(\eps^2),\quad t\geq0,\ x\in\bR,
\end{equation}
where $z=z(x,t)$ solves a nonlinear Schr\"{o}dinger equation independent of $\eps$. In space, the solution $u$ of NKGE (\ref{KG model eps 1d}) is rather smooth with wavelength independent of $\eps$. As an illustrative example, we take
\begin{equation}\label{example eps}
u_0=5\fe^{-x^2},\quad v_0=\frac{1}{2}\textrm{sech}(x^2),\quad \lambda=\frac{1}{2},
\end{equation}
and we show the profile of the exact solution $u(x,t)$ at $t=4$ in the physical domain $x\in(-4.5,-4)$ under different $\eps$ in Figure \ref{fig:pml2insidelayer}.
While in the PML-II formulation (\ref{KG limit pml2}), due to the fast reflected waves (or the incoming waves due to the periodic boundary condition), the $\eps$-dependent oscillations will be induced to the space variable $x$ of the solution. To see this, we consider in  (\ref{KG limit pml2}) the physical domain $I=(-4,4)$ with the layer setup $\sigma_0=3,\,\delta=0.5,\,R=\eps^{-2}$ in $\sigma_P$ (\ref{sigma}), and we show in Figure \ref{fig:pml2insidelayer} the PML solution $u_{\textrm{pml}}(x,t)$ of (\ref{KG limit pml2}) at $t=4$ inside the layer $-L-\delta\leq x\leq -L$. The PML solution $u_{\textrm{pml}}(x,t)$ under the same parameters values but with $R=1$ and $\sigma_{B_2}$ (\ref{Bermudez}) is also given in  Figure \ref{fig:pml2insidelayer}. It can be seen from Figure \ref{fig:pml2insidelayer} that, as $\eps$ decreases, the PML solution $u_{\textrm{pml}}(x,t)$ exhibits more  and more spatial oscillations inside the layer. Such phenomenon will cause  numerical difficulties to the spatial discretization of the PML-II formulation (\ref{KG limit pml2}). As we shall report in the numerical tests, the accuracy of the Fourier pseudo-spectral discretization for (\ref{KG limit pml2})  depends on $\eps.$ 

For the FD-FP method, the finite-difference discretization of (\ref{KG limit pml2}) with $R=R^\eps$ 
directly reads:
$$ \eps^2\frac{u^{n+1}-2u^n+u^{n-1}}{\tau^2}+\frac{A^\eps}{2}
  (u^{n+1}+u^{n-1})+\frac{1}{2\eps^2}
  (u^{n+1}+u^{n-1})+\lambda (u^n)^3=0,$$
which gives the FD-FP scheme
\begin{align}\label{FD-FP eps}
  u^{n+1}=-u^{n-1}
  +(G^\eps)^{-1}\left[\frac{2\eps^2}{\tau^2}u^n
  -\lambda (u^n)^3\right],\quad n\geq1,
\end{align}
with
\begin{align*}
 G^\eps=\left[\left(\frac{\eps^2}{\tau^{2}}+\frac{1}{2\eps^2}\right)\mathbb{I}
  +\frac{1}{2}A^\eps\right],\quad A^\eps=-d_0D_1d_0D_1,\quad d_0=\mathrm{diag}(1/(1+R^\eps\sigma)).
\end{align*}
For the above matrix $G^\eps$ that depends on $\eps$ in the scheme,
we will show later by numerical results that its condition number remains small and stays uniformly bounded as $\eps\to0$. Similarly as (\ref{FD-FP c}), it can also be implemented efficiently by the GMRES solver with a preconditioner in the non-relativistic scaling:
\begin{equation}\label{preconditioner eps}
\left\{
  \begin{array}{l}
  \displaystyle
  u^{n+1}=-u^{n-1}
  +w^n,\quad n\geq1,\\
   \displaystyle \mathcal{P}^\eps G^\eps w^n   =\mathcal{P}^\eps\left[\frac{2\eps^2}{\tau^2}u^n -\lambda (u^n)^3
   \right],
  \end{array}
  \right.\quad \mathcal{P}^\eps=\left(\frac{\eps^2}{\tau^{2}}
  +\frac{1}{2\eps^2}-\frac{\partial_{xx}}{2}\right)^{-1}.
  \end{equation}
The starting value for (\ref{FD-FP eps}) if  directly using the construction from the Taylor expansion reads
$$u^1=u_0+\tau v_0-\frac{\tau^2}{2\eps^2}\left[A^\eps u_0+\frac{u_0}{\eps^2}+\lambda u_0^3\right],$$
where the $O(1/\eps^2)$ and $O(1/\eps^4)$ terms in the above will induce instability to the scheme, since $u^1=O(1/\eps^4)$ could be very large for small   $\eps$. To avoid this, we consider a filtered data:
\begin{equation}\label{u1 eps}u^1=u_0+\tau v_0-\frac{\tau}{2}\sin(\tau/\eps^2)\left[A^\eps u_0+\lambda u_0^3\right]-
\frac{\tau}{2}\sin(\tau/\eps^4)u_0,\end{equation} so that
$u^1=O(1)$.

Thanks to the filtered starting value (\ref{u1 eps}) as well as the time averaging, the FD-FP method (\ref{FD-FP eps}) is stable for all $\eps\in(0,1]$. This means that we have the boundedness of the numerical solution of FD-FP up to a  fixed time if $\tau\leq C$ for some $C>0$ independent of $\eps$ or the spatial mesh size $h$.  However,  owning to the $\eps$-dependent oscillation frequencies in both time and space directions of the PML solution of (\ref{KG limit pml2}), it is  expected that the FD-FP method (\ref{FD-FP}) is not uniformly accurate for approximating (\ref{KG limit pml2}) as $\eps\to0$. One then needs to use mesh sizes $\tau>0$ and $h>0$ that depend on $\eps$ to reach the same accuracy level at a fixed time.


\begin{remark}
It is beyond the scope of this paper to investigate whether the asymptotic expansion (\ref{mfe}) also holds for the PML-II equation (\ref{KG limit pml2}) or not. If it is still valid, then formally as $\eps\to0$, the PML-II (\ref{KG limit pml2}) is consistent with the limit: a nonlinear Schr\"odinger equation with PML presented in \cite{Antoine-review,PML-NLS}. We shall address this question in a future study.
\end{remark}

\subsection{PML-I}
Next, we consider the PML-I formulation under the non-relativistic scaling. To write the correct $\eps$-dependent PML  model for (\ref{KG model eps 1d}), we shall proceed the derivation as in \cite{PML-wave}.

We begin by considering a linear case of (\ref{KG model eps 1d}):
\begin{equation}\label{linear model}\eps^2\partial_{tt}u-\partial_{xx}u+\frac{1}{\eps^2}u+\lambda u=0.\end{equation}
By  taking the Laplace transform on both sides of the above equation, we get
$$\left(\eps^2s^2-\partial_{xx}+\frac{1}{\eps^2}+\lambda\right)
\widehat{u}(x,s)=0,\quad \widehat{u}(x,s)=\int_0^\infty u(x,t)
\fe^{-st}dt,$$
where $s$ is the Laplace variable.
By (\ref{ansatz}), we find formally, for some $\alpha\geq0$ and absorption function $\sigma(x)$,
\begin{equation}\label{pml1 eps derivation eq1}\left(\eps^2s^2+\frac{1}{\eps^2}+\lambda\right)
\widehat{u}=\frac{s+\alpha}{s+\alpha+\sigma}\partial_{x}
\left(1-\frac{\sigma}{s+\alpha+\sigma}\right)
\partial_{x}\widehat{u}.\end{equation}
Let $$\widehat{\eta}_1=-\frac{1}{s+\alpha+\sigma}
\partial_{x}\widehat{u},$$
then from (\ref{pml1 eps derivation eq1}) we have
\begin{align*}
    \eps^2\left(s^2+s\sigma\right)\widehat{u}
    -\eps^2\sigma\alpha\widehat{u}
    +\sigma\frac{\eps^2\alpha^2}{s+\alpha}\widehat{u}
    =\partial_{x}\left(\partial_{x}\widehat{u}+\sigma
    \widehat{\eta}_1\right)-\left(\frac{1}{\eps^2}
    +\lambda\right)\widehat{u}-\sigma\frac{\eps^{-2}
    +\lambda}{s+\alpha}\widehat{u}.
\end{align*}
By taking the inverse Laplace transform of the above equation and further letting
$$\widehat{\eta}_2=-\frac{\eps^2\alpha^2+\eps^{-2}+\lambda}
{s+\alpha}\widehat{u},$$ we obtain
\begin{align*}
&\eps^2\partial_{tt}u-\partial_{xx}u+\frac{1}{\eps^2}u
    +\lambda u=\partial_x(\sigma \eta_1)+\eps^2\sigma\alpha u -\eps^2\sigma\partial_tu+\sigma \eta_2,\\
   & \partial_t\eta_1+(\alpha+\sigma)\eta_1+\partial_x u=0,\quad \partial_t\eta_2+\alpha\eta_2+\left(\eps^2\alpha^2
+\frac{1}{\eps^2}\right)u+\lambda u=0,
\end{align*}
which is a PML system for the linear model problem (\ref{linear model}).
To proceed from the linear case to the nonlinear case, as proposed in \cite{PML-wave}, we directly replace the term $\lambda u$ by $\lambda u^3$ in the above expressions.
Then together with the domain truncation similarly as (\ref{KG model pml}), the PML-I formulation for NKGE (\ref{KG model eps 1d}) in the non-relativistic scaling reads
\begin{equation}\label{PML-I eps}\left\{\begin{split}
&\eps^2\partial_{tt}u-\partial_{xx}u+\frac{1}{\eps^2}u
    +\lambda u^3=\sigma\left(\eta_2+\eps^2\alpha u -\eps^2\partial_tu\right)+\partial_x(\sigma \eta_1),\quad t>0,\ x\in I^*,\\
&\partial_t\eta_1+(\alpha+\sigma)\eta_1+\partial_x u=0,\\
&\partial_t\eta_2+\alpha\eta_2+\left(\eps^2\alpha^2
+\frac{1}{\eps^2}\right)u+\lambda u^3=0,\\
&u(x,0)=u_0(x),\quad \partial_tu(x,0)=v_0(x),\quad\eta_1(x,0)=\eta_2(x,0)=0, \quad x\in I^*,\\
&u(-L^*,t)=u(L^*,t),\quad \eta_1(-L^*,t)=\eta_1(L^*,t),\quad \eta_2(-L^*,t)=\eta_2(L^*,t),\quad t\geq0.\end{split}\right.
\end{equation}

Asymptotically, the limit of the PML-I system (\ref{PML-I eps}) as $\eps\to0$ is not clear.
It can be seen from (\ref{PML-I eps}) that there is a $O(\eps^{-2})$ term in the equation of $\eta_2$, which acts as a  very stiff source term as $0<\eps\ll1$. Numerically, to propose a stable and accurate numerical algorithm for solving (\ref{PML-I eps}) is a challenge for small $\eps$. In this paper, we do not intend to address this issue, and we shall focus on the investigation of the two PML formulations and the two choices of absorption functions.  In the following, let us simply apply the  EWI-FP discretization from Section \ref{sec:classical PMLI}.


The EWI-FP method for (\ref{PML-I eps}) then reads: for $n\geq0,$
\begin{subequations}\label{EWI-FP eps}
\begin{align}
 u^{n+1}=&\cos(\langle\partial_x\rangle_\eps\tau)u^n+\frac{\sin(\langle\partial_x\rangle_\eps\tau)}{
 \langle\partial_x\rangle_\eps}v^n+\frac{
 \tau\sin(\langle\partial_x\rangle_\eps\tau)}{2\eps^2\langle\partial_x\rangle_\eps}f^n, \\
 \eta_1^{n+1}=&\fe^{-(\sigma+\alpha)\tau}\eta_1^n-\frac{\tau}{2}\left[\fe^{-(\sigma+\alpha)\tau}
 \partial_x u^{n}+\partial_x u^{n+1}\right],\\
 \eta_2^{n+1}=&\fe^{-\alpha\tau}\eta_2^n-\frac{\tau}{2}\left[\fe^{-\alpha\tau}\left((\eps^2\alpha^2+\eps^{-2})u^n
 +\lambda (u^n)^3\right)+(\eps^2\alpha^2+\eps^{-2})u^{n+1}
 +\lambda (u^{n+1})^3\right],\\
 v^{n+1}=&-\langle\partial_x\rangle_\eps\sin(\langle\partial_x\rangle_\eps\tau)u^n
 +\cos(\langle\partial_x\rangle_\eps\tau)v^n+
 \frac{\tau}{2\eps^2}\left[\cos(\langle\partial_x\rangle_\eps\tau)f^n+f^{n+1}\right],
\end{align}
\end{subequations}
with
$$\langle\partial_x\rangle_\eps=\frac{\sqrt{1-\eps^2
\partial_{xx}}}{\eps^2},\quad
f^n=\sigma\left(
 \eta_2^n-\eps^2v^n+\eps^2\alpha u^n\right)+\partial_x(\sigma \eta_1^n)-\lambda (u^n)^3.
$$

As we discussed in the previous section, for each fixed $\eps$, the EWI-FP has a CFL as $\tau<Ch$. Now, the scheme is much more constrained since
for each fixed $h$ as $\eps\to0$,  the EWI-FP scheme (\ref{EWI-FP eps})  has a severe stability constraint: $\tau\lesssim\eps^2$, due to the stiffness that we mentioned  above. 
Moreover, we are not able to directly apply the Berm\'udez absorption function in (\ref{EWI-FP eps}) due to the singularity in $\sigma_{B_k}$ (\ref{Bermudez}) and stability issue. So, we have to consider the classical choice (\ref{sigma}) for (\ref{PML-I eps}) or (\ref{EWI-FP eps}) and tune the parameters for every $\eps$. Concerning these drawbacks, we recommend to rather use the PML-II formulation (\ref{KG limit pml2}) than  PML-I (\ref{PML-I eps}) for NKGE in the non-relativistic scaling.

\subsection{Numerical results}\label{sec:result}
In this subsection, we present the numerical  results for the above two PML formulations for the NKGE in the non-relativistic scaling.  Let us take  the given data in (\ref{example eps}) for the NKGE (\ref{KG model eps 1d})  in the following tests.
We begin with the numerical investigations of PML-II, and we shall consider similarly as before its performance under the two kinds of absorption functions, and finally analyze shortly PML-I.

\begin{figure}[hbt!]
$$\begin{array}{c}
\psfig{figure=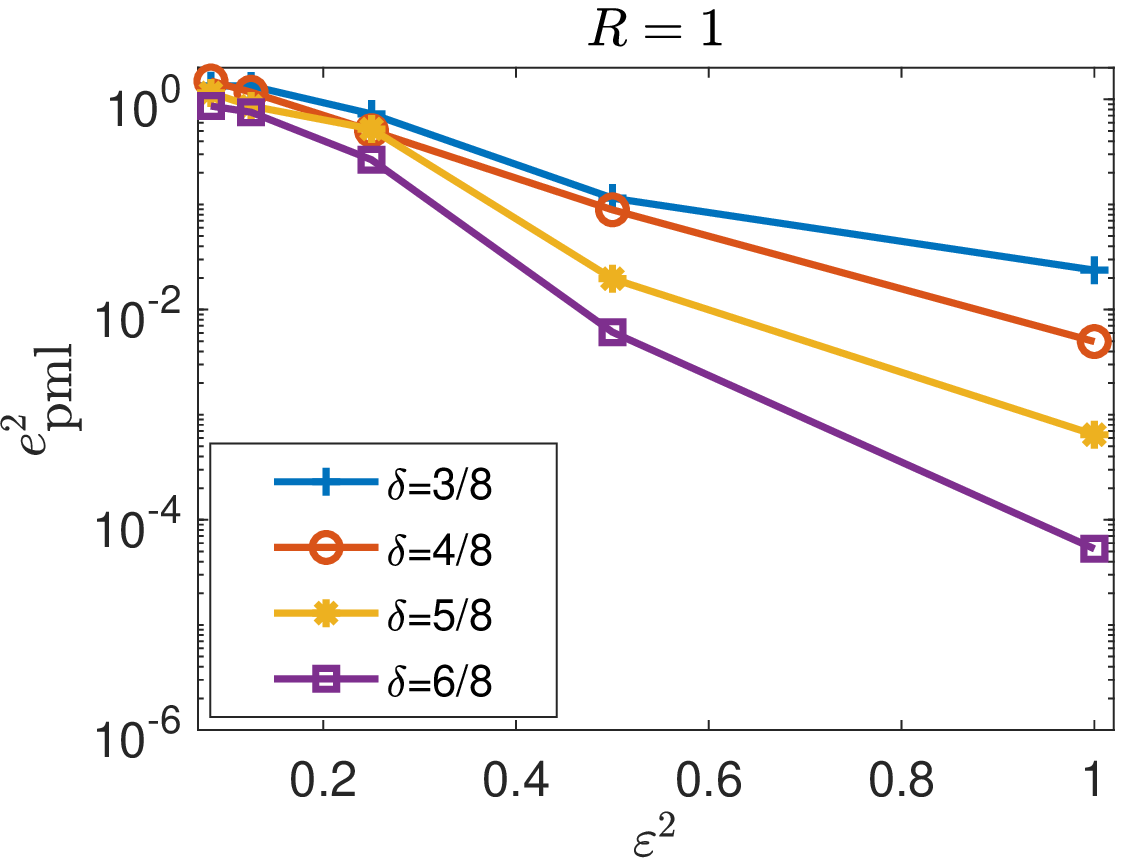,height=4cm,width=6cm}
\psfig{figure=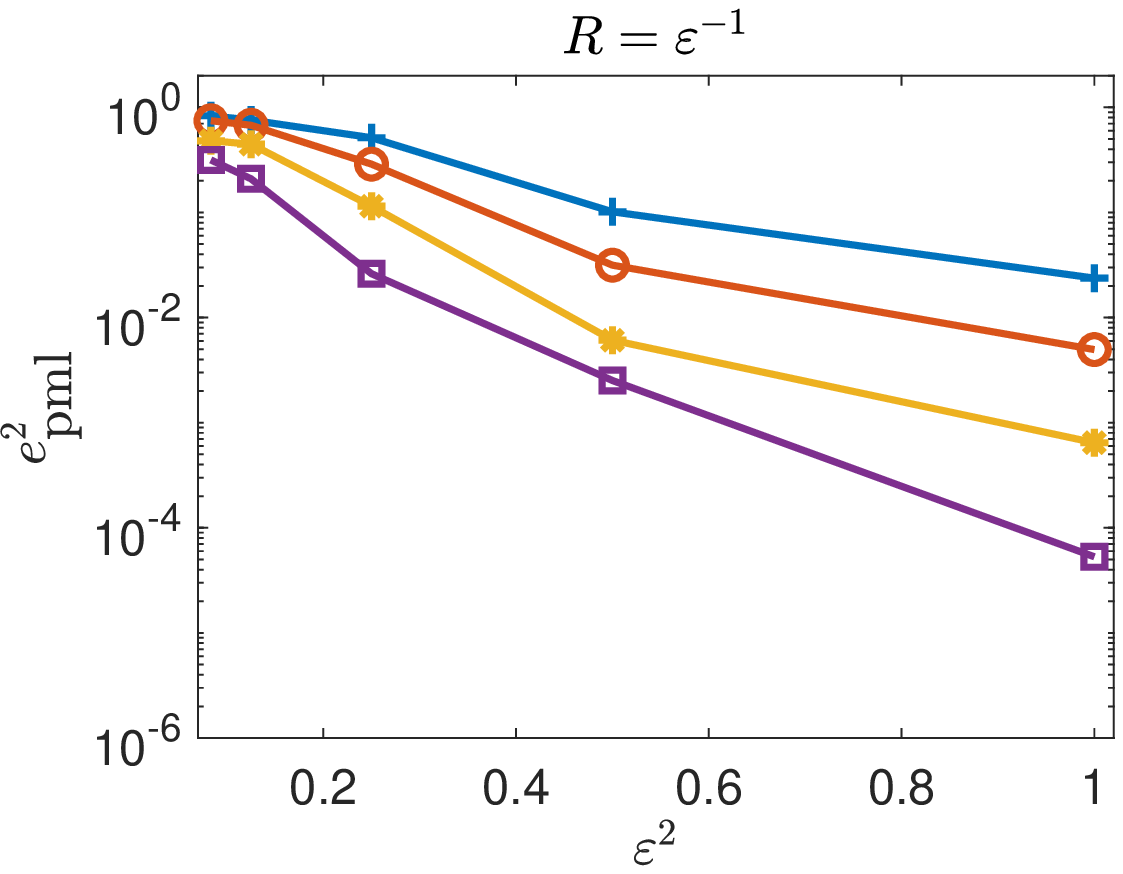,height=4cm,width=6cm}\\
\psfig{figure=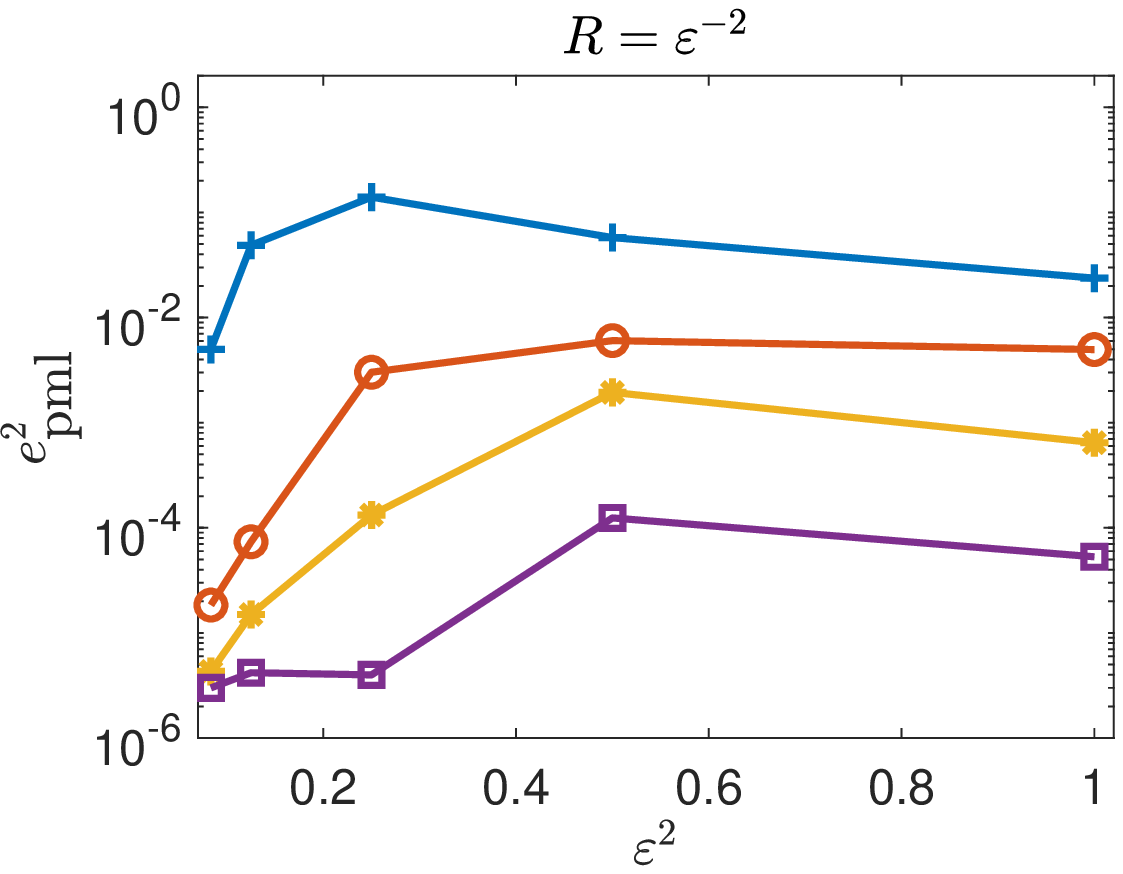,height=4cm,width=6cm}
\end{array}$$
\caption{PML-II $e^2_\textrm{pml}$-error  with respect to $\eps$ for fixed $\sigma_0=3$ and different parameter values  $\delta$: $\sigma=\sigma_P$  with strategy $R=1$ (top left), $R=\eps^{-1}$ (top right) and $R=\eps^{-2}$ (bottom).}
\label{fig:eps}
\end{figure}

\medskip
\noindent\textbf{PML-II with polynomial type absorption function.}
We start with the \textit{polynomial type} absorption function $\sigma(x)=\sigma_P(x)$ as (\ref{sigma}) for the PML-II (\ref{KG limit pml2}). First of all, we illustrate the performance of the PML-II formulation (\ref{KG limit pml2}) with (\ref{sigma}) under different choices of the parameter $R>0$. We measure the error (\ref{err}) on the physical domain $I=(-4,4)$ by choosing
$$R=1,\quad \mbox{or}\quad R=1/\eps,\quad \mbox{or}\quad R=1/\eps^2,$$ with $\sigma_0=3$ fixed and some values of the layer thickness $\delta$. We show in Figure \ref{fig:eps}  the error (\ref{err}) of the PML-II (\ref{KG limit pml2}) for approximating (\ref{KG model eps 1d}) at $t=4$ under  different $\eps$, where the last choice is clearly the most accurate.

\begin{figure}[hbt!]
$$\begin{array}{cc}
\psfig{figure=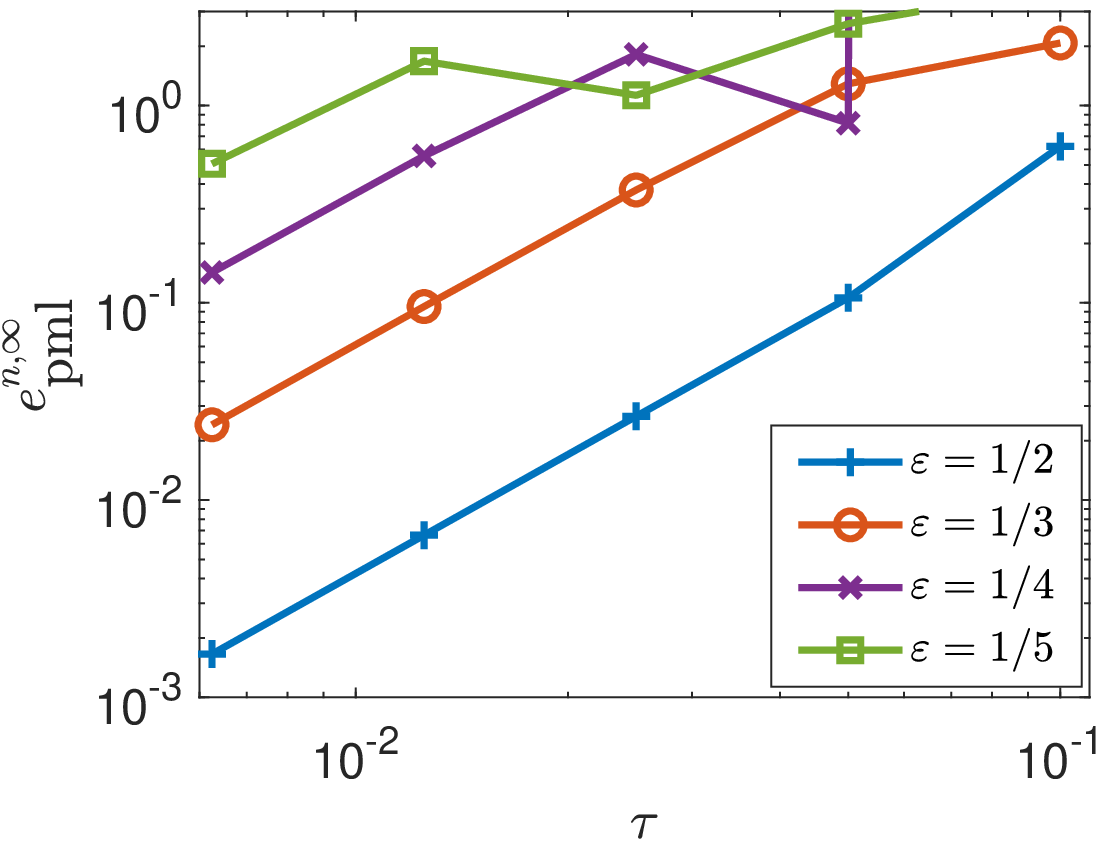,height=4cm,width=6cm}&
\psfig{figure=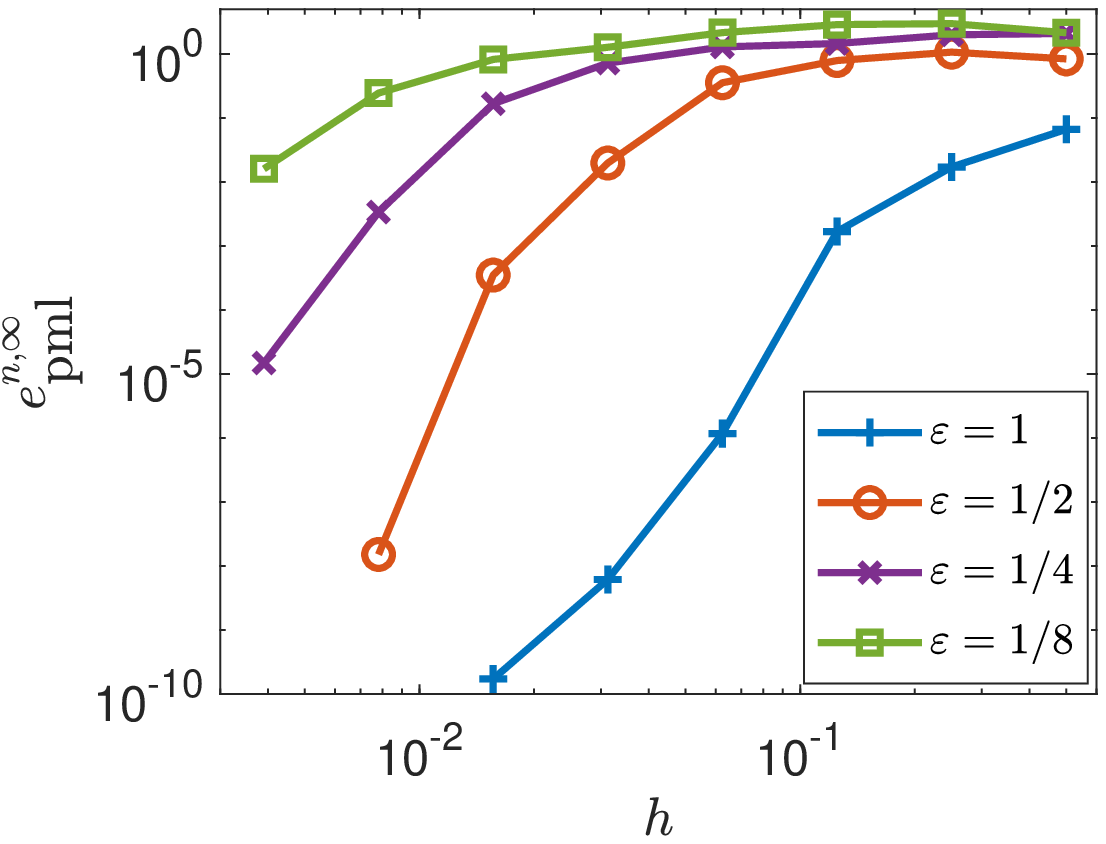,height=4cm,width=6cm}
\end{array}$$
\caption{Temporal (left) and spatial (right) accuracy of FD-FP:  $e^{n,\infty}_{\textrm{pml}}$-error for $\sigma=\sigma_P$  with different $\eps$.}
\label{fig:pml2eps accuracy}
\end{figure}

\begin{figure}[hbt!]
$$\begin{array}{c}
\psfig{figure=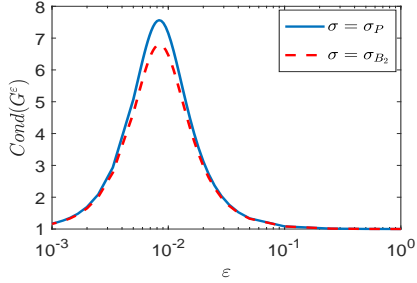,height=4cm,width=6cm}
\end{array}$$
\caption{The condition number of the iterative matrix $G^\eps$ of FD-FP for PML-II: polynomial case $\sigma=\sigma_P$ and
Berm\'udez type $\sigma=\sigma_{B_2}$.}
\label{fig:cond}
\end{figure}

Next, we test the performance of the FD-FP method (\ref{FD-FP eps}) for approximating the PML-II formulation (\ref{KG limit pml2}) for different $\eps$ under the choice $R=R^\eps=1/\eps^2$. In these tests, we fix $\sigma_0=3$ and $\delta=4/8$ for the layer. The reference solution is obtained by using the FD-FP method with a very fine mesh, e.g. $\tau=10^{-4},\,h=1/2048$ for each $\eps$. The temporal  and  spatial errors (\ref{errorinfty}) of the FD-FP method under different $\eps$ at $t=4$ are shown in Figure \ref{fig:pml2eps accuracy}.
Moreover, by taking $\tau=10^{-4},\,h=1/128$ as mesh sizes in this test,  we show in Figure \ref{fig:cond} the condition number (in spectral norm) of the iterative matrix $G^\eps$ as a function of $\eps$ for the FD-FP method (\ref{FD-FP eps}).


Based on the numerical results in Figures \ref{fig:eps}-\ref{fig:cond}, we can conclude on the following observations for the PML-II formulation and the FD-FP method with polynomial type absorption function $\sigma_P$:
\begin{itemize}
\item[1)] With the strategy $R=1/\eps^2$, the accuracy of the PML-II is uniform with respect to $\eps$ as far as we could test. In  contrast, if one takes $R=1$ or $R=1/\eps$, the error of PML-II increases dramatically as $\eps$ decreases (see Figure  \ref{fig:eps}). Overall, when $\eps$ is small, the accuracy of PML-II with $\sigma_P$ is  much more sensitive than for the classical scaling case to the parameters $R$ or saying equivalently the strength $R\sigma_0$  and the layer size $\delta$.

\item[2)] For  $0<\eps\leq1$, the condition number of the iterative matrix $G^\eps$ in the FD-FP scheme (\ref{FD-FP eps}) stays uniformly bounded and small (see Figure \ref{fig:cond}), so the linear system in (\ref{FD-FP eps}) can  be solved accurately and efficiently.  For each fixed $\eps$, we still have the second-order temporal accuracy and the near spectral accuracy in space for the FD-FP method (see Figure \ref{fig:pml2eps accuracy}).  When $\eps$ decreases with fixed mesh size $\tau$ and $h$, the approximation error of FD-FP in both time and space increases quickly. The $\eps$-dependence of the temporal accuracy is due to the $O(\eps^{-2})$ oscillation frequency in time, and the spatial accuracy is affected by the fast reflections at the boundary of the layer as we explained before. A future goal is therefore to construct a numerical algorithm that is uniformly accurate for $0<\eps\leq1$, which would certainly make the PML more efficient.
\end{itemize}

\begin{figure}[hbt!]
$$\begin{array}{c}
\psfig{figure=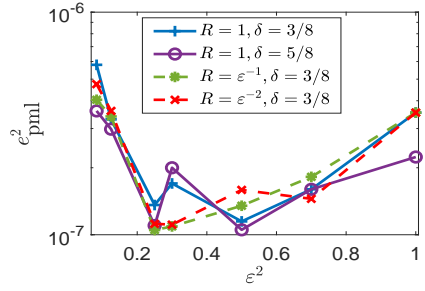,height=4cm,width=6cm}
\end{array}$$
\caption{PML-II $e^2_\textrm{pml}$-error  with respect to $\eps$ for fixed $\sigma_0=3$ and different $\delta$ or $R$ for  $\sigma=\sigma_{B_2}$.}
\label{fig:PML2 error Bermudez}
\end{figure}

\begin{figure}[hbt!]
$$\begin{array}{cc}
\psfig{figure=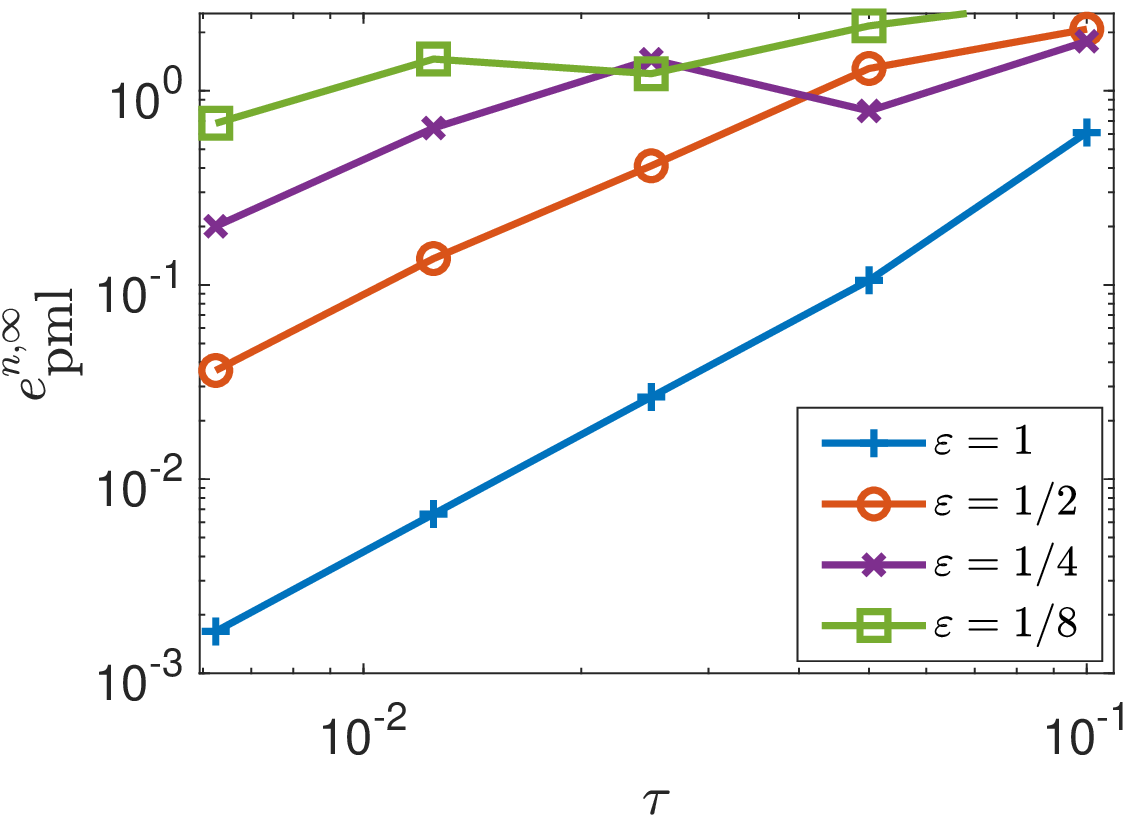,height=4cm,width=6cm}&
\psfig{figure=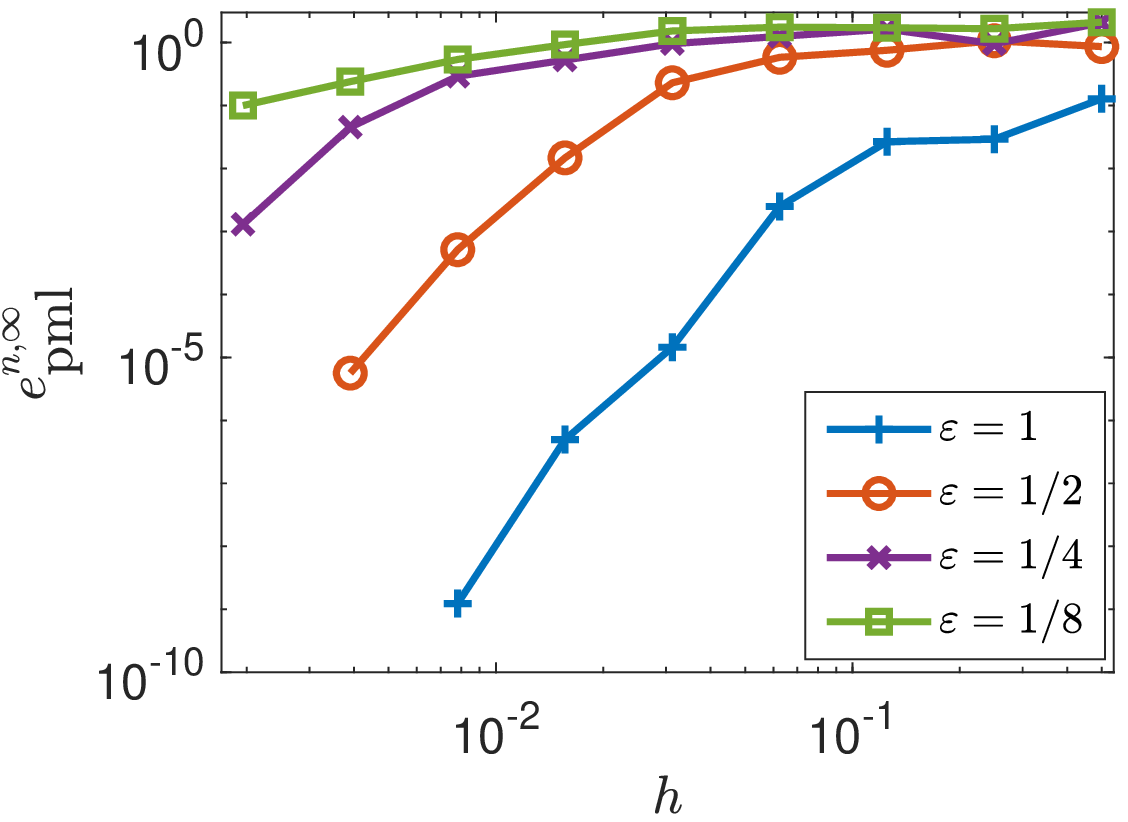,height=4cm,width=6cm}
\end{array}$$
\caption{Temporal (left) and spatial (right) accuracy of FD-FP:  $e^{n,\infty}_{\textrm{pml}}$-error  for $\sigma=\sigma_{B_2}$  with different $\eps$.}
\label{fig:pml2eps accuracy Bermudez}
\end{figure}

\begin{figure}[hbt!]
$$\begin{array}{cc}
\psfig{figure=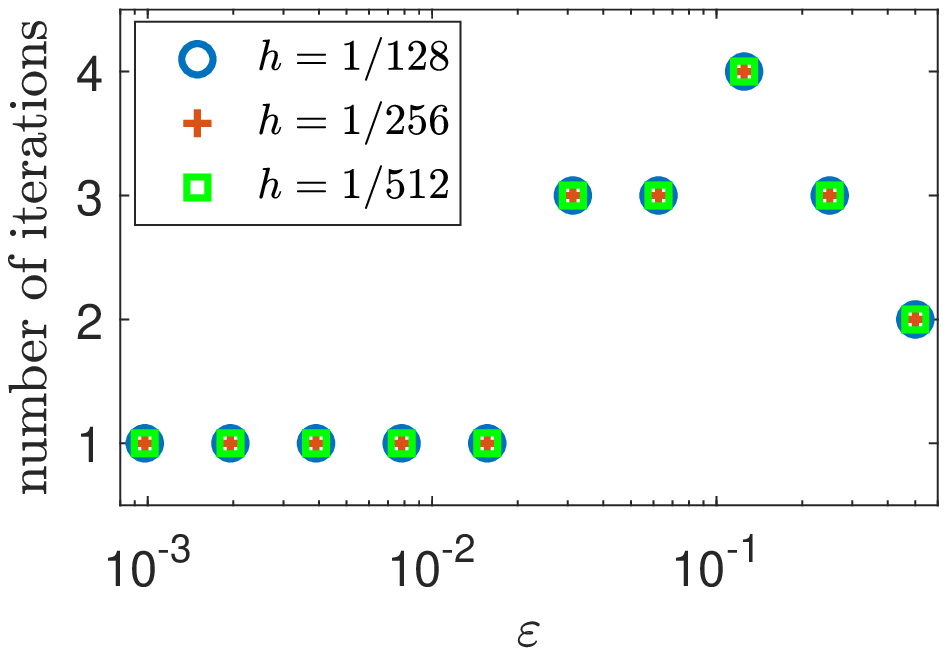,height=4cm,width=6cm}&
\psfig{figure=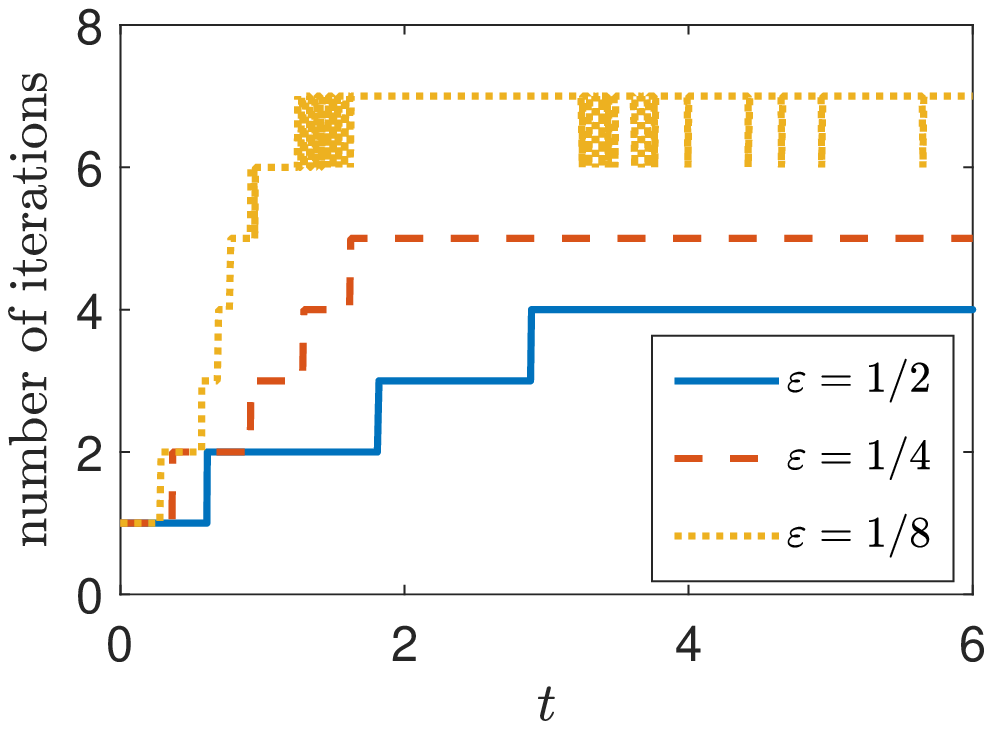,height=4cm,width=6cm}
\end{array}$$
\caption{Number of iterations needed by  GMRES (fixed threshold $\epsilon=10^{-10}$) with the preconditioner for (\ref{FD-FP eps}) at $n=1$ (left) and for (\ref{FD-FP eps}) until $t=6$ (right).}
\label{tab:gmres eps}
\end{figure}

\medskip
\noindent\textbf{PML-II with Berm\'udez type absorption function.} Let us consider
$\sigma(x)=\sigma_{B_k}(x)$ as defined in (\ref{Bermudez}) for the PML-II (\ref{KG limit pml2}).
To test the PML error, we
 take the PML-II (\ref{KG limit pml2}) with $\sigma=\sigma_{B_2}$, i.e. $k=2$ in (\ref{Bermudez}), for  $\sigma_0=3$ fixed. By using some different values of $R>0$ and $\delta>0$, we show in Figure \ref{fig:PML2 error Bermudez} the PML approximation error (\ref{err}) on the physical domain $I=(-4,4)$ at $t=4$.  With $\sigma=\sigma_{B_2},\,R=1,\,\delta=4/8$, we plot in Figure \ref{fig:pml2eps accuracy Bermudez} the numerical discretization error (\ref{errorinfty}) of the FD-FP method at $t=4$ in time and in space. Under $\tau=10^{-4},\,h=1/128$, the corresponding condition number of the matrix $G^\eps$ as a function of $\eps$ is also shown in Figure \ref{fig:cond}.  The number of iterations needed for the convergence of GMRES (for threshold $\epsilon=10^{-10}$ without restart)  with the preconditioner  (\ref{preconditioner eps}) for solving (\ref{FD-FP eps}) at $n=1$ under $\tau=0.02,\ \delta=3/8,\ R=1$ is given in Figure \ref{tab:gmres eps} (left) for different $\eps>0$ and values of $h$.  The number of iterations needed for solving (\ref{FD-FP eps}) under $\tau=2\times10^{-4},\,h=1/128,\,\delta=3/8,\ R=1$ until time $t=6$ is given in Figure \ref{tab:gmres eps} (right).
Let us remark that, without this preconditioner, the iterative solver may not converge when $\eps$ gets smaller.

 It can be seen from Figures \ref{fig:PML2 error Bermudez}-\ref{tab:gmres eps} that

\begin{itemize}
\item[1)] The numerical discretization error of the FD-FP method for PML-II (\ref{KG limit pml2})  with the Berm\'udez absorption function (\ref{Bermudez}) is close to the case of the polynomial choice above (cf. Figures \ref{fig:pml2eps accuracy} and  \ref{fig:pml2eps accuracy Bermudez}).  The temporal and spatial error grows dramatically as $\eps$ decreases, but with $\tau,\,h\to0$, the FD-FP method still converges and the condition number of $G^\eps$ is again uniformly bounded for $0<\eps\leq1$.   The preconditioner (\ref{preconditioner eps}) still works very well for $0<\eps<1$, which improves significantly the efficiency of the GMRES solver, since it exhibits a convergence rate independent of the mesh refinement $h$. Moreover, the number of iterations stays small and bounded in time,
and is relatively insensitive to $\varepsilon$.

 \item[2)] The PML-II (\ref{KG limit pml2}) with the Berm\'udez absorption function (\ref{Bermudez}) is much more accurate than that of the classical polynomial choice (\ref{sigma}) (cf. Figures \ref{fig:PML2 error Bermudez} and  \ref{fig:eps}).  With the Berm\'udez  function, the PML-II is still not sensitive to the parameters for the layer, i.e.  $R,\,\sigma_0,\,\delta$. More importantly, the error appears to be rather uniform with respect to $\eps\in (0,1]$ as far as we could test. Thus, we can simply take fixed $R,\,\sigma_0,\,\delta=O(1)$ for all $\eps \in(0,1]$.  While, we remark that for  smaller $\eps$, it is getting  harder to eliminate the impact from the numerical error of the FD-FP  solver due to error behaviour observed  in Figure \ref{fig:pml2eps accuracy Bermudez}.
\end{itemize}

\begin{figure}[hbt!]
$$\begin{array}{c}
\psfig{figure=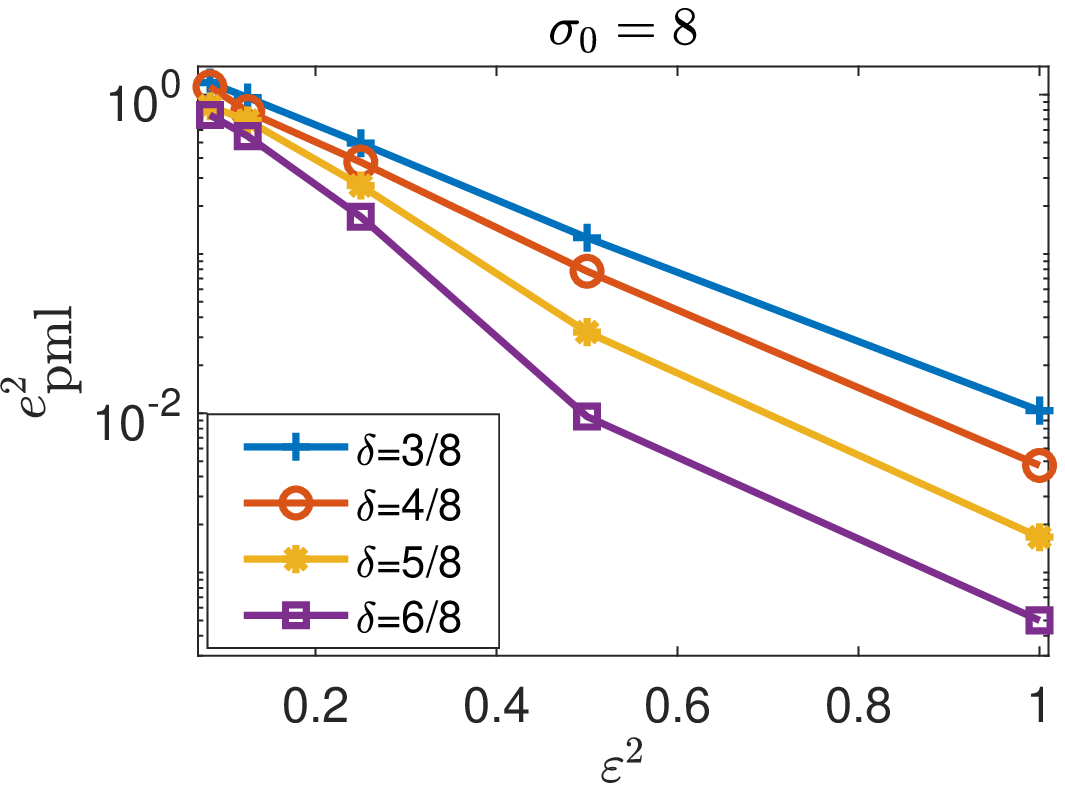,height=4cm,width=6cm}
\psfig{figure=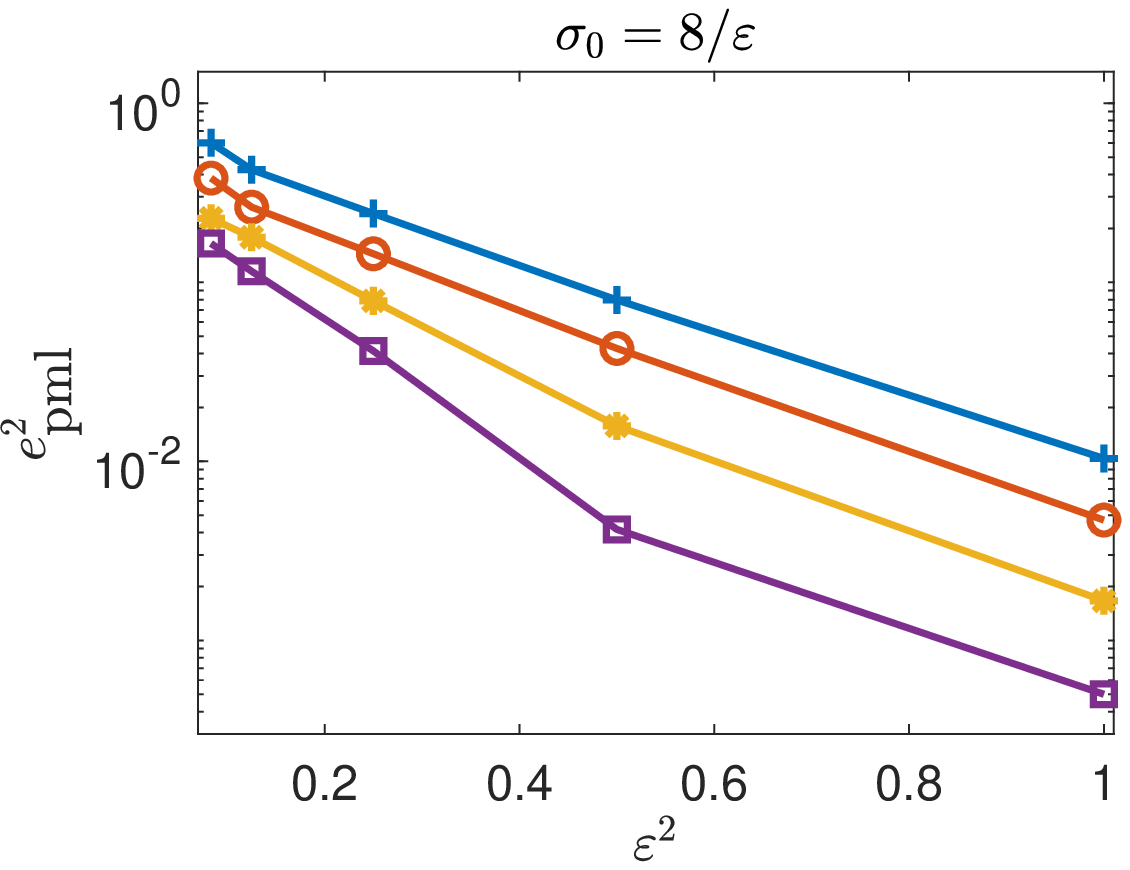,height=4cm,width=6cm}\\
\psfig{figure=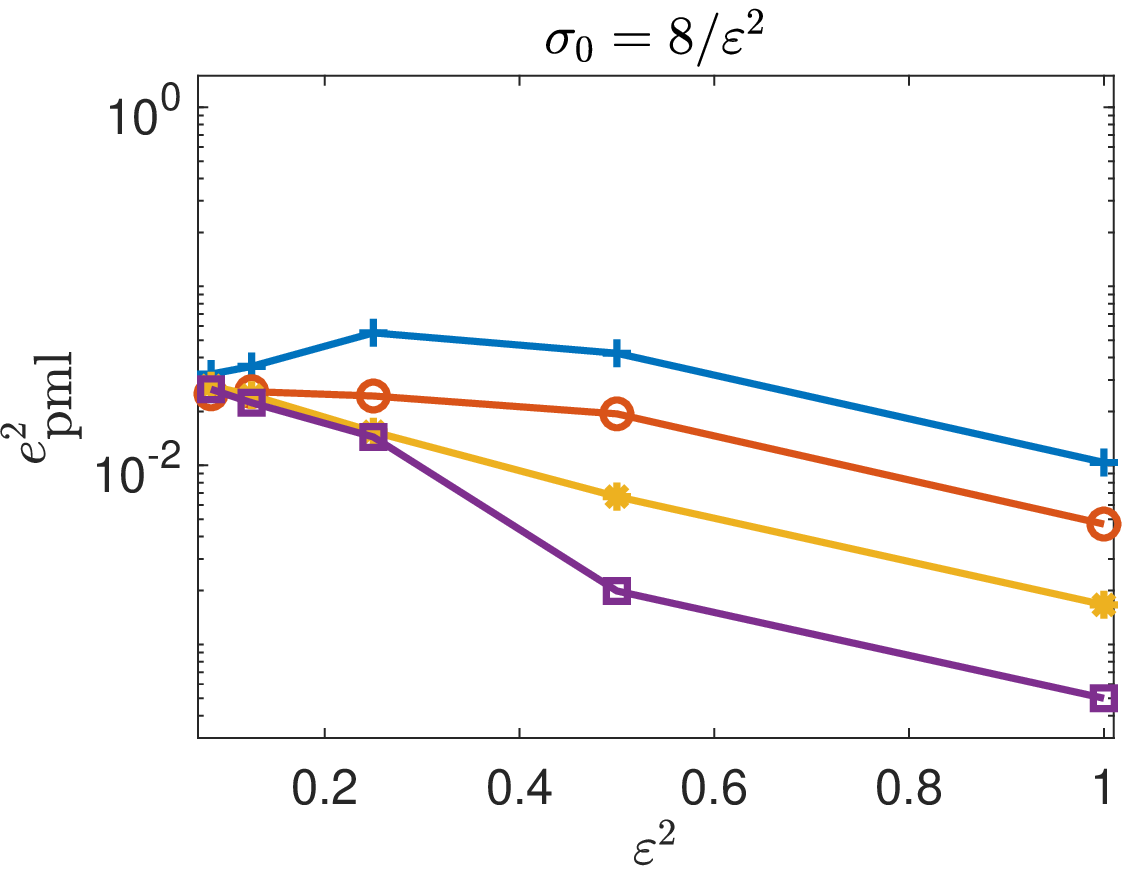,height=4cm,width=6cm}
\end{array}$$
\caption{PML-I $e^2_\textrm{pml}$-error  with respect to $\eps$ for different $\delta$:  $\sigma_0=8$ (top left), $\sigma_0=8/\eps$ (top right) and $\sigma_0=8/\eps^2$ (bottom).}
\label{fig:PML1 error eps}
\end{figure}

\noindent\textbf{Tests for PML-I.}
At last, we report some numerical results for the PML-I formulation (\ref{PML-I eps}) when approximating the NKGE (\ref{KG model eps 1d}) and illustrate its drawbacks
compared with PML-II (\ref{KG limit pml2}).  We use the same numerical example (\ref{example eps}) and  take $\alpha=0$ in (\ref{PML-I eps}). The PML error (\ref{err}) of the PML-I formulation (\ref{PML-I eps}) at $t=4$ with respect to $\eps$ is given in Figure \ref{fig:PML1 error eps}  under $\sigma_0=8$ and different parameter values $\delta$.
It is clear from the numerical results in Figure \ref{fig:PML1 error eps} that the PML-I formulation (\ref{PML-I eps}) with the classical  absorption function (\ref{sigma}) is very sensitive to the parameters: $\sigma_0,\,\delta$ and $\eps$. When $\eps$ becomes small, the error of PML-I is much larger than that of the PML-II (cf. Figure \ref{fig:PML1 error eps} and Figures \ref{fig:eps} and \ref{fig:PML2 error Bermudez}). In particular, for small values of $\eps$, the accuracy of PML-I can barely be improved by increasing $\delta$. One must keep enlarging the strength parameter $\sigma_0$, while the stiffness brought by large $\sigma$ in (\ref{PML-I eps}) will consequently increase the difficulty for the numerical solver.

\section{Extension to the 2D rotating NKGE}\label{Sec:rotating}
In this section, we shall apply the PML-II formulation to simulate the vortices dynamics in a two-dimensional rotating NKG model.
The rotating NKGE has been introduced to model the dynamics of the cosmic superfluid in a rotating frame, and in two space dimensions it reads \cite{RKG,RKG-PRD}:
\begin{align*}&\partial_{tt}\Psi-\Delta\Psi+\Psi
+\lambda|\Psi|^2\Psi-2i\Omega L_z\partial_t\Psi-\Omega^2L_z^2\Psi=0,\quad \bx=(x,y)^\intercal\in\bR^2,\ t>0,\\
&\Psi(\bx,0)=\psi_0(\bx),\quad \partial_t\Psi(\bx,0)=\psi_1(\bx),
\quad \bx\in\bR^2,
\end{align*}
where the solution $\Psi=\Psi(\bx,t)$  is a complex-valued scalar field, $\Omega>0$ is the given angular velocity and $L_z=-i(x\partial_y-y\partial_x)$ is the angular momentum operator.
As  proposed in \cite{RKG}, by introducing the rotating matrix
$$A(t):=\begin{pmatrix}\cos(\Omega t)&\sin(\Omega t)\\ -\sin(\Omega t)&\cos(\Omega t)\end{pmatrix},\quad
t\geq0,$$
and a change of variable to the \emph{rotating Lagrangian coordinates}:
\begin{equation*}u(\tilde{\bx},t)=\Psi(\bx,t),\quad \bx=A(t)\tilde{\bx},
\end{equation*}
the two angular momentum terms can be eliminated and the model becomes the NKGE (removing the $\tilde{\,}$ for conciseness) \cite{RKG}:
\begin{equation}\label{KG 2dmodel}\left\{\begin{split}
&\partial_{tt}u(\bx,t)-\Delta u(\bx,t)+u(\bx,t)+
\lambda|u(\bx,t)|^2u(\bx,t)=0,\quad \bx\in \bR^2,\ t>0,\\
&u(\bx,0)=\psi_0(\bx),\quad \partial_tu(\bx,0)=\Omega\nabla\psi_0(\bx)\cdot
\binom{y}{-x}+\psi_1(\bx),\quad \bx\in\bR^2.\end{split}\right.
\end{equation}
Then, we can directly extend the PML-II formulation to the above initial value problem which leads to
\begin{equation}\label{KG 2dmodel pml2}
\left\{\begin{split}
&\partial_{tt} u(\bx,t)-\frac{1}{1+R\sigma(x)}\partial_x\left(\frac{1}{1+R\sigma(x)}
\partial_xu(\bx,t)\right)
-\frac{1}{1+R\sigma(y)}\partial_y\left(\frac{1}{1+R\sigma(y)}\partial_yu(\bx,t)\right)\\ &\hspace{3cm}+u(\bx,t)+\lambda |u(\bx,t)|^2u(\bx,t)=0,\quad t>0,\ \bx\in I^*=(-L^*,L^*)^2,\\
&u(\bx,0)=\psi_0(\bx),\quad \partial_tu(\bx,0)=\Omega\nabla\psi_0(\bx)\cdot
\binom{y}{-x}+\psi_1(\bx), \quad \bx\in I^*,\\
&u(-L^*,y,t)=u(L^*,y,t),\quad t\geq0,\quad y\in[-L^*,L^*],\\
&u(x,-L^*,t)=u(x,L^*,t),\quad t\geq0,\quad x\in[-L^*,L^*],
\end{split}\right.
\end{equation}
where $L^*=L+\delta$.
With
$$S_x=-\frac{1}{1+R\sigma(x)}\partial_x\left(\frac{1}{1+R\sigma(x)}\partial_x\right),
\quad S_y=-\frac{1}{1+R\sigma(y)}\partial_y\left(\frac{1}{1+R\sigma(y)}\partial_y\right),$$ for short,
the finite-difference time discretization for (\ref{KG 2dmodel pml2}) then reads
\begin{align*}
  &\frac{u^{n+1}-2u^n+u^{n-1}}{\tau^2}+\frac{S_x+S_y}{2}
  (u^{n+1}+u^{n-1})+\frac{1}{2}
  (u^{n+1}+u^{n-1})+\lambda |u^n|^2u^n=0,\quad n\geq1,\\
  &u^1=u_0+\tau v_0-\frac{\tau^2}{2}\left[(S_x+S_y) u_0+u_0+\lambda |u_0|^2u_0\right],
\end{align*}
where $u^n=u^n(\bx)\approx u(\bx,t)$ for $n\geq1$ and $u_0=u(\bx,0),\ v_0=\partial_tu(\bx,0)$.
By applying the Fourier pseudo-spectral discretization in the $x$- and $y$-directions to the equations, we obtain the FD-FP scheme. By further incorporating the previous equations with the preconditioned GMRES and FFT in 2D, we can implement the FD-FP scheme in the same manner as before:
\begin{equation*}
\left\{
  \begin{array}{l}
  \displaystyle
  u^{n+1}=-u^{n-1}
  +w^n,\quad n\geq1,\\
   \displaystyle  \mathcal{P}G w^n   =\mathcal{P}\left[\frac{2}{\tau^2}u^n -\lambda |u^n|^2u^n\right],
  \end{array}
  \right.\quad
  \end{equation*}
where
$$G=\left[\left(\frac{1}{\tau^{2}}+\frac{1}{2}\right)\mathbb{I}
  +\frac{S_x+S_y}{2}\right],\quad\mathcal{P}=\left(\frac{1}{\tau^{2}}
  +\frac{1}{2}-\frac{\partial_{xx}}{2}-\frac{\partial_{yy}}{2}\right)^{-1}.$$

\begin{figure}[hbt!]
$$\begin{array}{cc}
\psfig{figure=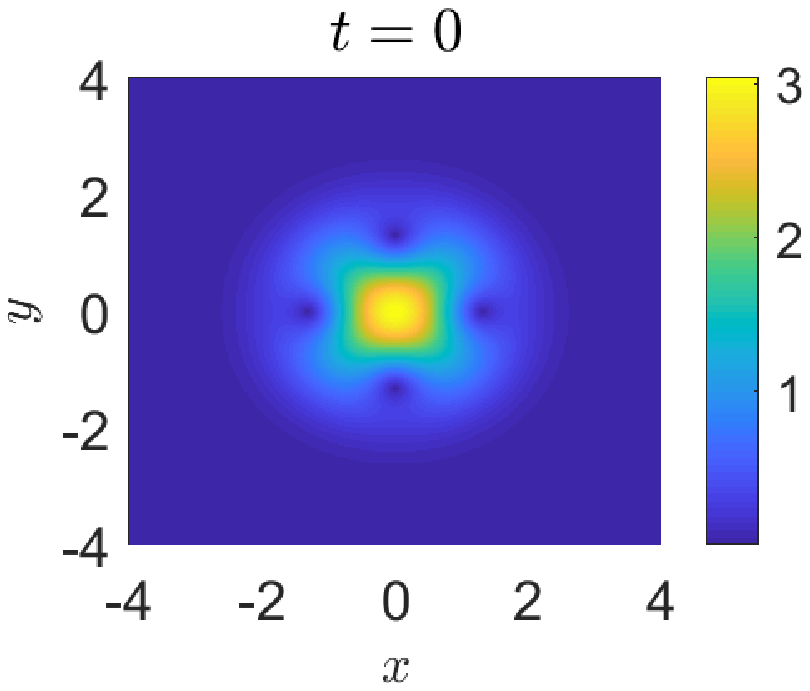,height=4.5cm,width=4.5cm}&
\psfig{figure=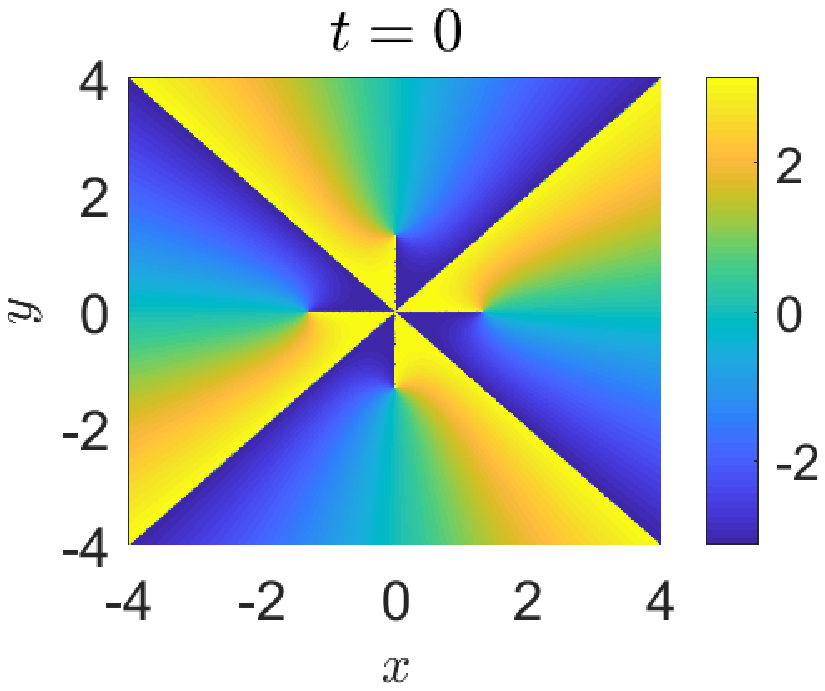,height=4.5cm,width=4.5cm}
\end{array}$$
\caption{Contour plot of the initial data: $|\psi_0(\bx)|$ (left) and $\mathrm{Arg}(\psi_0(\bx))$ (right)  for (\ref{KG 2dmodel}) or (\ref{KG 2dmodel pml2}).}
\label{fig:2dt0}
\end{figure}

To illustrate the accuracy of the PML-II approach, we consider the following example. We take in (\ref{KG 2dmodel pml2}) the parameters
$$L=4,\quad \lambda=3,\quad \Omega=2,\quad R=1,$$
and four initially separated vortices
$$\psi_0(\bx)=\psi_1(\bx)=(x-c_0+iy)(x+c_0+iy)(x+i(y-c_0))(x+i(y+c_0))
\fe^{-(x^2+y^2)/2},$$
with $c_0=1.32$. The profile of the initial data for the NKGE is shown in Figure \ref{fig:2dt0}.

\begin{figure}[hbt!]
\centerline{Exact solution: $|u(\bx,t)|$ (top) and $\mathrm{Arg}(u(\bx,t))$ (bottom)}
$$\begin{array}{ccc}
\psfig{figure=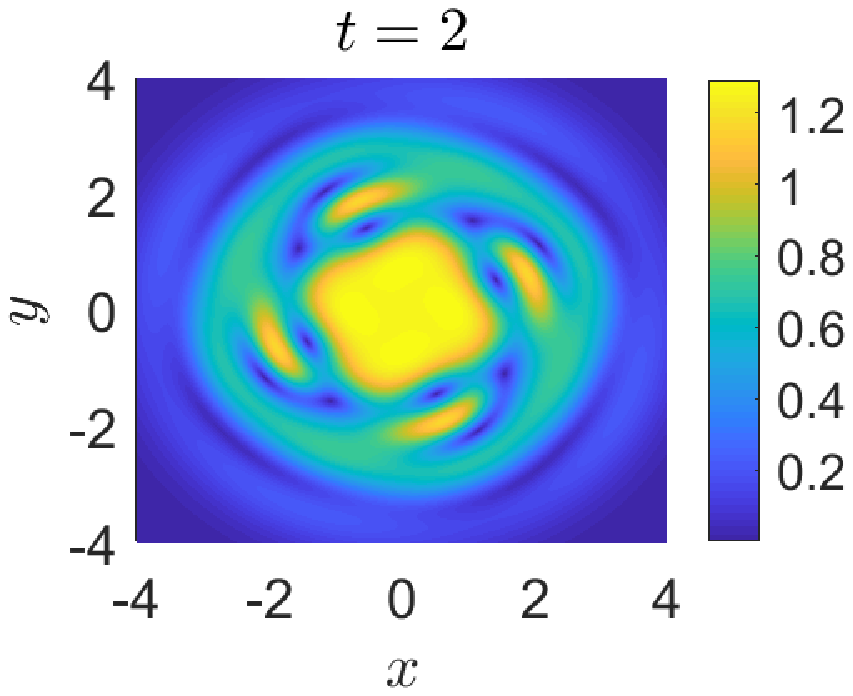,height=4.5cm,width=4.5cm}&
\psfig{figure=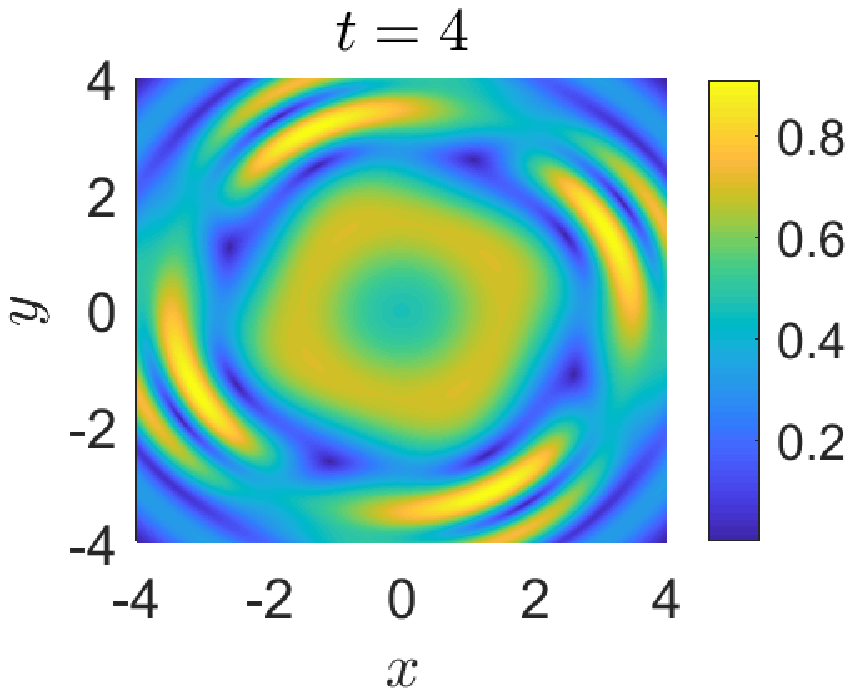,height=4.5cm,width=4.5cm}&
\psfig{figure=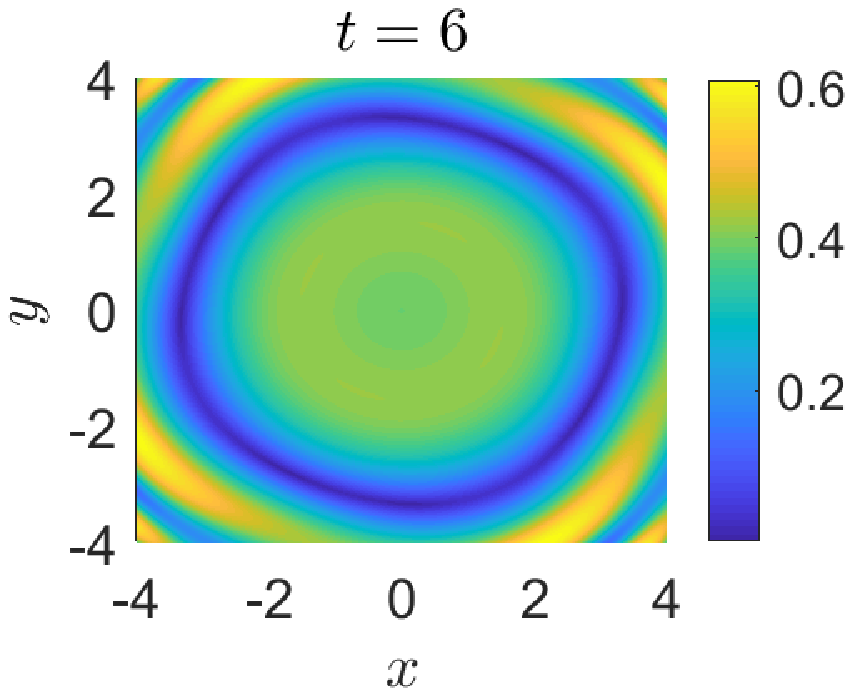,height=4.5cm,width=4.5cm}\\
\psfig{figure=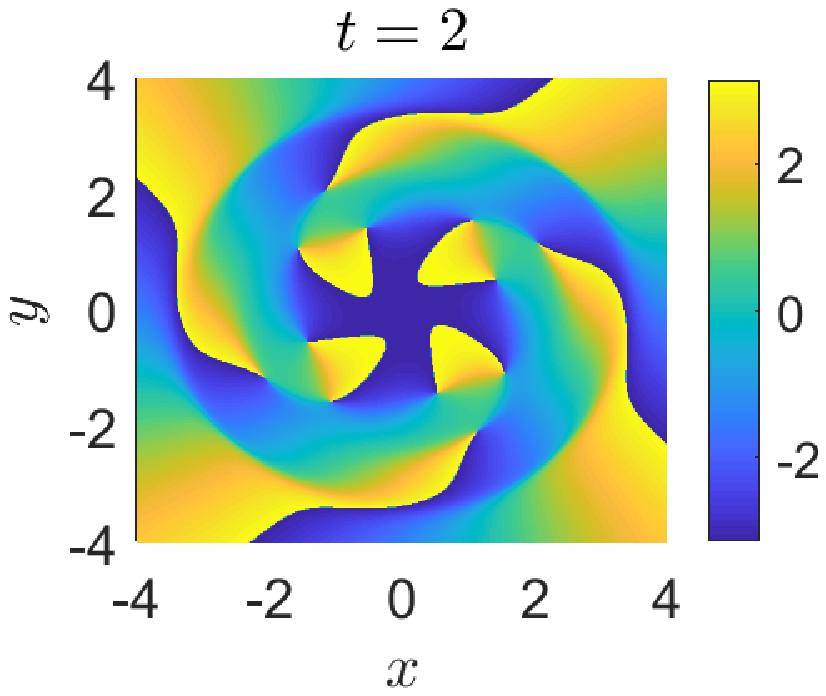,height=4.5cm,width=4.5cm}&
\psfig{figure=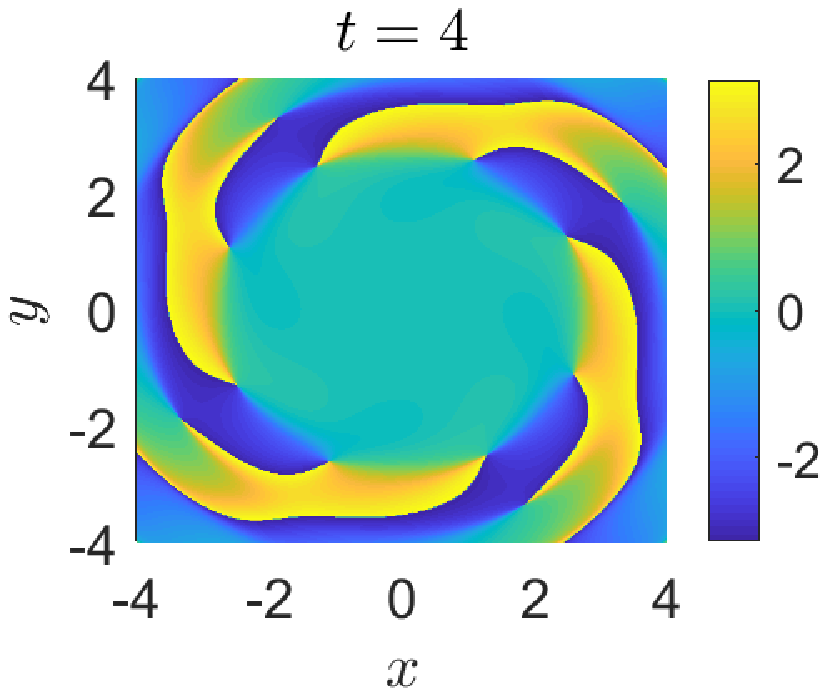,height=4.5cm,width=4.5cm}&
\psfig{figure=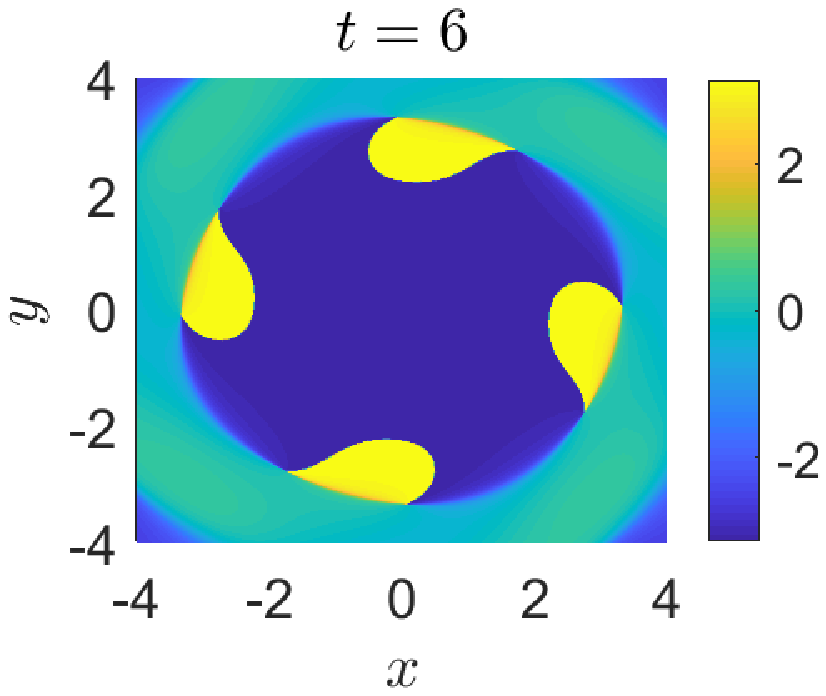,height=4.5cm,width=4.5cm}
\end{array}$$
\centerline{PML solution: $|u_{\textrm{pml}}(\bx,t)|$ (top) and $\mathrm{Arg}(u_{\textrm{pml}}(\bx,t))$ (bottom)}
$$\begin{array}{ccc}
\psfig{figure=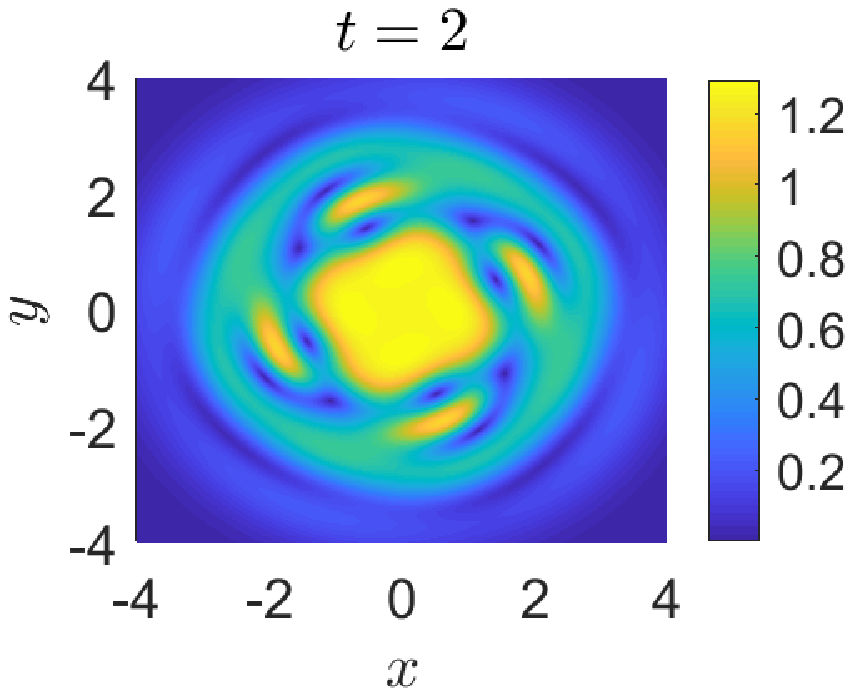,height=4.5cm,width=4.5cm}&
\psfig{figure=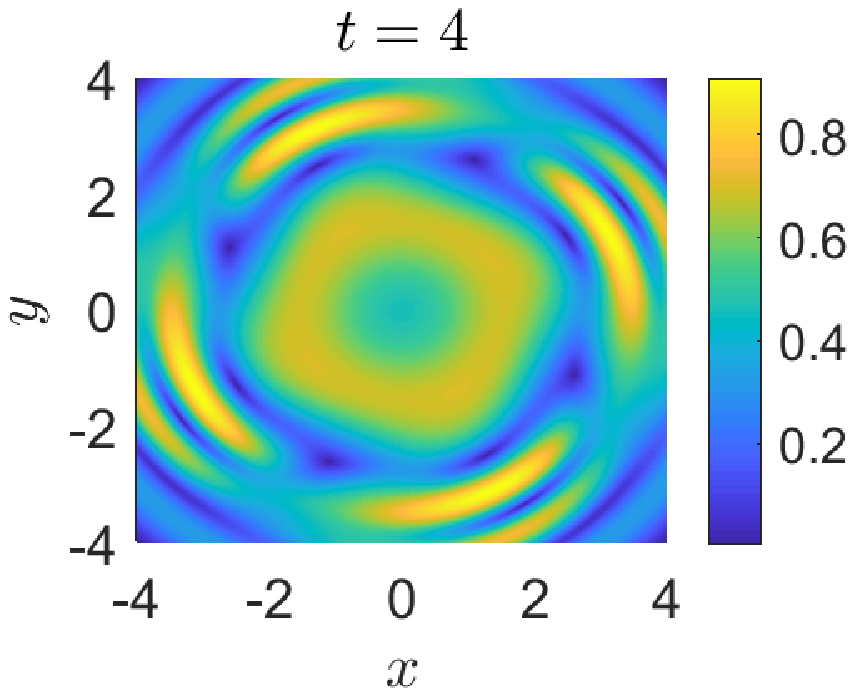,height=4.5cm,width=4.5cm}&
\psfig{figure=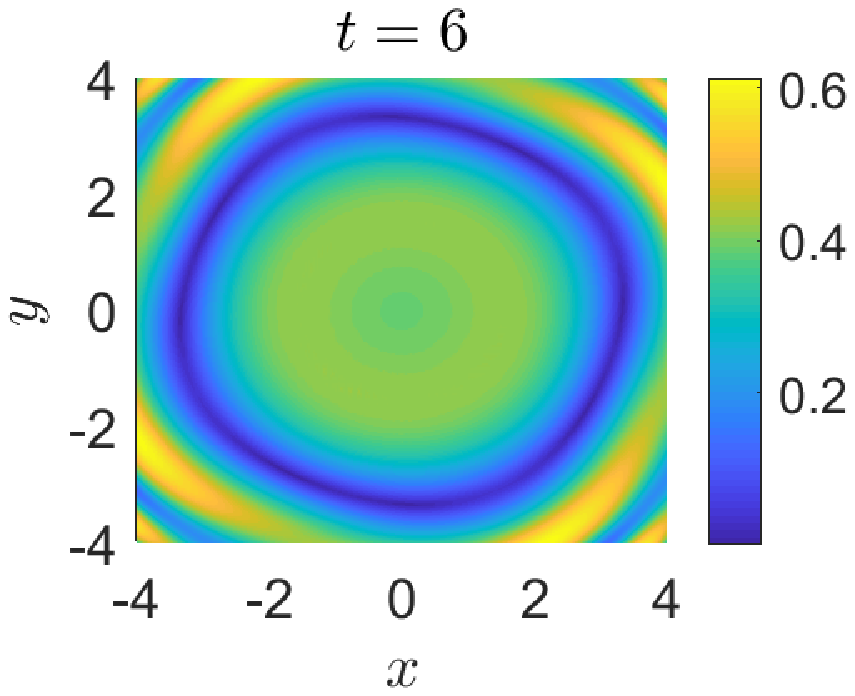,height=4.5cm,width=4.5cm}\\
\psfig{figure=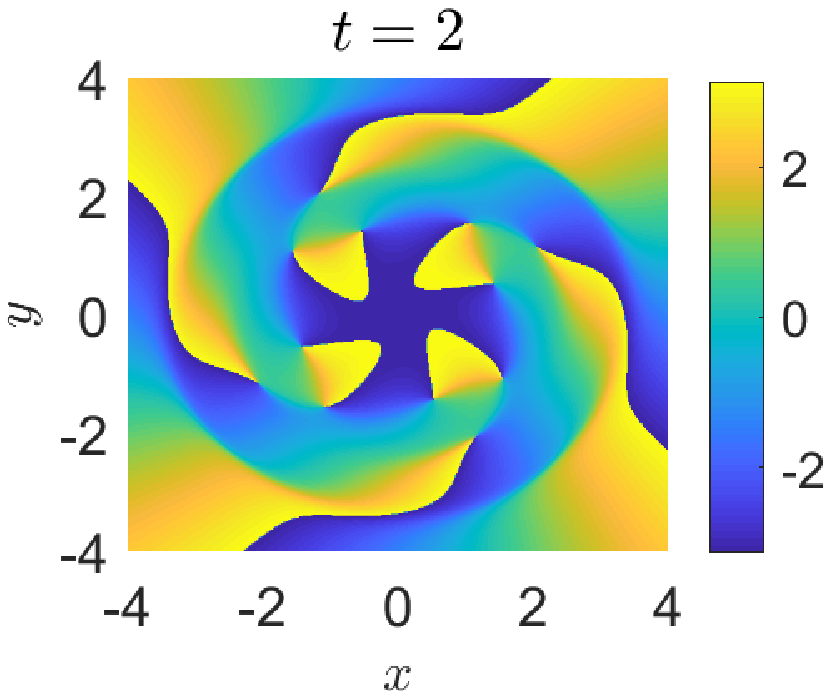,height=4.5cm,width=4.5cm}&
\psfig{figure=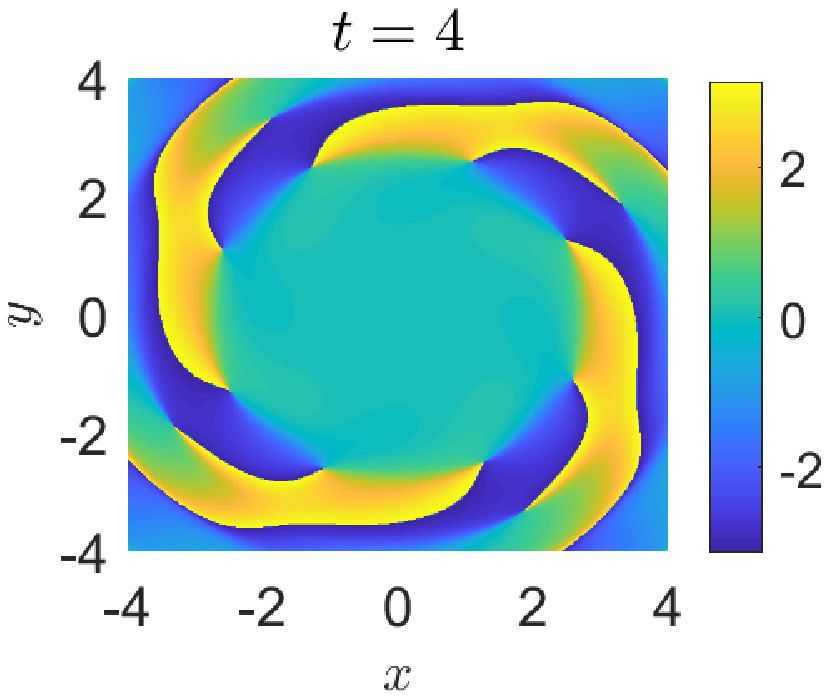,height=4.5cm,width=4.5cm}&
\psfig{figure=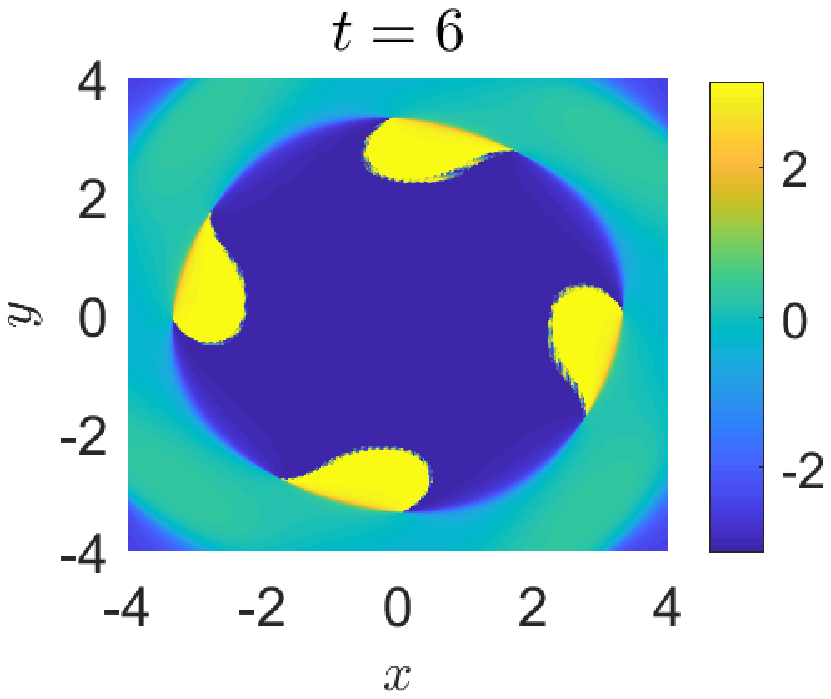,height=4.5cm,width=4.5cm}
\end{array}$$
\caption{Contour plots of exact solution $|u(\bx,t)|$ and $\mathrm{Arg}(u(\bx,t))$ of (\ref{KG 2dmodel}) and the PML solution $|u_{\textrm{pml}}(\bx,t)|$ and $\mathrm{Arg}(u_{\textrm{pml}}(\bx,t))$ of (\ref{KG 2dmodel pml2}) at  times $t=2, 4, 6$ on the domain $I$.}
\label{fig:2d}
\end{figure}

\begin{figure}[hbt!]
$$\begin{array}{c}
\psfig{figure=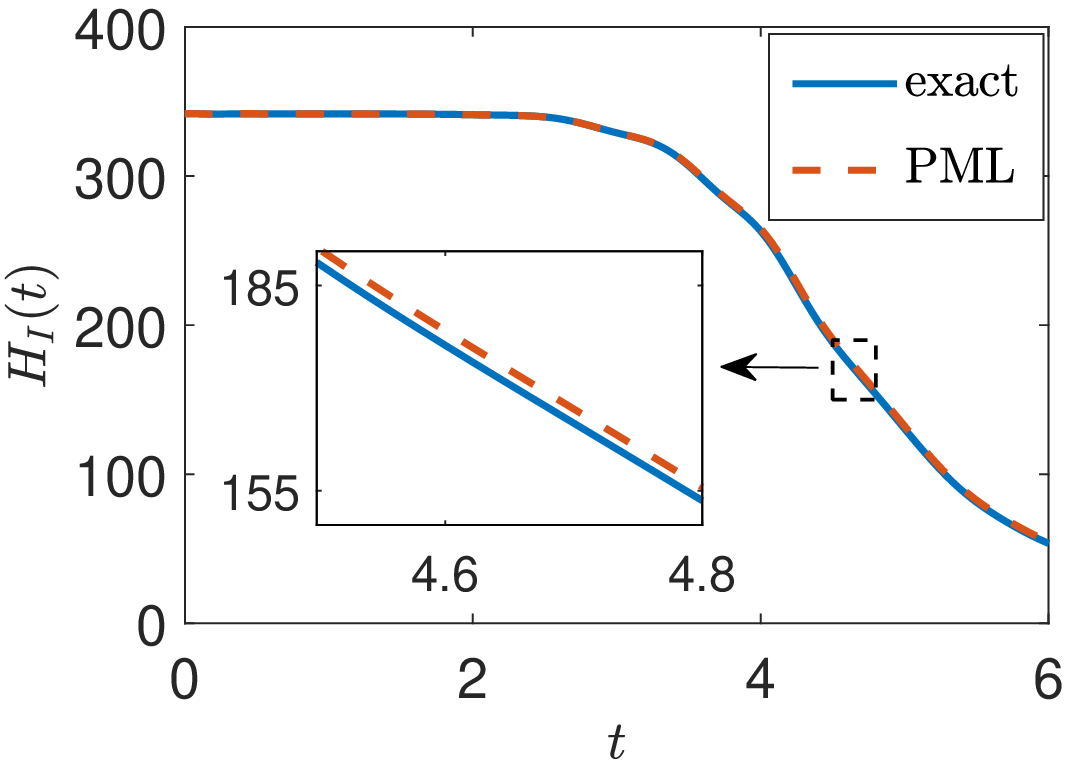,height=4cm,width=6cm}
\end{array}$$
\caption{The energy decay of the 2D example: $H_I(t)$ with respect to time.}
\label{fig:2d energy}
\end{figure}

We use the Berm\'udez function $\sigma=\sigma_{B_2}$ (\ref{Bermudez}) as the absorption function and choose $\sigma_0=3$ and $\delta=0.5$.
By applying the FD-FP similarly as in Section \ref{sec2 PML2}, we solve the PML equation (\ref{KG 2dmodel pml2}) accurately and look at the solution $u_{\textrm{pml}}$ on the physical domain $I=(-L,L)^2$. As a reference solution, we solve (\ref{KG 2dmodel}) directly on the sufficiently large domain $(-16,16)^2$ and compare the exact solution $u$ with the PML solution. The dynamics of the exact solution $u(\bx,t)$ and the PML solution $u_{\textrm{pml}}(\bx,t)$ are shown on the domain $I$ in Figure \ref{fig:2d}. The corresponding trajectory of the decaying energy $H_I(t; w)$:
$$ H_I(t;w):=\int_{I}\left[|\partial_tw(x,t)|^2+|\nabla w(x,t)|^2
    +|w(x,t)|^2+\frac{\lambda}{2}|w(x,t)|^4\right]dx,\quad t\geq0,$$
with $w=u$ or $u_{\textrm{pml}}$ is reported in Figure \ref{fig:2d energy}. For both the solution and the energy,
 we see an excellent agreement between both solutions that cannot be reached when using standard periodic boundary conditions, unless taking a huge spatial domain.

\section{Conclusion}\label{sec:conclusion}
In this work, we investigated numerically two kinds of PML formulations for approximating the nonlinear Klein-Gordon equation: a first-order formulation borrowed from the study for general wave equations (PML-I), and a second-order formulation proposed as analogy for the Helmholtz and nonlinear Schr\"odinger equations (PML-II). For both PML formulations, we developed  efficient and accurate pseudo-spectral schemes for the numerical discretizations: i) a pseudo-spectral exponential integrator scheme for PML-I and ii) a linearly implicit preconditioned pseudo-spectral
finite-difference scheme for PML-II. To obtain a spectral accuracy,  
the PML absorption functions were locally regularized.  For the NKGE in classical scaling, the PML-II with smoothed Berm\'udez type absorption functions and the resulting pseudo-spectral finite-difference scheme are proved to be the most accurate and efficient,  owning to the observation that it is less subject to the parameter tuning problems for the layer. Such advantage also distinguishes PML-II for approximating NKGE in the non-relativistic scaling, where the PML solution
is found to be not sensitive to the small parameter $\varepsilon$ arising 
in the non-relativistic limit regime.  
We also extend the method to the accurate computation of vortex dynamics for a rotating NKGE.

\appendix
\section*{Acknowledgements}
Research conducted within the context of the Sino-French International Associated Laboratory for Applied Mathematics - LIASFMA. X. Antoine thanks the LIASFMA funding support of the Universit\'e de Lorraine. X. Zhao is partially supported by the Natural Science Foundation of Hubei Province No. 2019CFA007 and the NSFC 11901440.

\bibliographystyle{model1-num-names}

\end{document}